\documentclass[authoryear,preprint,review,11pt]{elsarticle}
\biboptions{authoryear}

\usepackage{hyperref}
\hypersetup{
  colorlinks   = true, %Colours links instead of ugly boxes
  urlcolor     = blue, %Colour for external hyperlinks
  linkcolor    = blue, %Colour of internal links
  citecolor   = blue!70 %Colour of citations
}
\hypersetup{
  pdfauthor={Vítor A. Barbosa, Sunil Tiwari, Rafael A. Melo},
  pdftitle={The pickup and delivery problem with time windows and scheduling on the edges}
}

\usepackage[numbered]{bookmark}

\usepackage{enumitem}
\setlist[enumerate,1]{label=(\alph*), ref=(\alph*)}
\setlist[enumerate,2]{label=\roman*., ref=\roman*.}

\usepackage{float} %To use 'H' in the Table and Figure float option and fix the position where it is.
\floatstyle{plaintop} %To re-position the Table caption on top of the table
\restylefloat{table}

\usepackage{xfrac}
\usepackage{color, colortbl}
\usepackage[dvipsnames, table]{xcolor}
\usepackage{cancel}
\usepackage{bbding}
\usepackage{rotating}
\usepackage{array}[=2016-10-06] % arxiv suggestion for submissions https://info.arxiv.org/help/submit_tex.html#issues-with-the-array-package
\usepackage{animate}
\usepackage{tikz}
\usepackage{pgfplots}
\usepackage{mathtools}
\usepackage{pifont}% http://ctan.org/pkg/pifont
\usetikzlibrary{fit, shapes.misc}
\usepackage{lscape}

\usepackage{amsmath,amssymb,amsthm,amsfonts,dsfont,mathrsfs,latexsym,rotate,epic,eepic,epsfig,graphicx}

\usepackage{comment}
\usepackage{setspace}
\usepackage{geometry}
\DeclarePairedDelimiter\ceil{\lceil}{\rceil}
\DeclarePairedDelimiter\floor{\lfloor}{\rfloor}
\renewcommand\appendix{\par
  \setcounter{section}{0}
  \setcounter{subsection}{0}
  \setcounter{figure}{0}
  \setcounter{table}{0}
  \renewcommand\thesection{Appendix \Alph{section}}
  \renewcommand\thefigure{\Alph{section}\arabic{figure}}
  \renewcommand\thetable{\Alph{section}\arabic{table}}
  \renewcommand\thesubfigure{(\alph{subfigure})} % used only for figure label index
}
\theoremstyle{plain}
\theoremstyle{plain}
\theoremstyle{remark}
\theoremstyle{definition}

\pgfplotsset{compat=1.14}

\usetikzlibrary{patterns}
\newcommand{\okC}[1]{\textcolor{black}{#1}}

\usepackage{algorithmic}
\usepackage[ruled,commentsnumbered,linesnumbered,vlined]{algorithm2e}
\usepackage{varwidth}% http://ctan.org/pkg/varwidth

\usepackage{pdflscape}

\usepackage{longtable}
\usepackage{picture}
\usepackage{placeins} % \FloatBarrier
\usepackage{adjustbox}
\usepackage{booktabs, multirow} % for tables
\usepackage{changepage,threeparttable} % for wide tables

\usepackage{lmodern}
\usepackage[skip=2pt,font=scriptsize]{caption}

\usepackage{subfigure}

\usepackage{soul} % \st command to strikeout text
\usepackage[normalem]{ulem}
\makeatletter
\newcommand{\captshort}{\let\blx@imc@ifciteseen\@firstoftwo}

\usepackage{siunitx}

\begin{document}
\begin{frontmatter}
\title{The pickup and delivery problem with time windows and scheduling on the edges}

\author[1st_address]{Vítor A. Barbosa}
\ead{vitor.barbosa@ufba.br}

\author[2nd_address]{Sunil Tiwari}
\ead{sunil.tiwari@bristol.ac.uk}

\author[1st_address]{Rafael A. Melo\corref{cor1}}\cortext[cor1]{Corresponding author}
\ead{rafael.melo@ufba.br}

\address[1st_address]{Institute of Computing, Universidade Federal da Bahia, Salvador, BA 40170-115, Brazil}

\address[2nd_address]{Technology \& Operations Group, University of Bristol Business School, Bristol, United Kingdom}

% \address[3rd_address]{Institute of Computing, Universidade Federal da Bahia, Salvador, BA 40170-115, Brazil}

\begin{abstract}
We introduce the Pickup and Delivery Problem with Time Windows and Scheduling on the Edges (PDPTW-SE), a generalization of the PDPTW that integrates vehicle routing and machine scheduling. The problem involves defining routes for transportation requests with specific pickup and delivery locations using a heterogeneous vehicle fleet, while machines must be scheduled to traverse certain edges. The objective is to minimize the total completion time subject to capacity, time window, and precedence constraints. We propose a mixed-integer linear programming (MIP) formulation, including preprocessing and valid inequalities, and a multi-start heuristic with a linear programming (LP) improvement procedure. A benchmark set with two instance families is also introduced: (i) coordination of pickups and deliveries across islands requiring cargo ships, and (ii) transport across multiple floors, as in hospitals, requiring elevator scheduling. Computational experiments show that the solver on the MIP formulation solves instances with up to 12 requests and finds feasible solutions for 95.0\% of the 320 instances with up to 12 requests. 
For these, the heuristic consistently provides feasible solutions with low deviations, often matching or outperforming the MIP results. For the remaining 160 instances with 40 and 60 requests, only the heuristic finds feasible solutions. We thus recommend the MIP for short-horizon instances (up to 12 requests) and the heuristic for larger or long-horizon instances. Results also highlight the LP improvement procedure’s relevance, reducing solution values by at least 5\% on average in general. For larger, tightly constrained instances, an additional machine helps with feasibility and solution quality.
\end{abstract}


\end{frontmatter}

%%%%%%%%%%%%%%%%%%%%%%%%%%%%%%%%%%%%%%%%%%%%%%%%%%%%%%%%%%%%

% %%Graphical abstract
% \begin{graphicalabstract}
% % \includegraphics{grabs}
% \end{graphicalabstract}

% %%Research highlights
% \begin{highlights}
% \item Research highlight 1
% \item Research highlight 2
% \end{highlights}

\section{Introduction}
\label{sec:introduction}
The pickup and delivery problem with time windows (PDPTW) is a widely studied optimization problem in the field of vehicle routing. In the PDPTW, a homogeneous fleet of vehicles must depart from the depot, fulfill pickup and delivery requests at various nodes, and return to the depot. A solution to the problem must satisfy several constraints, including precedence between pickups and deliveries, time windows, and the maximum capacity of the vehicles, which should never be exceeded during transportation. Numerous extensions of the PDPTW have been proposed in the literature (see Section \ref{sec:relatedworks}). In this work, we present a generalization of the PDPTW, denoted as \textit{pickup and delivery problem with time windows and scheduling on the edges} (PDPTW-SE). The PDPTW-SE extends the PDPTW by allowing pickup and delivery nodes to be located in different regions inaccessible to regular vehicles. To fulfill the requests, a specific type of machine must be employed to transport a vehicle between regions. However, due to limitations on the availability of these machines, they must be scheduled. The PDPTW-SE is a computationally challenging problem that integrates routing and scheduling decisions as it generalizes the PDPTW, which is known to be NP-hard \citep{LenKan81, DumDesSou91}.

\subsection{Potential applications}

Applications of the PDPTW-SE arise in a variety of practical situations. In what follows, we provide three examples of such potential applications.

\textbf{Application A:} 
One application of the PDPTW-SE is in the logistical planning of pickups and deliveries involving multiple islands. 
In this scenario, requests are placed on islands, and cargo ships are required to transport the vehicles between the islands. 
\cite{Calderwood2011} highlighted the main challenges involved in servicing the needs of an island community, such as the supply of short-shelf life products, not only on the main island but also on the surrounding island communities. 
They found that the distance to market and limited transport infrastructure increased product handling, leading to longer lead times and reducing the shelf-life of fresh items. 
Given that most goods are delivered by sea, ferry crossings are subject to prevailing weather conditions, which can further affect quality and availability. 
Therefore, trading on an island incurs higher operational costs. 
As stated in \cite{FREATHY2020180}, deficiencies in the transport infrastructure, spatial disparities in grocery provision, and mainland retailers' strategies all represent market impediments. 
Therefore, considering new operational approaches seems necessary to deal with such disadvantages and improve the economies of these fragile areas.

\textbf{Application B:} Another application of the PDPTW-SE is in the logistical planning of pickups and deliveries by automated guided vehicles (AGVs) in multi-floor buildings, such as hospitals or hotels. In these scenarios, requests are placed on different floors, and it may be necessary to use elevators to transport the vehicles between the floors. 
In the context of healthcare, for instance, AGVs can replace labor-intensive, repetitive,  low-value-added, and mundane rule-based tasks.
These tasks include the delivery and replacement of linen, food services to the wards during mealtimes, and the delivery of consumables such as adult diapers \citep{Ashrafian2015, Kim2021, Riek2017}.
Besides that, traditional daily medicine delivery involves a human delivery team supplying patients by handcarts. 
Due to the multiple steps in this medicine delivery process, efficiency is diminished, and the risk of contamination is heightened \citep{Chen2021b}. 
The COVID-19 pandemic further highlighted this scenario, stressing the importance of contactless delivery services \citep{Chen2021a}.
Therefore, hospitals would be better equipped to fight infections and the propagation of infectious diseases by adopting AGVs in their daily routine. Some prototypes for hospitals have been proposed in the literature \citep{BacDurBirKysPerPad17, Antony2020}.
\cite{Jones2021} references some examples of AGV usage in libraries, hotels, hospitals, and aged care centers. For instance, the Royal Adelaide Hospital in Australia employed a fleet of 25 AGVs to collect waste and deliver food, mail, laundry, and medical supplies throughout the hospital.

\textbf{Application C:} A third possible application of the PDPTW-SE is in the logistical planning of pickup and delivery of goods in areas affected by man-made or natural disasters.
Emergency logistics operations following large-scale disasters must consider additional constraints, as infrastructure can be severely damaged \citep{Jiang2012}. 
In the context of damaged transportation infrastructure, one can have blocked roads, trapped people, and delayed medical aid \citep{Shiri2020}.
According to \cite{Jiang2012}, it could be beneficial to explore new alternative solutions, such as investigating the combinatorial choice of multimode transportation routing to cope with infrastructure damage/availability.
In a scenario where all roads between regions are closed for regular vehicles, a possible alternative solution would be to use machines that can traverse between regions carrying regular vehicles. This is the case considered in the PDPTW-SE proposed in this work. 
\cite{FIKAR2016104} studied an alternative approach encompassing the disaster relief coordination between private and relief organizations to distribute goods by scheduling and routing trucks, off-road vehicles, and unmanned aerial vehicles. 
Off-road vehicles can traverse closed roads due to special permissions or equipment. 
To deliver goods in their studied problem, either a detour on intact roads is driven, or the shipment is coordinated with relief organizations at one of the given transfer points. 
Due to organizational factors, only a limited number of points can be operated.
The goods are delivered from the supply point to a transfer point by a private organization. 
The relief organization ships the goods to a second receiving transfer point.
From there, the final delivery to the demand point is performed by the private organization, which is assumed to have a private vehicle available for this final shipment. 
The main differences between this problem and the PDPTW-SE are that in the PDPTW-SE it is not possible to take a detour on intact roads, and the delivery of goods does not consider transshipment.
That is, the same vehicle is responsible for picking up and delivering the request directly.

\subsection{Main contributions and organization}

We now outline the main contributions of our work.
First, we introduce the PDPTW-SE, which novelly integrates routing and scheduling decisions, finding potential applications in various domains.
Second, we model the problem as a mixed-integer programming (MIP) formulation\okC{, including a preprocessing step and several valid inequalities}.
Third, we propose a multi-start heuristic with a linear programming (LP) improvement procedure.
Finally, we propose a benchmark set of instances to compare new approaches to the problem. This benchmark set characterizes two different potential applications of the problem. In the first one, the request nodes are spread across multiple islands. In the second one, these nodes are distributed among several floors.

The remainder of this paper is organized as follows. 
\okC{Section~\ref{sec:formalization} formally defines the PDPTW-SE.
Section~\ref{sec:relatedworks} revises the related literature.}
Section~\ref{sec:formulation} describes the proposed compact MIP formulation\okC{, preprocessing step, and valid inequalities}.
Section~\ref{sec:msheurlp} describes our multi-start heuristic with an LP improvement procedure.
Section~\ref{sec:experiments} discusses the performed computational experiments.
Section~\ref{sec:finalremarks} draws some concluding remarks.

\section{Problem statement}
\label{sec:formalization}

The PDPTW-SE can be formally defined using an undirected graph $G=(V,E)$ as follows. 
Let $V=V_p\cup V_d \cup \{0\}$, where $V_p = \{1,\ldots,n\}$ is a set of pickup nodes, $V_d = \{n+1,\ldots,2n\}$ is a set of delivery nodes, $n$ is the number of pickup/delivery nodes, and node $0$ represents the depot. 
In addition, let $E = E^s \cup E^m$ be the set of edges, where $E^s$ is the set of edges that can be directly traversed by a vehicle and $E^m$ is the set of edges whose traversals involve using a machine that has to be scheduled.  
Let $R = \{r_1,r_2,\ldots,r_n\}$ be a set of unsplittable requests, each denoted by a tuple $r_i~=~\langle v_i, v_{n+i}, q_i \rangle$, where $v_i \in V_p$ denotes its origin (pickup or loading node), $v_{n+i} \in V_d$ its destination (delivery or unloading node), and $q_i$ its weight. It is assumed that \okC{$q_i \in \mathbb{Z}_{\geq 0}$} for $i \in V_p$ and \okC{$q_i \in \mathbb{Z}_{\leq 0}$} for $i \in V_d$, with $q_i = -q_{n+i}$ for $i \in V_p$.
Each node $i \in V_p\cup V_d$ has a strict time window $[e_i,l_i]$ in which the vehicles should start the corresponding operation.
There is also a time window $[e_0,l_0]$ in which the vehicles must depart and arrive at the depot.
The time for loading/unloading a vehicle at node $i\in V_p\cup V_d$ is denoted by $s_i$. 

There is a set of heterogeneous vehicles $K$, with each $k\in K$ having a capacity $Q_k$. 
Let $F = \{0, 1, ..., z-1\}$ be the set of vehicle accessible regions, where $z$ is the number of regions. 
If two nodes $i,j \in V$ belong to the same region $f \in F$, then the edge $ij \in E^s$.
Consider a set of machines $H$ dedicated exclusively to carrying the vehicles between intermediate locations, called stations, to transpose certain obstacles.
\okC{We remark that machine stations correspond to locations in the geographical space, but they are not explicitly formalized as nodes \okC{of the graph $G$}.}
Besides, let $H_{ij} \subseteq H$ be the set of machines that can be used in the traversal of edge $e \in E^m, \ e = \{i,j\}$.
Let $F_h \subseteq F$ be the set of regions in which machine $h \in H$ has a station. 
Each machine $h \in H$ has at most one station in each region $f \in F$. 
Notice that $|F_h| \geq 2, \; \forall h \in H$, so that each machine \okC{serves} at least two regions. 
We consider the following assumptions: there is at least one machine in $H$ that has a station in all the regions, a machine can only carry one vehicle at a time, no preemption is allowed, and a single machine is necessary to cross an edge $ij \in E^m$.
Let $f^h_i$ be the station of machine $h\in H$ corresponding to node $i \in V$. 
The parameter $f^h_i$ maps $h \in H, i \in V \rightarrow \{0,...,|F_h|-1\}$, that is, the station points of machine $h \in H$. 
If the machine doesn't have a station in the region of node $i$, $f^h_i$ is not defined.
Denote by $\bar{d}^k_{ih}$ the time to travel between node $i\in V$ and the station $f^h_i$ of machine $h$ corresponding to node $i$, using the vehicle $k\in K$. This is only defined if $f^h_i$ is specified.

The time for the machine $h\in H$ to traverse the obstacle corresponding to edge $ij \in E^m$ is denoted by $O^h_{f^h_if^h_j}$. 
It is assumed that each machine $h \in H$ begins its journey at its machine station indexed by $0$, i.e., the first machine station point.
In case the machine $h \in H$ is not positioned at the boarding station $f^h_{i'}$ of its next traversal corresponding to edge $i'j' \in E^m$, the machine must traverse the edge with dead freight from its current station to the station $f^h_{i'}$.
This current station for machine $h$ can be either the station indexed by $0$ (travel duration: $O^h_{0 f^h_{i'}}$), if the machine $h$ is starting its journey, or $f^h_{j}$ (travel duration: $O^h_{f^h_{j} f^h_{i'}}$), if the machine $h$ just finished traversing the obstacle corresponding to edge $ij \in E^m$.
\okC{Besides, note that despite each machine $h \in H_{ij}$ is used to traverse an obstacle corresponding to edge $(i,j) \in E^m$ at most one time, it can travel multiple times over the same obstacle, either transporting or not a vehicle $k \in K$. 
Therefore, the segment representing the obstacle between two machine stations is not explicitly modeled as an edge in graph $G$. 
Rather, it represents part of the traversal associated with an edge in $E^m$ that connects two customer locations.}

Finally, traversing the edge $ij \in E^s$ with a vehicle $k\in K$ takes time $\hat{d}^k_{ij}$.
The time $\hat{d}^k_{ij}$ for an edge $ij \in E^m$ defines a lower bound in real-time to traverse the edge, as it does not consider the schedule of the machine. 
This lower bound is set to $\hat{d}^k_{ij} = \texttt{min}_{h \in H_{ij}} (\bar{d}^k_{ih} + O^h_{f^h_if^h_j} + \bar{d}^k_{jh})$.
It is assumed that these traversal times respect the triangular inequalities.

Denote by $C_k$ the completion time of vehicle $k \in K$ and let $C_{total} = \sum_{k\in K}C_k$ denote the total completion time. The problem thus consists of defining a feasible routing and scheduling \okC{plan} of the vehicles and machines that minimizes the total completion time.

\okC{
The problem assumptions are summarized in Table \ref{tab:problem_assumptions}.
Table \ref{tab:parameters} (see \ref{app:notations}) summarizes the problem parameters.
}

\vspace{-.8cm}
\okC{
\begin{adjustwidth}{-1.0 cm}{-1.0 cm}
\centering
\begin{threeparttable}[htb]
\footnotesize
\caption{\okC{Summary of the problem assumptions}\label{tab:problem_assumptions}}
\begin{tabular}{p{0.02\textwidth} p{0.75\textwidth}}
\toprule
\# & Assumption \\
\midrule
1  & Each request is assigned to exactly one vehicle (no transshipment). \\
2  & Vehicles traverse edges within a region autonomously. \\
3  & Inter-region edges require scheduling a machine to transport the vehicle over the obstacle. \\
4  & Services must start within a time window; finishing later is allowed. \\
5  & Vehicles can wait at a customer's location until its time window opens. \\
6  & Vehicles can load multiple requests if they satisfy capacity constraints. \\
7  & Each machine has at most one exclusive station per region. \\
8  & Machines are eligible on inter-region edges if their stations exist in both regions. \\
9  & At least one machine has a station in all regions. \\
10 & Machines carry only one vehicle at a time. \\
11 & Each vehicle uses only one machine to cross an inter-region edge. \\
12 & Each machine has an initial station; return not required. \\
13 & Traversal times satisfy the triangle inequality. \\
14 & No fixed vehicle-machine assignment; vehicles can use different machines. \\
15 & Machines may move with dead freight between stations. \\
16 & Machines traverse each edge at most once, but may travel multiple times between regions. \\
17 & Stations are not associated to nodes in the graph, and their connections are not edges. \\
\bottomrule
\end{tabular}
\label{tab:problem_constraints}
\end{threeparttable}
\end{adjustwidth}
}

\begin{figure}[!ht]
    \centering
    \includegraphics[width=.6\linewidth]{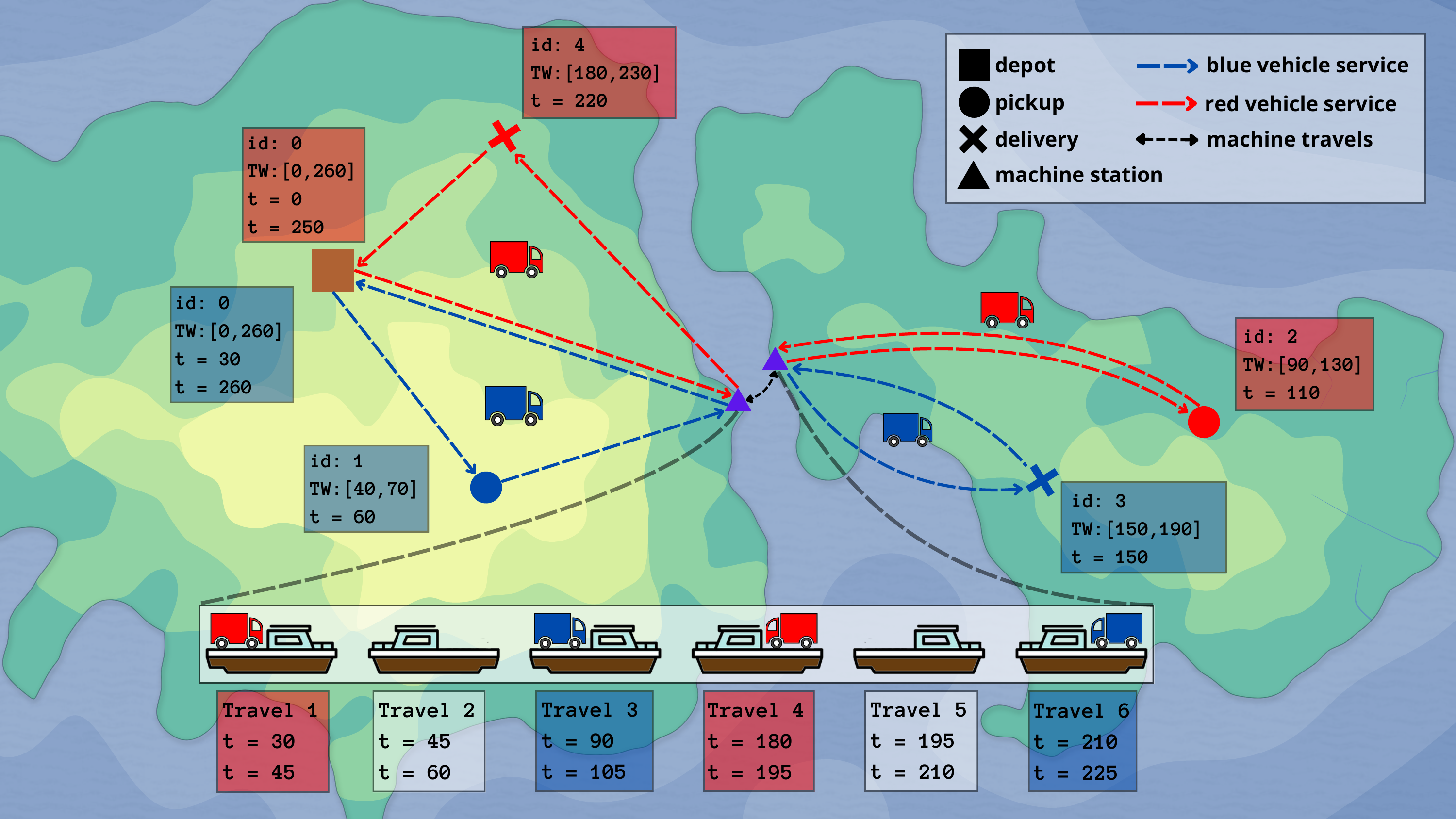}
    \caption{PDPTW-SE solution example.}
    \label{fig:pdptw-se_example_sol}
\end{figure}

Figure \ref{fig:pdptw-se_example_sol} illustrates a solution for an instance with two requests ($n = 2$), two vehicles ($|K|=2$), two regions ($z=2$), and one machine ($|H|=1$).
The red vehicle starts its journey at time 0, moving toward pickup node 2, which is located in another region.
Therefore, it arrives at the machine station at time 30, and the machine transports the vehicle to the other region, arriving at time 45.
Then, the red vehicle continues its journey to node 2, while the machine returns to its previous station to wait for the blue vehicle.
The blue vehicle leaves the depot at time 30, picks up the goods in node 1 at time 60, and meets the machine at its station at time 90.
Next, the machine transports the blue vehicle to the other region, arriving at time 105.
While the blue vehicle moves to node 3, the red vehicle picks up the goods at node 2 at time 110 and returns to the machine station.
The blue vehicle arrives before the TW, but it waits to deliver the goods at node 3 at time 150, and while it returns to the machine station, the machine transports the red vehicle to the other region.
Afterward, the machine returns to its previous station to transport the blue vehicle that was already waiting at the machine station, while the red vehicle continues its journey to deliver the goods at node 4.
Finally, both vehicles return to the depot.
The total completion time of this solution is 480, as the blue and red vehicles complete their journey in 230 and 250, respectively.

\section{\okC{Related works}}
\label{sec:relatedworks}

Pickup and delivery problems (PDPs) have been extensively discussed in the literature. \cite{SavSol95} provided the first survey on this topic, examining the characteristics that differentiate the PDP from traditional vehicle routing problems and describing the problem types and current solution methods, focusing on deterministic models.
We refer to \cite{BerCorGriLap07} and \cite{ParDoeHar08} for literature surveys of the solution methodologies for PDPs.

A particular case of the PDP is the dial-a-ride problem (DARP), which considers the transportation of people \citep{BodBer79,CulJarRat81}. \cite{Ahuja2003} considered the convex cost integer dual network flow problem that can be applied for the DARP. \cite{Cordeau2007} summarized the most important models and algorithms developed for the DARP up to 2007. \cite{HO2018395} surveyed DARPs and discussed problem variants, solution methodologies, benchmark instances, and further research developments.

\cite{DumDesSou91} introduced the PDP with time windows (PDPTW). 
Since the first work on this topic, several others have been developed to improve the solution approaches and provide new benchmark instances. The current state-of-the-art methods for the PDPTW are the works of \cite{Sartori2020}, which extended the metaheuristic algorithm proposed by \cite{CURTOIS2018151}, and that of \cite{Christiaens2020}. The first one presented the use of well-known heuristic components embedded in an iterated local search (ILS) framework. The components include the adaptive guided ejection search (AGES) of \cite{CURTOIS2018151}, large neighborhood search, and a mathematical programming phase using a set partitioning formulation of the PDPTW, combined with a new perturbation mechanism to improve the results of the AGES. The second work proposed a ruin-and-recreate method for some vehicle routing problems, including the PDPTW.
\okC{Regarding exact methods for the PDPTW, \cite{FURTADO2017334} extended the classical two-index formulation of the VRPTW for the PDPTW with a homogeneous fleet of vehicles.}
The standard benchmarks of the PDPTW were proposed by \cite{LiLim03} and \cite{Ropke2009}, which consist of locations defined in a two-dimensional Cartesian plane with symmetrical travel times and artificial time windows.
Recently, \cite{Sartori2020} proposed a new set of instances with more realistic configurations.

\okC{Several problems related to the PDPTW have also been studied in the literature, with many considering the use of different vehicle types. \cite{Dahle2019} considered a crowd-shipping solution in which a transportation company not only has its own fleet of vehicles to service requests but may also use the services of occasional drivers. These drivers receive a small compensation for taking a detour to serve one or more transportation requests.} \okC{\cite{Ghilas2016b} introduced the PDPTW with scheduled lines (PDPTW-SL), which integrates freight and public transportation. The PDPTW-SL includes two transport options: direct and indirect shipments. In a direct shipment, a request's origin and destination points are visited using only one pickup and delivery (PD) vehicle. On the other hand, the indirect shipment implies that a request is picked up by a PD vehicle and transferred to a transfer node, which is assumed to be located nearby. From there, the request continues its journey on fixed scheduled lines (SLs). Afterward, the request is picked up again by another PD vehicle to be delivered to its destination point.} \okC{\cite{Soriano2018} studied the two-region multi-depot pickup and delivery problem (2R-MDPDP), which is a two-level and multi-period transportation problem. 
In this problem, heavy-duty vehicles (HDVs) handle inter-region transfers between depots, while light-duty trucks (LDTs) perform local pickups and deliveries within each region. 
Both inter- and intra-region requests are considered, and the routing and scheduling decisions are interdependent, as the depots used by HDVs determine the LDT operations. 
The objective is to minimize transportation costs and the number of short-haul routes.}

Notice that the 2R-MDPDP and the PDPTW-SL are similar to the PDPTW-SE because they also consider that pickups and deliveries can be separated in different regions. The difference is that the PDPTW-SE does not consider transshipment, i.e., each request (pickup and delivery) in the PDPTW-SE is \okC{served} by only one vehicle, as in the PDPTW. Besides, some machines are responsible for transporting the vehicles between two regions.
\cite{Mourad2021} moves in this direction by considering a similar scenario from the PDPTW-SL, where the PD vehicles are PD robots, which allow them to be transported through the SLs, extending their geographical reach.
However, in the PDPTW-SE, these machines are limited and \okC{have a flexible schedule}, i.e., they do not operate in a fixed schedule. Besides, each PD robot can only serve one freight request at a time.

According to \cite{Drexl2012}, classical vehicle routing problems require synchronization between vehicles concerning which vehicle visits which customer.
However, integrating more than one vehicle to fulfill a task incorporates additional synchronization requirements regarding spatial, temporal, and load aspects.
This leads to an \textit{interdependence problem}, where a change in one route may affect the feasibility of other routes that are ``linked'' to it.
\cite{SOARES2024817} summarizes the five categories of synchronization of \cite{Drexl2012} into two main categories: operation and movement synchronization.
The first refers to the interdependencies between tasks on different routes that need to be performed within some temporal offset.
The second refers to the interdependencies between sequences of tasks among different routes.
The PDPTW-SE can be categorized under \textit{movement synchronization en route}, as a vehicle may join a machine at its departure station to traverse some edge and separate at the arrival station to continue its journey.
\okC{Pickup and delivery problems} with movement synchronization en route can be found, for instance, in the Active-Passive VRP \citep{Meisel2012, Tilk2018}, the PDPTW-SL with PD robots \citep{Mourad2021}, the PDP with drones \citep{MULUMBA2024103377}, \okC{and the fleet sizing and routing problem with synchronization for AGVs \citep{AziCotCoe22}}.
Interested readers are referred to \cite{Drexl2012} and \cite{SOARES2024817} for literature surveys of vehicle routing problems with synchronization.

\okC{
The PDPTW-SE is also related to the Active-Passive VRP (APVRP) studied by \cite{Meisel2012} and \cite{Tilk2018}.
In the APVRP, the transport requests are also comprised of pickup and delivery tasks.
To accomplish these tasks, there are two types of transport resources, namely, passive and active means of transport.
The passive means cannot move autonomously, thus requiring an active means to transport them between two locations.
Therefore, the active means is responsible for carrying a passive means to a pickup location, where it is loaded, or to carry a loaded passive means to a delivery location, where it is dropped off and then unloaded.
Besides, there is no fixed assignment of active means to passive means, i.e., the passive means can be carried by different active means over time.
Consequently, each travel of a passive means is synchronized with a travel of a vehicle along the same arc of the logistics network.
This enables a higher vehicle utilization, although it makes the planning more complex.
In the PDPTW-SE, we can view the machines as the active means, carrying the PD vehicles between two regions.
However, the PD vehicle is not a passive means, given that it can move autonomously in the same region and perform multiple services independently. Therefore, we define it as a partially active means of transport.}

The PDPTW-SE \okC{is also somewhat related} to van-based robot delivery problems, such as the two-echelon van-based robot hybrid pickup and delivery system \okC{(2E-VRHPD)} \citep{Yu2021}\okC{, and to PDPs with drones, such as the drone-assisted pickup and delivery problem (DAPDP) \citep{MULUMBA2024103377}}.
In the 2E-VRHPD, both vans and robots can pick up and deliver the goods, but only the robots can reach all areas due to van access restrictions.
Thus, to attend to these restricted areas, the van carries a robot and stops at parking nodes to drop off and/or pick up the robot.
\okC{Analogously, in the DAPDP, the trucks and the drones are respectively similar to the vans and the robots, except that the truck can serve any customer and the drone is limited by the amount of cargo it can carry, and to requests that do not require a signature.}
However, there are many differences from the PDPTW-SE, and some are highlighted in sequence.
In the 2E-VRHPD\okC{/DAPDP}, each van\okC{/truck must carry one and only} one robot\okC{/drone}, and each robot\okC{/drone} depends on a van\okC{/truck} to leave the depot and serve the customers.
Conversely, in the PDPTW-SE, the machines are limited; any vehicle can be assigned to a machine, and a PD vehicle is only synchronized with machines when traversing obstacles, i.e., it is independent when traversing edges \okC{within} a region.
Besides, each robot\okC{/drone} in the 2E-VRHPD\okC{/DAPDP} is picked up by the same van\okC{/truck} that dropped it off, while each PD vehicle in the PDPTW-SE can switch between machines.
Moreover, in the PDPTW-SE, the machines are dedicated exclusively to transporting the PD vehicles between two regions\okC{, not being responsible for serving any customer}.
Finally, the machines do not depart or return to the depot, as they can only traverse between predefined points outside the regions, which we call machine stations.
Instead, each machine departs from its initial station, can visit any station multiple times, and is not required to return to the initial station.

\okC{
Routing problems involving pickup and delivery can also be found in healthcare systems. 
To design AGV systems, strategic (system design), tactical (fleet sizing, flow-path layout design), and operational issues must be taken into consideration \citep{Vis06}.
\cite{AziCotCoe22} studied the usage of AGVs, proposing the fleet sizing and routing problem with synchronization for AGVs with dynamic demands (FSRPS-AGV), which combines tactical and operational decisions in the same setting. 
In the FSRPS-AGV, AGVs carry out automated transportation tasks \okC{corresponding to different requests}.
Moreover, an AGV must transport a specific type of cart for each task\okC{, and the tasks of a request may be distributed among multiple AGVs.
Notably, the FSRPS-AGV and the APVRP share similarities in terms of the means of transportation: a cart is a passive means, while the AGV is an active means.
Consequently, the differences observed between the APVRP and the PDPTW-SE also apply to the~FSRPS-AGV.}}

\okC{
Another PDP variation in healthcare systems is the dynamic pickup and delivery problem with soft time windows (DPDPSTW) proposed by \cite{Aziez2023}, which considered the transport of patients between different healthcare units in large hospitals.
Those hospitals have multiple buildings, each with a set of floors.  
The patients are transported by a heterogeneous fleet of vehicles composed of two vehicle types: ambulances, limited to one patient, and electric vehicles, which can load up to four patients at a time.
In this problem, the authors considered a dynamic version of the PDP, where there is a deterministic set of requests at the beginning of the day and other requests that can arrive throughout the day. 
They also consider soft time windows, that is, it's allowed to complete the request delivery after the end of the time window, although there is a penalty for the lateness, which is weighted based on the urgency level of the request. 
The objective is to minimize all vehicles' total weighted lateness and travel times.} There are some key differences between the PDPTW-SE and the DPDPSTW. 
First, in the PDPTW-SE, we consider the transportation of goods, not people.
Second, in the PDPTW-SE, we also want to handle the scheduling time for transporting the goods between different floors. 
The DPDPSTW incorporates the time of traveling between floors into the total travel time between two hospital locations, so it becomes the traditional pickup and delivery~problem.

\okC{
Table \ref{tab:literature_gaps} summarizes the literature gaps considering related pickup and delivery problems.
We classified the major problems discussed in this section in six main categories.
Category ``region'' indicates if the nodes of the problem are located in a single region or in multiple regions. 
We classify a problem as single region if all \okC{pickup and} delivery vehicles can move autonomously between any two locations.
Otherwise, it is classified as multiple region.
Category ``ship. type'' designates whether a request is \okC{always} assigned to only one vehicle (direct shipment) or \okC{it may be assigned} to more than one vehicle (transshipment).
Category ``sched. lines'' shows if some vehicles have a fixed/flexible schedule line.
We classify the problem as \okC{having fixed schedule lines} in the case a vehicle has fixed departure times.
Otherwise, \okC{it has flexible schedule lines}, i.e., the vehicles only departure according to a demand.
Category ``sync. type'' indicates whether the problem \okC{is classified under operation or movement synchronization}.
Category ``means of transp.'' designates how vehicles move, i.e., whether autonomously (active), partially autonomous (partially active), or dependently (passive).
Finally, category ``pair. type'' shows the way different vehicles can be paired.
We classify the problem as a fixed pairing type when there is a fixed assignment of a vehicle to another vehicle.
Otherwise, it is classified as flexible pairing type, because vehicles can switch pairs.
}

\textcolor{black}{
\begin{adjustwidth}{-1.0 cm}{-1.0 cm}
\centering
\scriptsize
\begin{threeparttable}[!htb]
\setlength{\tabcolsep}{4pt}
\caption{\textcolor{black}{Literature gaps considering related pickup and delivery problems}}\label{tab:literature_gaps}
\begin{tabular}{l cc cc cc cc ccc cc}
\toprule
&\multicolumn{2}{c}{region} &\multicolumn{2}{c}{ship. type} &\multicolumn{2}{c}{sched. lines} &\multicolumn{2}{c}{sync. type} &\multicolumn{3}{c}{means of transp.} &\multicolumn{2}{c}{pair. type}\\
\cmidrule(lr){2-3} \cmidrule(lr){4-5} \cmidrule(lr){6-7} \cmidrule(lr){8-9} \cmidrule(lr){10-12} \cmidrule(lr){13-14} 
Work &single &multi &direct &trans. &fixed &flex. &oper. &mov. &act. &part. act. &pass. &fixed &flex. \\
\midrule
\cite{SavSol95}             &\checkmark &- &\checkmark &- &- &- &- &- &\checkmark &- &- &- &- \\
\arrayrulecolor{black!20}\hline
\cite{BodBer79}             &\checkmark &- &\checkmark &- &- &- &- &- &\checkmark &- &- &- &- \\
\arrayrulecolor{black!20}\hline
\cite{DumDesSou91}          &\checkmark &- &\checkmark &- &- &- &- &- &\checkmark &- &- &- &- \\
\arrayrulecolor{black!20}\hline
\cite{Dahle2019}            &\checkmark &- &\checkmark &- &- &- &- &- &\checkmark &- &- &- &- \\
\arrayrulecolor{black!20}\hline
\cite{Ghilas2016b}          &- &\checkmark &- &\checkmark &\checkmark &- &\checkmark &- &\checkmark &- &- &\checkmark &- \\
\arrayrulecolor{black!20}\hline
\cite{Soriano2018}          &- &\checkmark &- &\checkmark &- &- &\checkmark &- &\checkmark &- &- &- &- \\
\arrayrulecolor{black!20}\hline
\cite{Mourad2021}           &- &\checkmark &\checkmark &- &\checkmark &- &- &\checkmark &\checkmark &\checkmark &- &\checkmark &- \\
\arrayrulecolor{black!20}\hline
\cite{Meisel2012}           &- &\checkmark &- &\checkmark &- &\checkmark &\checkmark &\checkmark &\checkmark &- &\checkmark &- &\checkmark \\
\arrayrulecolor{black!20}\hline
\cite{MULUMBA2024103377}    &- &\checkmark &\checkmark &- &- &- &\checkmark &\checkmark &\checkmark &\checkmark &- &\checkmark &- \\
\arrayrulecolor{black!20}\hline
\cite{AziCotCoe22}          &- &\checkmark &\checkmark &- &- &- &\checkmark &\checkmark &\checkmark &- &\checkmark &- &- \\
\arrayrulecolor{black!20}\hline
\cite{Yu2021}               &- &\checkmark &- &\checkmark &- &- &\checkmark &\checkmark &\checkmark &\checkmark &- &\checkmark &- \\
\arrayrulecolor{black!20}\hline
\cite{Aziez2023}            &- &\checkmark &\checkmark &- &- &- &- &- &\checkmark &- &- &- &- \\
\arrayrulecolor{black!20}\hline
\textbf{This work}                   &- &\checkmark &\checkmark &- &- &\checkmark &- &\checkmark &\checkmark &\checkmark &- &- &\checkmark \\
\arrayrulecolor{black}
\bottomrule
\end{tabular}
\begin{tablenotes}
    \item \underline{Abbreviations:} ship. = shipment; trans. = transhipment; flex. = flexible; sched. = schedule; sync. = synchronization; oper. = operation; mov. = movement, trans. = transport; act. = active; part. act. = partial active; pass. = passive; pair. = pairing.
\end{tablenotes}
\end{threeparttable}
\end{adjustwidth}
}

\section{MIP formulation}
\label{sec:formulation}

In this section, we present an arc-based MIP formulation for the PDPTW-SE.
Consider the directed graph $G'= (V',A)$, where $V' = V \cup \{2n+1\}$ and node $2n+1$ is a copy of the depot (node 0).
\okC{The set of arcs is defined as $A=\{(i,j)\ | \ i,j\in V', \ i \neq j\}$.}
For every $k\in K$, define: 
$d^k_{ij}=d^k_{ji}=\hat{d}^k_{ij}$ for $ij \in E$ with $i,j \in V_p \cup V_d$; $d^k_{i,2n+1} = \hat{d}^k_{0i}$ for $i \in V_d$; 
$d^k_{0j} = \hat{d}^k_{0j} $ for $j \in V_p$;
and $d^k_{0,2n+1} = 0$.
Additionally, let $A^s$ and $A^m$ be the arcs associated with the edges in $E^s$ and $E^m$, respectively.

In what follows, we divide the presentation of the formulation into its routing and scheduling parts. 
Section~\ref{sec:miprouting} describes the routing part, while Section~\ref{sec:mipscheduling} details the scheduling part.
Section~\ref{sec:completeformulation} presents the complete formulation.
\okC{Section~\ref{sec:mippreprocessing} shows the preprocessing step.
Section~\ref{sec:mip_tightening} describes the proposed valid inequalities.}
\subsection{Routing constraints}
\label{sec:miprouting}

We define the following routing-related variables.
The binary variable $x^k_{ij}$ is equal to one if the vehicle $k$ traverses the arc $(i,j)\in A$ \okC{(i.e., goes directly from node $i \in V$ to node $j \in V' \setminus \{0\}$)}, zero otherwise.
\okC{We state that when the binary variable $x^k_{0, 2n+1}$ is equal to one, the vehicle $k$ is inactive, i.e., the vehicle is not used.}
The continuous variable $z^k_i$ denotes the weight of the vehicle $k\in K$ just after leaving node $i\in V'$. 
The routing part of the problem can thus be formulated as:
\begingroup
\allowdisplaybreaks
\begin{align}
& \sum_{j \in V_p \cup \{2n+1\}} x^k_{0j} = 1, \qquad \forall \ k \in K. \label{mip:r:10} \\
& \sum_{j:(j,i) \in A} x^k_{ji} - \sum_{j:(i,j)\in A} x^k_{ij} = 0, \qquad \forall \ k \in K, \ i \in V_p \cup V_d. \label{mip:r:15} \\
& \sum_{j \in V_d \cup \{0\}} x^k_{j,2n+1} = 1, \qquad \forall \ k \in K. \label{mip:r:20} \\
& \sum_{k \in K} \sum_{j:(j,i) \in A} x^k_{ji} = 1, \qquad \forall \ i \in  V_p \cup V_d. \label{mip:r:30} \\
& \sum_{j:(j,i) \in A} x^k_{ji} = \sum_{j:(j,n+i) \in A} x^k_{j,n+i}, \qquad \forall \ k \in K, \ i \in V_p. \label{mip:r:35} \\
& z^k_0 = 0, \qquad \forall \ k \in K. \label{mip:r:40} \\
& z^k_j \geq z^k_i + q_{j} - M_1(1-x^k_{ij}), \qquad \forall \ k \in K, \ (i,j) \in A. \label{mip:r:45} \\
& z^k_j \leq z^k_i + q_j + M_1(1-x^k_{ij}), \qquad \forall \ k \in K, \ (i,j) \in A. \label{mip:r:47}  \\
& z^k_i \leq \min\{Q^k, \max\{0,Q^k+ {q_i}\}\} \sum_{j:(j,i) \in A} x^k_{ji}, \qquad \forall \ k \in K, \ i \in V_p \cup V_d. \label{mip:r:50} \\
& z^k_i \geq {q_i} \sum_{j:(j,i) \in A} x^k_{ji}, \qquad \forall \ k \in K, \ i \in V_p. \label{mip:r:55} \\
& x^k_{ij} \in \{0,1\}, \qquad \forall \  k \in K, \ (i,j) \in A. \label{mip:r:60} \\
& z^k_{i} \geq 0 , \qquad \forall \  k \in K, \ i \in V'.  \label{mip:r:70}
\end{align}
\endgroup

Constraints~\eqref{mip:r:10}-\eqref{mip:r:20} are balance constraints for the vehicles.
They establish that every vehicle departs from the depot (constraints \eqref{mip:r:10}), leaves a node whenever they arrive at it (constraints \eqref{mip:r:15}), and goes back to the depot (constraints \eqref{mip:r:20}).
Constraints~\eqref{mip:r:30} ensure that every pickup and delivery node is visited.
Constraints~\eqref{mip:r:35} guarantee that the pickup and delivery nodes corresponding to a given request are visited by the same vehicle.
Constraints~\eqref{mip:r:40}-\eqref{mip:r:55} enforce that the capacities of the vehicles are respected.
They imply that the vehicles are empty at the depot (constraints \eqref{mip:r:40}), their cargo weights are updated whenever they go from one node to another (constraints \eqref{mip:r:45}-\eqref{mip:r:47}), their cargo weight never exceeds their capacities (constraints \eqref{mip:r:50}), and the weights of every visited node are taken into account (constraints \eqref{mip:r:55}).
Constraints~\eqref{mip:r:60}-\eqref{mip:r:70} restrain the domains of the $\boldsymbol{x}$ and $\boldsymbol{z}$ variables. 
$M_1$ is a large enough number that can be set as $\texttt{max}_{k \in K}\{Q_k\} + \texttt{max}_{i \in V_p}\{q_i\}+1$.
\subsection{Scheduling constraints}
\label{sec:mipscheduling}

We introduce the following vehicle scheduling-related variables.
The continuous variable $t_i$ provides the starting time to serve node $i\in V_p \cup V_d$.
The continuous variables $t^k_0$ and $t^k_{2n+1}$ represent, respectively, the departing and arrival time of the vehicle $k\in K$ at the depot.
The continuous variable $C_k$ gives the completion time of the vehicle $k \in K$.
Furthermore, consider the following machine scheduling-related variables. 
The binary variable $\phi^h_{ij}$ is equal to one if machine $h \in H_{ij}$ is used to cross the arc $(i,j)\in A^m$, zero otherwise.
The binary variable $\gamma^h_{iji'j'}$ is equal to one if the use of machine $h \in H_{ij} \cap H_{i'j'}$ to cross arc $(i,j)\in A^m$ precedes its use to cross arc $(i',j')\in A^m$, with $(i',j') \neq (i,j)$, zero otherwise.
The continuous variable $\alpha^h_{ij}$ denotes the start time to use machine $h\in H_{ij}$ to cross arc $(i,j)\in A^m$.
In this way, the scheduling part of the problem can be formulated as:

\begingroup
\allowdisplaybreaks
\small
\begin{align}
& t_j \geq t_i + s_i + d^k_{ij} - M_2 (1-x^k_{ij}), \quad \forall \ k \in K,\ (i,j)\in A, \ \textrm{with }i\neq 0 \textrm{ and }j\in V_p \cup V_d. \label{mip:s:10} \\
& t_j \geq t^k_0 + d^k_{0j} - M_3 (1-x^k_{0j}), \quad \forall \ k \in K,\ (0,j)\in A, \ \textrm{with }j\in {V_p}. \label{mip:s:10b} \\
& t_{n+i} \geq t_i + s_i + \sum_{k \in K}\sum_{l:(l,i) \in A} d^k_{i,n+i}x^k_{li} , \quad \forall \ i \in V_p. \label{mip:s:15} \\
& \sum_{h \in H_{ij}}\phi^h_{ij} = \sum_{k \in K} x^k_{ij}, \quad \forall \ (i,j) \in A^m.  \label{mip:sl:05} \\
& \alpha^h_{ij} \geq t_i + s_i + \bar{d}^k_{ih} - M_4 (2-\phi^h_{ij}- x^k_{ij}),\quad \forall \ h \in H_{ij},\ k \in K,\ (i,j) \in A^m, \ i \neq 0. \label{mip:sl:10} \\
& \alpha^h_{0j} \geq t^k_0 + \bar{d}^k_{0h} - M_5 (2-\phi^h_{0j} - x^k_{0j}), \quad \forall \ h \in H_{0j},\ k \in K,\ (0,j) \in A^m, \ j \in V_p.\label{mip:sl:12} \\
& t_j \geq \alpha^h_{ij} + O^h_{f^h_if^h_j} + \bar{d}^k_{jh} - M_6 (2-\phi^h_{ij}-x^k_{ij}), \quad \forall \ h \in H_{ij},\ k \in K,\ (i,j) \in A^m,\ j\in V_p \cup V_d. \label{mip:sl:15} \\
& \gamma^h_{iji'j'} + \gamma^h_{i'j'ij} \geq  \phi^h_{ij} + \phi^h_{i'j'} - 1, \quad \forall \ h \in H_{ij} \cap H_{i'j'},\ (i,j),(i',j') \in A^m, (i,j) < (i',j'). \label{mip:sl:20} \\
& \gamma^h_{iji'j'} \leq  \phi^h_{ij}, \quad \forall \ h \in H_{ij} \cap H_{i'j'},\ (i,j),(i',j') \in A^m,  (i,j) \neq (i',j'). \label{mip:sl:25}\\
& \gamma^h_{i'j'ij} \leq  \phi^h_{ij}, \quad \forall \ h \in H_{ij} \cap H_{i'j'},\ (i,j),(i',j') \in A^m,  (i,j) \neq (i',j'). \label{mip:sl:30} \\ 
\begin{split}
    & \alpha^h_{i'j'} \geq \alpha^h_{ij} + O^h_{f^h_if^h_j} + O^h_{f^h_{j}f^h_{i'}} - {M_7}(1-\gamma^h_{iji'j'}), \\
    & \phantom{SPAAAAACE} \forall \ h \in H_{ij} \cap H_{i'j'},\ (i,j),(i',j') \in A^m, (i,j) \neq (i',j').\label{mip:sl:35}\\
\end{split}\\
& \alpha^h_{ij} \geq O^h_{0f^h_i} - {M_8}(1-\phi^h_{ij}), \quad \forall \ h \in H_{ij} \cap H_{0i},\ (i,j) \in A^m. \label{mip:sl:37} \\
& t^k_{2n+1} \geq t_i + s_i + d^k_{i,2n+1} - M_2 (1 - x^k_{i,2n+1}), \quad \forall \ k \in K,\ (i,2n+1) \in A, \ {i \in V_d}. \label{mip:sl:55} \\
\begin{split} 
    & t^k_{2n+1} \geq \alpha^h_{i,2n+1} + O^h_{f^h_if^h_{2n+1}} + \bar{d}^k_{{2n+1}h} - M_4 (2-\phi^h_{i,2n+1}-x^k_{i,2n+1}), \\
    & \phantom{SPAAAAACE} \forall \ h \in H_{i, 2n+1},\ k \in K,\ (i,2n+1) \in A^m, \ i \in V_d.\label{mip:sl:57}\\
\end{split}\\
& C_k \geq t^k_{2n+1} - t^k_0, \quad \forall \ k \in K. \label{mip:sl:70} \\
& e_i \leq t_i \leq l_i, \quad \forall \ i \in V_p \cup V_d. \label{mip:s:20} \\
& e_0 \leq t^k_0 \leq t^k_{2n+1} \leq l_0 , \quad \forall \ k \in K .\label{mip:s:25} \\
& C_k \geq 0, \quad \forall \ k \in K. \label{mip:sl:27} \\
& \alpha^h_{ij} \geq 0 , \quad \forall \  h \in H_{ij}, \ (i,j) \in A^m. \label{mip:sl:45} \\
& \phi^h_{ij} \in \{0,1\} , \quad \forall \  h \in H_{ij}, \ (i,j) \in A^m. \label{mip:sl:40} \\
& \gamma^h_{iji'j'} \in \{0, 1\} , \quad \forall \ h \in H_{ij} \cap H_{i'j'},\ (i,j), (i',j') \in A^m. \label{mip:sl:50}
\end{align}
\endgroup

Constraints~\eqref{mip:s:10}-\eqref{mip:s:10b} ensure that whenever an arc is traversed by a vehicle, a lower bound on the time to serve its destination node is determined by the times of its origin node. Constraints~\eqref{mip:s:15} guarantee that a pickup node of a request is visited before its delivery node.
Constraints~\eqref{mip:sl:05} imply that a machine is used for an arc if and only if it is traversed by a vehicle.
Constraints~\eqref{mip:sl:10}-\eqref{mip:sl:12} define that whenever a machine is used to traverse an arc, a lower bound on its starting time on the machine is determined by the times corresponding to its origin node.
Constraints~\eqref{mip:sl:15} determine that whenever a machine is used to traverse an arc, the machine times define a lower bound on the starting time of its destination node.
Constraints~\eqref{mip:sl:20}-\eqref{mip:sl:30} define that there is an order between the traversal of two distinct arcs if and only if they are both scheduled on the same machine.
Constraints~\eqref{mip:sl:35} define a lower bound on the time to traverse an arc whenever it is preceded by another arc. 
Note that the machine $h$ moves with dead freight from its current station $f^h_{j}$ to its next boarding station $f^h_{i'}$, in case $f^h_{j} \neq f^h_{i'}$.
Constraints~\eqref{mip:sl:37} define a lower bound on the time to traverse an arc (whenever it is traversed) based on the time that the machine takes between its initial station and the initial travel station.
Constraints~\eqref{mip:sl:55}-\eqref{mip:sl:57} determine lower bounds on the times the vehicles arrive at the depot.
Constraints~\eqref{mip:sl:70} set lower bounds on the completion times of the vehicles.
Constraints~\eqref{mip:s:20}-\eqref{mip:sl:50} restrain the domains of the $\boldsymbol{t}$,  $\boldsymbol{t_0}$, $\boldsymbol{t_{2n+1}}$, $\boldsymbol{C}$, $\boldsymbol{\alpha}$, $\boldsymbol{\phi}$,  and $\boldsymbol{\gamma}$ variables. 
\okC{$M_i, \ i \in \{2,\ldots, 8\}$,} are sufficiently large numbers\okC{, and p}ossible values are shown in \mbox{\ref{app:big_m_bounds}}.
\subsection{The complete formulation}
\label{sec:completeformulation}

A formulation for \okC{the} PDPTW-SE can thus be cast as:
{\small
\begin{align}
 \min \qquad & \sum_{k\in K} C_k.  \label{mip:obj} \\
& \eqref{mip:r:10}-\eqref{mip:sl:50}. \nonumber
\end{align}}
The objective function \eqref{mip:obj} minimizes the total completion time.
\subsection{\okC{Preprocessing}}
\label{sec:mippreprocessing}

\okC{
Given that the PDPTW-SE is a generalization of the PDPTW, we can perform a preprocessing step to eliminate several inadmissible arcs similar to \cite{DumDesSou91}.
The main difference is that the fleet in the PDPTW-SE is heterogeneous.
} 

\okC{
Denote by $P$ the set of valid paths of the graph $G'$, i.e., every two consecutive nodes of the path must be associated with an arc in $A$.
Moreover, denote by $P_w \subseteq P$ the subset of paths $p \in P$ such that the length of the path $p$ is $w \in \mathbb{N}$, i.e., the path has $w$ arcs. 
Given $k \in K$ and $p \in P_w, \ p = (v_1, \ldots, v_{w+1}), \ (v_i, v_{i+1}) \in A, \ \forall i \in [1, w],$ we denote $\Delta^k_{p}$ as the earliest arrival time of vehicle $k$ in $v_{w+1}$ when $t_{v_1} = e_{v_1}$.
We can define a general partial path $p = (v_1, v_2, \ldots, v_{w+1})$ as infeasible by solving the following recurrence and checking whether $\Delta^k_p = t_{v_{w+1}} > l_{v_{w+1}}$ for every $k \in K$:
\begin{equation*}
t_{v_i} = 
\begin{cases}
e_{v_i}, & \text{if } i = 1, \\
\max\{e_{v_i},\ t_{v_{i-1}} + s_{v_{i-1}} + d^k_{v_{i-1}v_i}\}, & \text{if } 2 \leq i \leq w+1.
\end{cases}
\end{equation*}
Remark that $d^k_{v_{i-1}v_i}$ for $(v_{i-1}, v_i) \in A^m$ in this recurrence is the lower bound time set in Section~\ref{sec:formalization}.}

\okC{
Before determining the rules to eliminate arcs, a preliminary step proposed by \cite{DumDesSou91} involves the shrinking of the time windows associated with the pickup and delivery nodes.
Every request must be serviced by at least one vehicle, which implies that for $i \in V_p$ the partial paths $(0, i, n+i)$ and $(i, n+i, 2n+1)$ are feasible for all values of $t_i \in [e_i,l_i]$ and $t_{n+i} \in [e_{n+i},l_{n+i}]$.
Therefore, we denote a new strict time window for a node $i \in V_p \cup V_d$ as $[e'_i,\ l'_i]$, whose bounds can be successively defined as follows, adapted from \cite{DumDesSou91}:
\begin{itemize}
    \small
    \item $l'_{n+i} = \texttt{min}\{l_{n+i}, l_{0} - \texttt{min}_{k \in K}\{d^k_{n+i, 2n+1}\} - s_{n+i}\}$ and $l'_{i} = \texttt{min}\{l_{i}, l_{n+i} - \texttt{min}_{k \in K}\{d^k_{i,n+i}\} - s_{i}\}$, $\forall i \in V_p$;
    \item $e'_{i} = \texttt{max}\{e_{i}, e_0 + \texttt{min}_{k \in K}\{d^k_{0, i}\}$\} and $e'_{n+i} = \texttt{max}\{e_{n+i}, e_i + s_i + \texttt{min}_{k \in K}\{d^k_{i, n+i}\}$\}, $\forall i \in V_p$.
\end{itemize}
}

\okC{
Henceforth, the PDPTW-SE constraints are used to eliminate the following infeasible arcs:
\begin{enumerate}[noitemsep, topsep=1.15pt]
    \item \textit{priority}: arcs {\small$(0,n+i), \ (n+i, i), \ (i, 0), \ (n+i, 0),$ $(2n+1, 0), \ (2n+1, i), \ (2n+1, n+i), \ \forall i \in V_p$}; \label{arc_removals:10}
    \item \textit{pairing}: arcs {\small$(i,2n+1), \ \forall i \in V_p$}; \label{arc_removals:15}
    \item \textit{vehicle capacity}: arcs {\small$(i,j), \ (j,i), \ (i, n+j), \ (j, n+i), \ (n+i, n+j), \ (n+j, n+i),$ $\forall i, j \in V_p, \linebreak i \neq j, \ q_i + q_j > \texttt{max}_{k \in K}\{Q_k\}$};\label{arc_removals:20}
    \item \textit{time windows}: arcs {\small$(i,j)$} when {\small$e'_i + s_i + \texttt{min}_{k \in K}\{d^k_{ij}\} > l'_j$, \ $i, j \in V_p \cup V_d$}; \label{arc_removals:25}
    \item \textit{time windows and pairing of requests}: given that the travel times satisfy the triangle inequality, pairing some requests can be infeasible at any vehicle route due to their time windows and precedence constraints; thus, arcs combining the pickup and delivery nodes of such requests are eliminated; for $i,j \in V_p, \ i \neq j$:\label{arc_removals:30} 
    \begin{enumerate}[noitemsep, topsep=1.15pt]
        \item arc $(i,n+j)$ is eliminated if the path $(j, i, n+j, n+i)$ is infeasible for every $k \in K$ when setting $t_j = e'_j$;
        \item arc $(n+i, j)$ is eliminated if the path $(i, n+i, j, n+j)$ is infeasible for every $k \in K$ when setting $t_i = e'_i$;
        \item arc $(i,j)$ is eliminated if the paths $(i, j, n+i, n+j)$ and $(i, j, n+j, n+i)$ are infeasible for every $k \in K$ when setting $t_i = e'_i$;
        \item arc $(n+i, n+j)$ is eliminated if the paths $(i, j, n+i, n+j)$ and $(j, i, n+i, n+j)$ are infeasible for every $k \in K$ when setting $t_i = e'_i$ and $t_j = e'_j$, respectively.
    \end{enumerate}
    \item \textit{time windows and indirect request service}: given that the travel times satisfy the triangle inequality, then for $i \in V_p$, arcs $(i,j), \ j \in V_p \cup V_d,\ i \neq j \wedge j \neq n+i$, are eliminated if path $(i, j, n+i)$ is infeasible when setting $t_i = e'_i$.\label{arc_removals:35}
\end{enumerate}
}

\okC{
Note that items \ref{arc_removals:10}-\ref{arc_removals:30} are adapted from \cite{DumDesSou91}, while item \ref{arc_removals:35} is novel.
}

\okC{
Additionally, given {\small$(i,j) \in A^m, \ (i',j') \in A^m, \ h \in H_{ij} \cap H_{i'j'}$}, we can drop a subset of the variables {\small $\gamma^h_{iji'j'}$} whenever {\small $e'_i + s_i + \texttt{min}_{k \in K}\{\bar{d}^k_{ih}\} + O^h_{f^h_i f^h_j} + O^h_{f^h_j f^h_{i'}} + O^h_{f^h_{i'} f^h_{j'}} + \texttt{min}_{k \in K}\{\bar{d}^k_{j'h}\} > l'_{j'}$} because it would be set as zero.
Consequently, we remove some constraints from (\ref{mip:sl:25})-(\ref{mip:sl:35}), and substitute the variable to zero in constraints (\ref{mip:sl:20}).
Finally, given an arc $(i,j) \in A^m$, we can also eliminate machines $h \in H_{ij}$ such that $e'_i + s_i + \texttt{min}_{k \in K}\{\bar{d}^k_{ih}\} + O^h_{f^h_i f^h_j} + \texttt{min}_{k \in K}\{\bar{d}^k_{jh}\} > l'_j$.
This step enables reducing the number of variables $\boldsymbol{\phi}$, $\boldsymbol{\gamma}$ and $\boldsymbol{\alpha}$.
}
\subsection{\okC{Formulation tightening}}
\label{sec:mip_tightening}

\okC{
In this section, some valid inequalities are proposed to cut off undesirable fractional solutions and potentially improve the performance of a branch-and-bound system when solving the formulation.
}

\subsubsection{\okC{Routing cuts}}

\okC{
Constraints \eqref{mip:vi:00} state that each vehicle cannot be inactive and also serve a pickup/delivery node.
\begin{align}
    & x^k_{0,2n+1} + \sum_{j: (j,i) \in A} x^k_{ji} \leq 1, \ k \in K, \ i \in V_p \cup V_d. \label{mip:vi:00}
\end{align}
}

\okC{
Constraints \eqref{mip:vi:05}-\eqref{mip:vi:10} enforce that each edge can only be traversed in one direction, while constraints \eqref{mip:vi:15} ensure that only one arc can precede the other.
\begin{align}
    & \sum_{k \in K}x^k_{i j} + \sum_{k \in K}x^k_{j i} \leq 1, \ (i,j), (j,i) \in A, \ i < j. \label{mip:vi:05}\\
    & \sum_{h \in H_{ij}}\phi^h_{i j} + \sum_{h \in H_{ji}}\phi^h_{j i} \leq 1, \ (i,j), (j,i) \in A^m, \ i < j. \label{mip:vi:10}\\
    & \gamma^h_{iji'j'} + \gamma^h_{i'j'ij} \leq 1,\quad \forall \ h \in H_{ij} \cap H_{i'j'},\ (i,j),(i',j') \in A^m, (i,j) < (i',j'). \label{mip:vi:15}
\end{align}
}

\okC{
Constraints \eqref{mip:vi:20} guarantee that no path of length $w$ is assigned to a vehicle when the destination node is not reachable within its time window.
\begin{align}
    & \sum^{w}_{i =1}x^k_{v_{i} v_{i+1}} \leq w-1, \ k \in K, \ p \in P_w, \ p = (v_1, \ldots, v_{w+1}), \ \Delta^{k}_{p} > l'_{v_{w+1}}. \label{mip:vi:20}
\end{align}
}

\okC{
Constraints \eqref{mip:vi:25} state that no path of length $w$ composed of arcs in $A^m$ has all its arcs assigned to any machine when the destination node is not reachable within its time window for every $k \in K$.
\small
\begin{align}
    \begin{split}
        & \sum^{w}_{i=1}\sum_{h \in H_{v_{i} v_{i+1}}}\phi^h_{v_{i}v_{i+1}} \leq w-1, \\
        & \phantom{PHANTOM} p \in P_w, \ p = (v_1, \ldots, v_{w+1}), \ (v_i, v_{i+1}) \in A^m, \ i \in [1,w], \text{ when} \ \Delta^{k}_p > l'_{v_{w+1}}, \ \forall k \in K. \label{mip:vi:25}
    \end{split}
\end{align}
}

\subsubsection{\okC{Scheduling cuts}}

\okC{
Constraints \eqref{mip:vi:30} set a lower bound on the start time of the service at a pickup/delivery node.
\begin{align}
    & t_j \geq \sum_{\substack{i: (i,j) \in A}}\sum_{k \in K} x^{k}_{ij}(e'_i + s_i + d^k_{ij}), \ j \in V_p \cup V_d. \label{mip:vi:30}
\end{align}
}

\okC{
Constraints \eqref{mip:vi:35}-\eqref{mip:vi:40} fix bounds on the start time of a machine travel based on the time windows of the arc nodes.
Notice that, given an arc $(i,j) \in A^m$ and a machine $h \in H_{ij}$, the variable $\alpha^h_{ij}$ can assume any value when $\phi^h_{ij} = 0$, because it is not taken into account in a feasible solution.
Note that these inequalities are only valid after the preprocessing step on set $H_{ij}$ (see Section \ref{sec:mippreprocessing}).
}
\okC{
\begin{align}
    & \alpha^h_{ij} \geq e'_i + s_i + \texttt{min}_{k \in K}\{\bar{d}^k_{ih}\}, \ (i,j) \in A^m, \ h \in H_{ij}. \label{mip:vi:35}\\
    & \alpha^h_{ij} \leq l'_j - \texttt{min}_{k \in K}\{\bar{d}^k_{jh}\} - O^h_{f^h_i f^h_j}, \ (i,j) \in A^m, \ h \in H_{ij}. \label{mip:vi:40}
\end{align}
}

\okC{
Constraints \eqref{mip:vi:45}-\eqref{mip:vi:50} are similar to constraints \eqref{mip:vi:35}-\eqref{mip:vi:40}, but now they include the distinct traversal times of the vehicles. 
For vehicles with equal traversal times, constraints \eqref{mip:vi:35} and \eqref{mip:vi:40} dominate constraints \eqref{mip:vi:45} and \eqref{mip:vi:50}, respectively.
However, this is not the case when dealing with vehicles with distinct traversal times.
}
\okC{
\begin{align}
    & \alpha^h_{ij} \geq e'_i + s_i + \sum_{k \in K} x^k_{ij} \bar{d}^k_{ih}, \ (i,j) \in A^m, \ h \in H_{ij}. \label{mip:vi:45}\\
    & \alpha^h_{ij} \leq l'_j - \sum_{k \in K} x^k_{ij} (\bar{d}^k_{jh} + O^h_{f^h_i f^h_j}), \ (i,j) \in A^m, \ h \in H_{ij}. \label{mip:vi:50}
\end{align}
}

\okC{
Denote by $\mu^h_{iji'} \in \{0,1\}$ the constant indicator that is equal to one if the minimum arrival time of a vehicle at machine $h$ coming from node $i'$ ({\small $e'_{i'} + s_{i'} + \texttt{min}_{k \in K}\{\bar{d}^k_{i'h}\}$}) is greater than or equal to the maximum arrival time of the machine $h$ after traversing arc $(i,j)$ ({\small $l'_{j} - \texttt{min}_{k \in K}\{\bar{d}^k_{jh}\} + O^h_{f^h_{j} f^h_{i'}}$}), otherwise, it is zero.
Hence, constraints \eqref{mip:vi:55} enforce that the start time to use a machine $h$ in a certain arc~$(i',j') \in A^m$ must always precede the start time to use a machine $h$ in another arc $(i,j) \in A^m$ when~$\mu^h_{iji'} = 1$.
\begin{align}
    \begin{split}
    & \alpha^h_{i'j'} \geq \alpha^h_{ij} + O^h_{f^h_{i}f^h_{j}} + O^h_{f^h_{j}f^h_{i'}}, \\
    & \phantom{SPACE} (i,j), (i',j') \in A^m, \ h \in H_{ij} \cap H_{i'j'}, \ (i,j) \neq (i',j'), \text{ when } \mu^h_{iji'} = 1. \label{mip:vi:55}
    \end{split}
\end{align}
}

\section{Multi-start heuristic with an LP-based improvement procedure}
\label{sec:msheurlp}

In this section, we present a multi-start heuristic for the PDPTW-SE.
The heuristic, which we denote as MSLP, employs a semi-greedy multi-start construction method embedded with an LP-based improvement procedure. The latter attempts to efficiently optimize the schedules obtained with the routes and sequences of machine traversals defined by the former. 
In what follows, Section \ref{sec:heur_notations_and_data_structs} describes the notations and data structures used in the algorithms.
After that, Section \ref{sec:heur_gen_framework} presents the general framework of the heuristic.
Section \ref{sec:greedy_insertion_heuristic} shows the greedy insertion heuristic used in the multi-start heuristic.
Section \ref{sec:semi-greedy-ins-heur} details how the greedy insertion heuristic of Section \ref{sec:greedy_insertion_heuristic} is randomized to create the semi-greedy insertion heuristic.
Finally, Section \ref{sec:lp_form} details the LP-based improvement procedure.
More details of the multi-start heuristic's auxiliary procedures are presented in \ref{app:heur_details}.
The appendix describes more deepfully the semi-greedy insertion heuristic, the insertion feasibility analysis, the vehicle travel time calculation, and the solution update.

\subsection{Notations and main data structures}
\label{sec:heur_notations_and_data_structs}

The multi-start heuristic with an LP-based improvement procedure uses the following notations and data structures.
Let $S$ be a solution for the PDPTW-SE. 
Each solution $S$ contains the following structure:  $S.feasible$ defines solution feasibility; $S.vehicles$ is a list of routes containing one route for each $k \in K$;
$S.machines$ is a list of sequences of travels containing one sequence for each $h \in H$; 
and $S.cost$ defines the sum of completion times of the vehicle routes in $S$.
The length of the vehicle route $S.vehicles[k], \ k \in K$, is denoted by $S.vehicles[k].length$.
Each position of the vehicle route $S.vehicles[k]$, $ k \in K$, is called a $visit$ 
and contains the following structure: 
$visit.i, \ i \in V'$, informs the visited node; 
$visit.st$ defines the service start time, which may also represent the departure time from node $0$ or the arrival time at node $2n+1$; 
and $visit.load$ defines the load of the vehicle.
Each position of the machine sequence of travels $S.machines[h]$, $h \in H$, is called a $travel$ and contains the following structure:
$travel.i$ and $travel.j$ define the arc $(i,j) \in A^m$ in which the machine travel is held;
$travel.k$ defines which vehicle is transported in the arc $(i,j)$;
and $travel.st$ defines the start time of the machine travel.
Note that $S.machines[h], \ h \in H,$ only stores the machine travels when transporting a vehicle, i.e., it does not consider machine travels with dead freight.

In our multi-start heuristic, we first construct an initial solution with a purely greedy insertion constructive heuristic, and then a randomized version of this heuristic is executed multiple times.
In summary, each heuristic iteratively inserts unassigned request nodes between two positions in a vehicle's route.
Therefore, each feasible position to insert a request is called an insertion candidate.
This insertion candidate has a specific structure that is further used to update solution $S$.
Let $\mathcal{C}$ be an insertion candidate for a specific request.
Each candidate $\mathcal{C}$ contains the following structure:
$\mathcal{C}.feasible$ defines candidate feasibility;
$\mathcal{C}.cost$ defines the insertion cost;
$\mathcal{C}.k$ defines the vehicle route in which the request will be inserted;
$\mathcal{C}.p_{pos}$ and $\mathcal{C}.d_{pos}$ define in which positions the pickup and delivery nodes will be inserted at $S.vehicles[k]$, respectively;
$\mathcal{C}.mtrvs$ contains the machine travels to be inserted in $S.machines$.
Each $mtrv \in \mathcal{C}.mtrvs$ contains the following structure:
$mtrv.h$ defines which machine in $H$ is being used;
$mtrv.h_{pos}$ defines the position of the machine travel in $S.machines[h]$;
$mtrv.st$ defines the machine travel start time;
$mtrv.i$ and $mtrv.j$ define the arc $(i,j) \in A^m$ in which the machine travel is held;
and $mtrv.\Delta t$ defines the time between the departure time at node $i$ and the arrival time at node $j$.
\subsection{General framework}
\label{sec:heur_gen_framework}

The general framework for the MSLP is depicted in Algorithm \ref{alg:heur_mslp}.
The algorithm receives as input the instance data $I$ and the threshold parameter $\alpha \in [0,1]$, which defines the threshold value for the cost of an insertion candidate.
The insertion data $I$ contains the parameters presented in Section~\ref{sec:formalization}, which were summarized in Table \ref{tab:parameters}.
In lines \ref{alg:heur_mslp:010} and \ref{alg:heur_mslp:020}, we initialize the best solution $BS$ and its cost $BS.cost$ as infinity to serve as a sentinel, as this is a minimization problem.
Then, in line \ref{alg:heur_mslp:030}, a greedy solution $S$ is constructed using an insertion heuristic.
If the solution $S$ is feasible~(line~\ref{alg:heur_mslp:040}), we run an LP-based improvement procedure to optimize the solution's schedule and update the best solution $BS$~(line~\ref{alg:heur_mslp:050}).
In the while loop of lines \ref{alg:heur_mslp:060}-\ref{alg:heur_mslp:110}, while the stopping criterion is not satisfied~(line~\ref{alg:heur_mslp:060}), we iteratively run the semi-greedy insertion heuristic~(line~\ref{alg:heur_mslp:070}) to construct a solution~$S$.
If \mbox{$S.feasible = \textbf{true}$}~(line~\ref{alg:heur_mslp:080}), we run the LP-based improvement procedure on $S$ to optimize the schedule~(line~\ref{alg:heur_mslp:090}).
Then, if $S.cost$ is smaller than $BS.cost$~(line~\ref{alg:heur_mslp:100}), we update $BS$ with $S$~(line~\ref{alg:heur_mslp:110}).
After the stopping criterion is satisfied, the algorithm returns the best solution $BS$~(line~\ref{alg:heur_mslp:120}).

\begin{algorithm}[!ht]
    \scriptsize
    \caption{MULTI\_START\_HEURISTIC\_LP\_IMPROVEMENT($I$,$\alpha$)}
    \label{alg:heur_mslp}
    $BS \gets \emptyset$\; \label{alg:heur_mslp:010}
    $BS.cost \gets \infty$\; \label{alg:heur_mslp:020}
    $S \gets$ GREEDY\_INSERTION\_HEURISTIC($I$)\; \label{alg:heur_mslp:030}
    \If{$S.feasible$}{ \label{alg:heur_mslp:040}
        $BS \gets$ LP\_SCHEDULE($S$, $I$)\; \label{alg:heur_mslp:050}
    }
    \While{stopping criterion not satisfied}{ \label{alg:heur_mslp:060}
        $S \gets$ SEMI\_GREEDY\_INSERTION\_HEURISTIC($I$, $\alpha$)\; \label{alg:heur_mslp:070}
        \If{$S.feasible$}{ \label{alg:heur_mslp:080}
            $S \gets$ LP\_SCHEDULE($S$, $I$)\; \label{alg:heur_mslp:090}
            \If{$S.cost < BS.cost$}{ \label{alg:heur_mslp:100}
                $BS \gets S$\; \label{alg:heur_mslp:110}
            }
        }
    }
    \Return $BS$\; \label{alg:heur_mslp:120}
\end{algorithm}

\subsection{Greedy insertion heuristic}
\label{sec:greedy_insertion_heuristic}

Algorithm \ref{alg:ins_heur} describes the pseudocode for the GREEDY\_INSERTION\_HEURISTIC called in line~\ref{alg:heur_mslp:030} of Algorithm \ref{alg:heur_mslp}.
The algorithm receives as input the instance data $I$.
In line~\ref{alg:ins_heur:010}, an empty solution~$S$ is initialized with $feasible=\textbf{true}$, which indicates a partially feasible solution, and no request is inserted.
Each vehicle route in $S.vehicles$ begins with two $visit$ variables, starting and ending at the depot.
Besides that, each sequence of travels in $S.machines$ starts at its respective initial station.
We consider that each vehicle leaves the depot at time $e_0 = 0$, i.e., $t^k_0 = 0, \ \forall \ k \in K$.
After obtaining a feasible solution, the LP method will adjust the depot departure and arrival times for each vehicle to reduce any unnecessary waiting time.
In line~\ref{alg:ins_heur:020}, an insertion sequence $\mathcal{L}$ of pickup nodes is defined as the list of pickup nodes sorted in nondecreasing order of time window width. 
Then, in the for loop of lines \ref{alg:ins_heur:030}-\ref{alg:ins_heur:090}, we will iteratively insert in order each pickup and its respective delivery node of $\mathcal{L}$ in one of the $|K|$ vehicle routes.
Thus, for each $p \in \mathcal{L}$ (line~\ref{alg:ins_heur:030}) and delivery $d$ (line~\ref{alg:ins_heur:040}), the algorithm searches for the best insertion candidate for the request nodes $p$ and $d$ in the partially constructed solution~$S$ (line~\ref{alg:ins_heur:050}).
If the candidate $\mathcal{C}$ is feasible (line~\ref{alg:ins_heur:060}), the solution~$S$ is updated, inserting $p$ and $d$ following the candidate $\mathcal{C}$ structure (line~\ref{alg:ins_heur:070}).
Otherwise, the insertion fails, and the partially constructed solution~$S$, now infeasible (line~\ref{alg:ins_heur:085}), is returned (line~\ref{alg:ins_heur:090}).
After inserting all the requests, $S$ is returned (line~\ref{alg:ins_heur:110}).

\begin{algorithm}[H]
\scriptsize
\caption{\footnotesize GREEDY\_INSERTION\_HEURISTIC($I$)}
\label{alg:ins_heur}
    $S \gets $ Empty solution~$S$ with $S.feasible=\textbf{true}$ and no request inserted\; \label{alg:ins_heur:010}
    $\mathcal{L} \gets$ List of pickup nodes sorted in nondecreasing order of time window width\; \label{alg:ins_heur:020}
    \For{$p \in \mathcal{L}$\label{alg:ins_heur:030}}{ 
        $d \gets p+n$\; \label{alg:ins_heur:040}
        $\mathcal{C} \gets$ GET\_BEST\_INSERTION\_CANDIDATE($p$, $d$, $S$, $I$)\; \label{alg:ins_heur:050}
        \eIf{$\mathcal{C}.feasible$}{ \label{alg:ins_heur:060}
            $S \gets$ UPDATE\_SOLUTION($S$, $p$, $d$, $\mathcal{C}$, $I$)\; \label{alg:ins_heur:070}
        }{\label{alg:ins_heur:080}
            $S.feasible \gets \textbf{false}$\; \label{alg:ins_heur:085}
            \Return $S$\; \label{alg:ins_heur:090}
        }
    }
 \Return $S$\; \label{alg:ins_heur:110}
\end{algorithm}

Algorithm \ref{alg:best_ins_cand} describes the pseudocode for the GET\_BEST\_INSERTION\_CANDIDATE called in line~\ref{alg:ins_heur:050} of Algorithm \ref{alg:ins_heur}, which searches for the best insertion candidate for the nodes $p$ and $d$ in $S$.
In line~\ref{alg:best_ins_cand:010}, we initiate the best candidate $B\mathcal{C}$ with an infeasible dummy insertion candidate $\mathcal{C}$ with $\mathcal{C}.cost = \infty$ to serve as a sentinel.
Then, in lines \ref{alg:best_ins_cand:030}-\ref{alg:best_ins_cand:060}, we analyze the insertion of $p$ and $d$ in positions $p_{pos}$ and $d_{pos}$ of the route $S.vehicles[k]$ for vehicle $k \in K$, respectively.
If the insertion candidate $\mathcal{C}$ is feasible and its cost $\mathcal{C}.cost$ is smaller than $B\mathcal{C}.cost$ (line     \ref{alg:best_ins_cand:060}), we update $B\mathcal{C}$ by setting candidate $\mathcal{C}$.
Note that the evaluation order determines $B\mathcal{C}$ to be the first $\mathcal{C}$ evaluated with $\mathcal{C}.cost = B\mathcal{C}.cost$.
The evaluation order follows the given vehicle order defined in the data set, in addition to the values of $p_{pos}$ and $d_{pos}$ indicated in line~\ref{alg:best_ins_cand:030}.
After analyzing all candidates, we return $B\mathcal{C}$ (line~\ref{alg:best_ins_cand:070}).

\begin{algorithm}[H]
    \scriptsize
    \caption{\footnotesize GET\_BEST\_INSERTION\_CANDIDATE($p$, $d$, $S$, $I$)}
    \label{alg:best_ins_cand}
    $B\mathcal{C} \gets$ Dummy insertion candidate $\mathcal{C}$ with $\mathcal{C}.cost = \infty$ and $\mathcal{C}.feasible = \textbf{false}$\; \label{alg:best_ins_cand:010}
    \For{$k \in K$ \textbf{and} $p_{pos} \in 2:S.vehicles[k].length$ \textbf{and} $d_{pos} \in p_{pos}:S.vehicles[k].length$}{ \label{alg:best_ins_cand:030}
        $\mathcal{C} \gets$ ANALYZE\_INSERTION\_FEASIBILITY($k$, $p_{pos}$, $d_{pos}$, $p$, $d$, $I$)\; \label{alg:best_ins_cand:040}
        \If{$\mathcal{C}.feasible$ \textbf{and} $\mathcal{C}.cost < B\mathcal{C}.cost$}{ \label{alg:best_ins_cand:050}
            $B\mathcal{C} \gets \mathcal{C}$\; \label{alg:best_ins_cand:060}
        }
    }
    \Return $B\mathcal{C}$\; \label{alg:best_ins_cand:070}
\end{algorithm}

\subsection{Semi-greedy insertion heuristic}
\label{sec:semi-greedy-ins-heur}

Algorithm \ref{alg:semi_ins_heur} presents the pseudocode for the SEMI\_GREEDY\_INSERTION\_HEURISTIC called in line~\ref{alg:heur_mslp:070} of Algorithm~\ref{alg:heur_mslp}, which randomizes the Algorithm~\ref{alg:ins_heur}.
Thereafter, we highlight the differences between them.
In line \ref{alg:semi_ins_heur:020}, the insertion order $\mathcal{L}$ is now defined as a random sequence of pickup tasks.
In line \ref{alg:semi_ins_heur:050}, the algorithm constructs a candidate list $CL$ containing all feasible insertion candidates for request nodes $p$ and $d$ in solution $S$.
In line \ref{alg:semi_ins_heur:060}, a candidate is chosen from $CL$ by randomly selecting a candidate $\mathcal{C}$ from $CL$ restricted by quality, which requires the threshold parameter $\alpha$ (see Algorithm~\ref{alg:cand_by_qlty}).

\begin{algorithm}[H] 
\scriptsize
\caption{\footnotesize SEMI\_GREEDY\_INSERTION\_HEURISTIC($I$, $\alpha$)}\label{alg:semi_ins_heur}
    $S \gets $ Empty solution $S$ with $S.feasible=\textbf{true}$ and no request inserted\;
    $\mathcal{L} \gets$ Random sequence of pickup tasks\; \label{alg:semi_ins_heur:020}
    \For{$p \in \mathcal{L}$}{
        $d \gets p+n$\;
        $CL \gets$ GET\_INSERTION\_CANDIDATE\_LIST($p$, $d$, $S$, $I$)\; \label{alg:semi_ins_heur:050}
        $\mathcal{C} \gets$ CHOOSE\_CANDIDATE\_BY\_QUALITY($CL$, $\alpha$)\; \label{alg:semi_ins_heur:060}
        \eIf{$\mathcal{C}.feasible$}{
            $S \gets$ UPDATE\_SOLUTION($S$, $p$, $d$, $\mathcal{C}$, $I$)\;
        }{
            $S.feasible \gets \textbf{false}$\; 
            \Return $S$\;
        }
    }
    \Return $S$\;
\end{algorithm}

The algorithms for the GET\_INSERTION\_CANDIDATE\_LIST called In line \ref{alg:semi_ins_heur:050} of Algorithm~\ref{alg:semi_ins_heur} and the CHOOSE\_CANDIDATE\_BY\_QUALITY invoked In line \ref{alg:semi_ins_heur:060} of Algorithm~\ref{alg:semi_ins_heur} are given in \mbox{\ref{app:heur_details}}.

\FloatBarrier
\subsection{Linear programming based improvement procedure}
\label{sec:lp_form}

In this section, we present the linear programming (LP) formulation employed in the LP\_{\linebreak}SCHEDULE method invoked in lines \ref{alg:heur_mslp:050} and \ref{alg:heur_mslp:090} of Algorithm~\ref{alg:heur_mslp}. Given as input a solution represented by the vehicles' routes and the sequences of the machine traversals, it attempts to reduce the total completion time.
The LP formulation adapts the scheduling part of the MIP formulation presented in Section \ref{sec:formulation} by assuming that all the decisions associated with integer variables are predefined.
In what follows, Section \ref{subsec:lp_params} introduces the notation used for the input parameters.
Section \ref{subsec:lp_sched_const} details the LP formulation.

\subsubsection{LP parameters}
\label{subsec:lp_params}

Given a solution $S$, denote by $L^K_k = \{1,\ldots, \ell^K_k\}$ the route positions for each $visit \in S.vehicles[k]$, $k \in K$, $visit \in V_p \cup V_d$, where $\ell^K_k$ represents the number of these visits in the route.
Moreover, let \mbox{$\sigma^k = (\sigma^k_1, \sigma^k_2, \ldots,\sigma^k_{\ell^K_k})$} be the sequence of nodes $i \in V_p \cup V_d$ in $S.vehicles[k], k \in K$, i.e., $\forall \ l \in L^K_k, \ \sigma^k_l = i, \ i \in V_p \cup V_d$.

Furthermore, denote by $L^H_h = \{1,\ldots, \ell^H_h\}$ the sequence of travels of $S.machines[h], \ h \in H$, with $\ell^H_h$ denoting the sequence length.
Besides, let $\psi^h = (\psi^h_1, \psi^h_2, \ldots,\psi^h_{\ell^H_h})$ be the sequence of machine travels $(i,j,k), \ (i,j) \in A^m, \ k \in K$ in the sequence, i.e., $\forall \ l \in L^H_h, \ travel = S.machines[h][l], \ \psi^h_l = (travel.i,travel.j,travel.k)$, $(travel.i,travel.j) \in A^m, \ travel.k \in K$.

\subsubsection{The LP formulation}
\label{subsec:lp_sched_const}

The formulation uses the following continuous variables already defined in Section~\ref{sec:mipscheduling}: $t_i$ for $i\in V_p \cup V_d$; $t^k_0$, $t^k_{2n+1}$, and $C_k$ for $k\in K$; and 
$\alpha^h_{ij}$ for $(i,j)\in A^m$, $h\in H_{ij}$.
However, to reduce the number of variables, we only define $\alpha^h_{ij}$ when $\exists \ l \in L^H_h \ : \ \psi^h_{l} = (i,j,k), \ (i,j) \in A^m, \ h \in H$.
In this way, the LP formulation for the scheduling problem to minimize the total completion time with predefined sequences for the routes and machine travels is given by:
\begingroup
\small
\allowdisplaybreaks
\begin{align}
    & \min \ \sum_{k \in K} C_k. \label{lpfeas:00} \\
    & t_{\sigma^k_i} \geq t_{\sigma^k_{i-1}} + s_{\sigma^k_{i-1}} + d^k_{\sigma^k_{i-1}, \sigma^k_{i}}, \quad \forall \ k \in K, \ \ell^K_k > 0, \ i \in L^K_k \setminus \{1\}, \ (\sigma^k_{i-1}, \sigma^k_{i}) \in A^s.    \label{lpfeas:01} \\
    & t_{\sigma^k_1} \geq t^k_{0} + d^k_{0 \sigma^k_{1}}, \quad \forall \ k \in K, \ \ell^K_k > 0, \ (0, \sigma^k_{1}) \in A^s.    \label{lpfeas:02} \\
    & \alpha^h_{i j} \geq t_{i} + s_{i} + \bar{d}^{k}_{ih}, \quad \forall \ h \in H, \ \ell^H_h > 0, \ l \in L^H_h, \ \psi^h_l = (i,j,k), \ i \neq 0.   \label{lpfeas:03} \\
    & \alpha^h_{0 j} \geq t^{k}_{0} + \bar{d}^{k}_{0 h}, \quad \forall \ h \in H, \ \ell^H_h > 0, \ l \in L^H_h, \ \psi^h_l = (i,j,k) , \ i = 0.     \label{lpfeas:04} \\
    & t_{j} \geq \alpha^h_{i j} + O^h_{f^h_{i}f^h_{j}} + \bar{d}^{k}_{jh}, \quad \forall \ h \in H, \ \ell^H_h > 0, \ l \in L^H_h , \ \psi^h_l = (i,j,k), \ j \neq 2n+1.    \label{lpfeas:05} \\
    & \alpha^h_{i j} \geq \alpha^h_{i' j'} + O^h_{f^h_{i'} f^h_{j'}} + O^h_{f^h_{j'} f^h_{i}}, \quad \forall \ h \in H, \ \ell^H_h > 0, \ l \in L^H_h \setminus \{1\}, \ \psi^h_l = (i,j,k), \ \psi^h_{l-1} = (i',j',k'). \label{lpfeas:06}\\
    & \alpha^h_{i j} \geq O^h_{0 f^h_{i}}, \quad \forall \ h \in H, \ \ell^H_h > 0, \ \psi^h_1 = (i,j,k). \label{lpfeas:07} \\
    & t^k_{2n+1} \geq t_{\sigma^k_{\ell^K_k}} + s_{\sigma^k_{\ell^K_k}} + d^k_{\sigma^k_{\ell^K_k}, 2n+1}, \quad \forall \ k \in K, \ \ell^K_k > 0, \ (\sigma^k_{\ell^K_k}, 2n+1) \in A^s. \label{lpfeas:08} \\
    & t^{k}_{2n+1} \geq \alpha^h_{i, 2n+1} + O^h_{f^h_{i} f^h_{2n+1}} + \bar{d}^{k}_{2n+1,h}, \quad \forall \ h \in H, \ \ell^H_h > 0, \ l \in L^H_h, \ \psi^h_l = (i,j,k), \ j = 2n+1. \label{lpfeas:09} \\
    &  C_k \geq t^k_{2n+1} - t^k_0 , \quad \forall \ k \in K.    \label{lpfeas:10} \\
    &  e_i \leq t_i \leq l_i ,  \quad \forall \ k \in K.   \label{lpfeas:11} \\
    &  e_0 \leq t^k_0 \leq t^k_{2n+1} \leq l_0 ,  \quad \forall \ k \in K.    \label{lpfeas:12} \\
    & C_k \geq 0, \quad \forall \ k \in K. \label{lpfeas:13} \\
    & \alpha^h_{i j} \geq 0, \quad \forall \ h \in H, \ \ell^H_h > 0, \ l \in L^H_h, \ \psi^h_l = (i,j,k). \label{lpfeas:14} 
\end{align}
\endgroup

The objective function \eqref{lpfeas:00} is the same for the MIP formulation, i.e., minimizing the total completion time of the vehicles.
Constraints~\eqref{lpfeas:01}-\eqref{lpfeas:02} ensure that whenever an intra-region arc is traversed by a vehicle, the time to serve the destination node is determined by the times of its origin.
Constraints~\eqref{lpfeas:03}-\eqref{lpfeas:04} define that whenever a machine is used to traverse an arc, the starting time on the machine is determined by the times of its origin node.
Constraints~\eqref{lpfeas:05} ensure that whenever an inter-region arc is traversed by the vehicle, the time to serve the destination node is determined by the times of its machine travel origin station.
Constraints~\eqref{lpfeas:06} define that the machine can only start a machine travel after finishing the previous travel and traveling to the next machine travel origin station. 
Constraints~\eqref{lpfeas:07} ensure that the machine can only start the first machine travel after traveling from its initial station to the first machine travel origin station.
Constraints~\eqref{lpfeas:08}-~\eqref{lpfeas:09} determine the arrival times of the vehicles at the depot based on their origin nodes.
Constraints~\eqref{lpfeas:10} set the lower bounds on the completion times of the vehicles.
Constraints~\eqref{lpfeas:11}-\eqref{lpfeas:14} restrain the domains of the $\boldsymbol{t}$, $\boldsymbol{t_0}$, $\boldsymbol{t_{2n+1}}$, $\boldsymbol{C}$, and $\boldsymbol{\alpha}$ variables.

\section{Computational experiments}
\label{sec:experiments}

This section summarizes the computational experiments carried out to evaluate the performance of both the proposed MIP formulation and the MSLP heuristic and to analyze the characteristics of the solutions to the problem. The tests were run on a machine with an Intel(R) Core(TM) \mbox{i7-10700} processor CPU @ 2.90GHz, 16 GB of RAM, 16 logical processors (eight physical cores), and Ubuntu~22.04 operating system. The solver Gurobi version 12.0.2 was applied to solve the MIP and LP formulations using the default settings, except for the number of threads. The MIP formulation was implemented in Python 3.10.12, using the package gurobipy v12.0.2.
We also applied Hexaly version~14.0 to solve the MIP formulation, which was implemented in Python 3.10.12 using Hexaly's Python package v14.0.
The MSLP heuristic was implemented in Julia 1.11.3, using the following packages: JuMP v1.23.6, Gurobi v1.5.0, and Gurobi\_jll v12.0.2.
We remark that implementing the formulation for Gurobi and Hexaly in Python does not compromise the validity of the results, as the solvers run their internal algorithms using the same underlying mathematical formulation regardless of the programming interface.
For the MIP formulation, in both solvers, we imposed a time limit of 3600 seconds (1 hour) on each solver execution, and also limited the number of threads to eight.
For the MSLP heuristic, we set the stop criterion to be 60000 iterations or 3600 seconds on each independent execution, and also limited the number of threads to one.
Such an iteration limit is reasonable for small instances with up to 12 requests because it generally takes less than 600 seconds~to run.

\okC{
Section~\ref{subsec:inst_gen} presents the benchmark instances.
Section~\ref{subsec:mipresults_no_vis} discusses the results obtained by the two solvers on the \okC{baseline MIP formulation, which includes the preprocessing step (see Section~\ref{sec:mippreprocessing}), but not the valid inequalities (see Section~\ref{sec:mip_tightening})}.
Section~\ref{subsec:impact_valid_inequalities} analyzes the impact of the valid inequalities on the MIP formulation with preprocessing.
Section~\ref{subsec:heur_results} compares the results for the MIP formulation and the MSLP heuristic, respectively.
\okC{
Section~\ref{subsec:sol_chars_summary} summarizes the analysis of the solutions' characteristics discussed in \ref{app:solution_chars}.
}
}
\subsection{Benchmark instances}
\label{subsec:inst_gen}

We created a benchmark set of instances for the PDPTW-SE based on the randomized instances of \cite{LiLim03} for the PDPTW \citep{Sintefpdptw}. 
The latter correspond to a traditional benchmark set in the vehicle routing literature~\citep{Gunawan2021} and were obtained from the VRPTW instances~\citep{solomon87} by randomly pairing the customer locations within the routes of the solutions achieved by \cite{LiLim03}'s heuristic approach. 
Our instance generators were implemented in Python 3.9.24.

We generated two families of instances, denoted \textit{multi-island} and \textit{multi-floor}, that represent two different potential applications of the PDPTW-SE.
These differ in how the nodes and machine stations are placed into the environment and how the environment is partitioned into regions.
In what follows, Sections \ref{subsubsec:multi-isl} and \ref{subsubsec:multi-floor} describe the placement of nodes in regions and the generation of machines for the multi-island and multi-floor instances, respectively.
After that, we describe the information shared between the multi-island and multi-floor instances. 
Namely,
Section \ref{subsubsec:vehicle-gen} explains the process of generating the heterogeneous fleet of vehicles.
Section \ref{subsubsec:inst-feas} deals with how we guarantee instance feasibility.
Section \ref{subsubsec:inst-size} describes the criteria for selecting instance sizes.
Section \ref{subsubsec:inst-chars} presents the chosen instance characteristics.

\subsubsection{Multi-island instances}
\label{subsubsec:multi-isl}
    
\begin{figure}[!ht]
\centering
\subfigure[Multi-island instance]{
\includegraphics[scale=0.42]{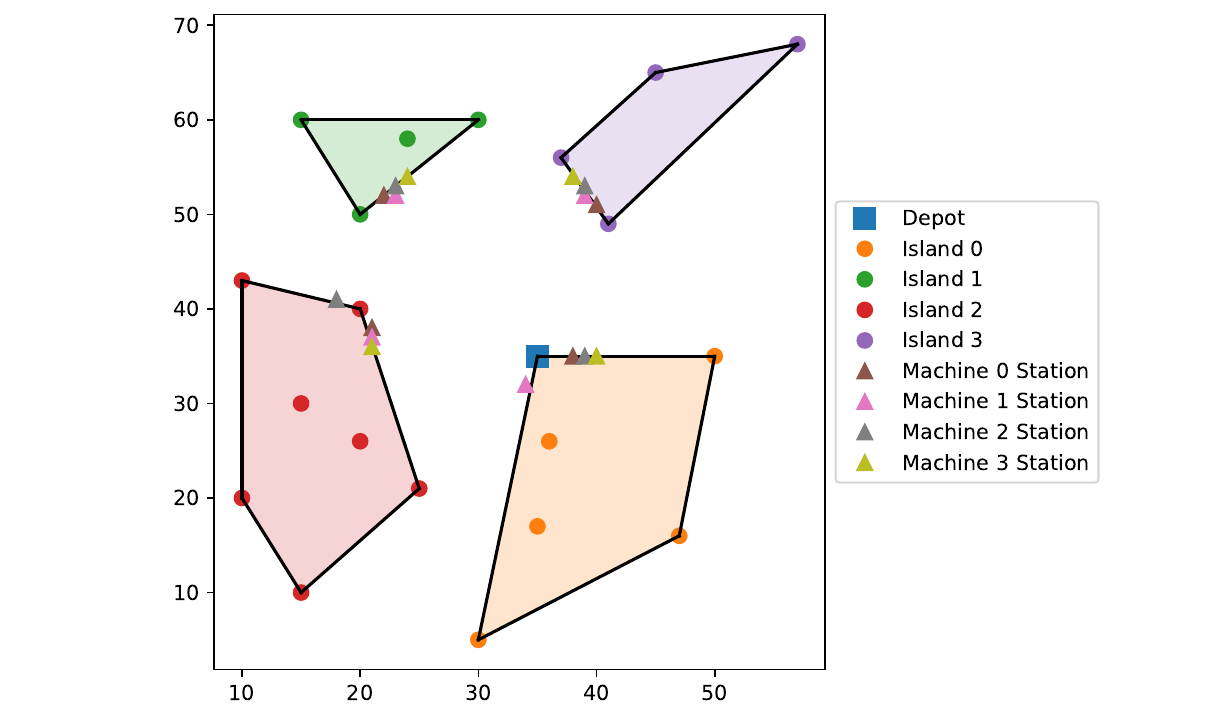}
\label{fig:ex_map_inst_isl}
}
\hspace{-1.5cm}
\subfigure[Multi-floor instance]{
\includegraphics[scale=.42]{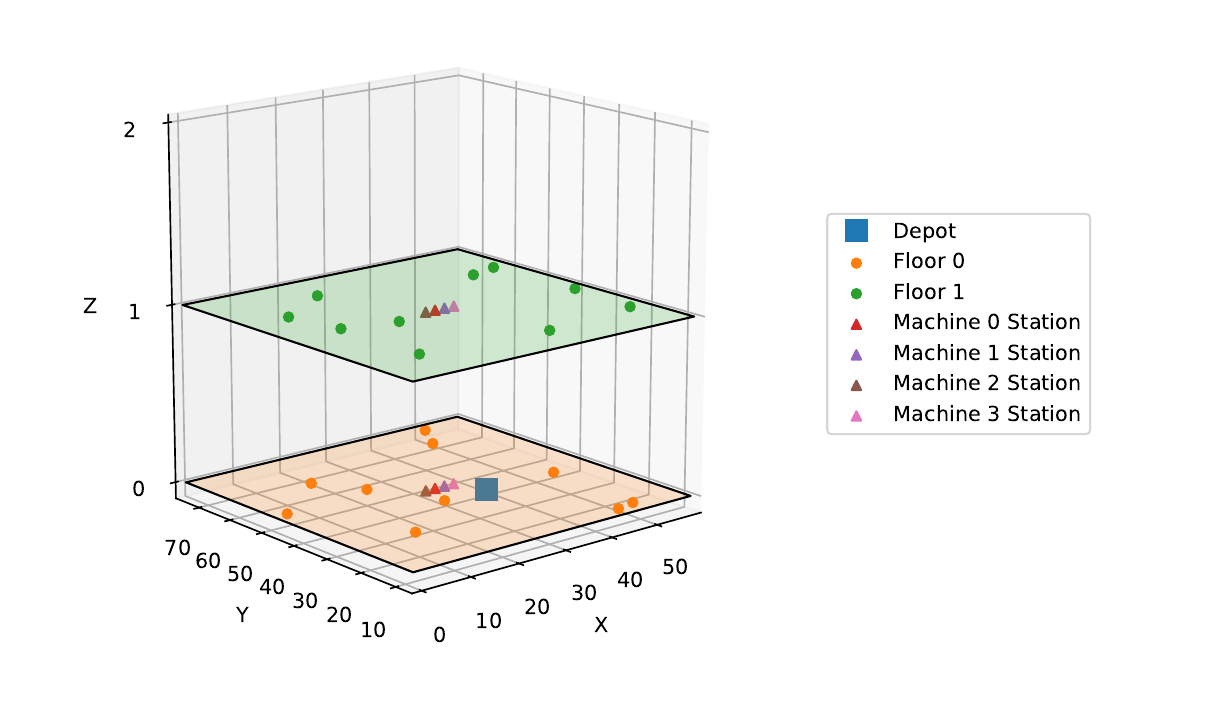}
\label{fig:ex_map_inst_isl_floor}
}
\caption{Illustration map for the two instance families}
\end{figure}

In this family, the nodes are subdivided into two-dimensional regions, denoted islands.
This subdivision uses the k-means clustering technique implemented in Python's Scikit-learn library v1.2.0~\citep{scikit-learn}.
The default parameters were used, except for the number of clusters, which were defined based on the number of islands required.
The machine stations are generated as follows.
We first find the convex hull for each island set of points using Python's SciPy library v1.9.3~\citep{2020SciPy-NMeth}. 
The convex hull serves as an approximation of the island's shape.
If there are less than three points or the points are collinear, we continually add a dummy point around the region's centroid until the algorithm can find a convex hull.
After generating the convex hulls, we apply a backtracking algorithm to find the set of vertices, one for each convex hull, that minimizes the total distance between any two points in this set, i.e., the minimum weighted clique.
As the number of requests is not \okC{too} large \okC{(see Section \ref{subsubsec:inst-size})}, this can be done in a reasonable computational time.
Then, for each machine, we generate the machine points by selecting for each vertex in the resulting backtracking set the best point in one of the two incident edges. 
The criteria for the best point is the minimum sum of the distances from the other chosen vertices.
This point is at least two units of Euclidean distance apart from the vertex and should be distinct from any other existing point.
To generate the point on the incident edge, we start closer to the chosen vertex and gradually increase the distance from it using the linear combination of the edge vertices' points until we find an available point.
These points are also rounded to obtain integer points. 
We define the machine velocity in this family as one unit of distance per unit of time. Figure \ref{fig:ex_map_inst_isl} illustrates a multi-island instance with ten requests, four islands, and four machines.

\subsubsection{Multi-floor instances}
\label{subsubsec:multi-floor}

In this family, the requests are subdivided into three-dimensional regions, denoted floors. 
Such a subdivision is made by randomly selecting a floor in $F$ for each task, based on the number of floors in the instance, setting the floor number as the $z$-axis value for each point. 
For example, requests in an instance with four floors will have $z$-axis values in $F = \{0,1,2,3\}$.
The distance between each floor is the Euclidean distance.

The machine stations are placed as a stack, i.e., they all have the same points for $x$ and $y$, but each floor ($z$) has one station.
The machines are centralized based on the ranges of point values for $x$ and $y$.
Starting at the center point, we alternately test the possibility of placing a point on the right or the left. 
If it is not possible, i.e., this point already contains a node or station, we continue to attempt to increase the distance from the center until a point is available.
The machine velocity in this family is $0.2$ distance unit per unit of time.
Figure \ref{fig:ex_map_inst_isl_floor} illustrates a multi-floor instance with ten requests, two floors, and three machines.

\subsubsection{Vehicle generation}
\label{subsubsec:vehicle-gen}

Given that the fleet is homogeneous in the PDPTW, the vehicles were generated by changing their capacities. 
As the number of requests is reduced, the original capacity becomes much larger than any possible load value on any possible route configuration. 
The PDPTW instances of \cite{LiLim03} didn't change the original vehicle capacity, even though the maximum load in a PDPTW solution is much smaller than that in a VRPTW solution.
The reason is that pairing the tasks in the route implies that the delivery task demand becomes the opposite of the respective pickup task demand.
Therefore, the vehicle capacity becomes irrelevant in obtaining the best solution for the PDPTW.
Then, we first define a base capacity $\mathcal{B}_Q$ that will be the capacity of the mid-size vehicle. 
The base capacity is inspired by \cite{Sartori2020}'s PDPTW instances. 
They generate the demands of the requests selecting uniformly from the interval $[10, 0.6 \cdot Q]$, where $Q$ is their capacity value. 
Similarly, we define $\mathcal{B}_Q$ as $\texttt{max}_{i \in V_p}(q_i)/0.6$.
Now, we set the capacity of vehicle $k$ as $Q_k=\mathcal{B}_Q + \texttt{random}(-\floor*{t/2}, \ceil*{t/2})\cdot \texttt{round}(\Delta{q}\cdot \mathcal{B}_Q)$, where $\texttt{random}(a,b)$ returns a random integer from discrete uniform distribution in the half-open interval $[a,b)$, $t$ is the number of different vehicle types, $\texttt{round}(a)$ rounds the number $a$ to the closest integer and $\Delta q$ is the capacity variation, in percentage.

\subsubsection{Instance feasibility}
\label{subsubsec:inst-feas}

Note that the PDPTW-SE becomes more complex than PDPTW since it encompasses obstacles that can only be traversed with machines. 
Consequently, not all adapted instances are naturally feasible. 
\okC{
Therefore, a procedure that can eventually modify the instance to ensure feasibility was devised.
This procedure, described in~\ref{app:alg_to_ensure_feas}, is an adaptation of the greedy constructive heuristic described in Section \ref{sec:greedy_insertion_heuristic}.
}
The modifications cover two aspects: shifting the time window and increasing vehicle capacity.
\okC{An overview of the instance changes is discussed in more detail in \ref{app:overview_inst_changes}.}

\subsubsection{Instance size}
\label{subsubsec:inst-size}

MIP formulations can generally only solve to optimality smaller instances within a reasonable amount of time for problems similar to the PDPTW-SE.
Besides, the original PDPTW instances from \cite{LiLim03} have at least 50 requests (100 pickup/delivery nodes).
For this reason, we consider \okC{a first set of} instances with up to 12 requests (24 pickup/delivery nodes)\okC{, which were generated from the original PDPTW instances with 50 requests (100 pickup/delivery nodes)}.
\okC{From this number onwards, finding even an initial feasible solution for \okC{an instance} becomes almost prohibitive for the solvers used in this work.
However, the heuristic can find feasible solutions for larger instances.
Therefore, we also developed a second \okC{set of larger instances} with up to 60 requests (120 pickup/delivery nodes) only for the heuristic's performance evaluation, which were generated from the original PDPTW instances with 100 requests (200 pickup/delivery nodes).} The chosen criterion used to reduce the number of requests was to select the first $n$ pickup tasks with their respective delivery tasks.
We also remove the dummy tasks created for the PDPTW instances in \cite{LiLim03}, as their goal was only to pair unpaired tasks in a vehicle route with an odd number of tasks.
\cite{solomon87}'s VRPTW 100-customer instances can be reduced to their first 25, 50, and 75 customers to obtain smaller instances, but this doesn't apply to \cite{LiLim03}'s PDPTW instances because they didn't create adapted instances for these smaller instances.

\subsubsection{Characteristics of the instances}
\label{subsubsec:inst-chars}

The characteristics presented herein are shared by the multi-island and multi-floor instances.
The randomized instances of \citet{LiLim03} are classified into two types.
Type 1 instances have short scheduling horizons, while Type 2 instances have long ones.
\okC{The first set comprises the first five instances of each type: lr101-lr105 (Type 1) and lr201-lr205 (Type 2).
We vary the number of requests ($n \in \{6,8,10,12\}$), vehicles ($|K| = n$), regions ($z \in \{2,4\}$), and machines ($|H| \in \{3,4\}$). 
The second set encompasses the first five instances of each type: LR1\_2\_1-LR1\_2\_5 (Type 1) and LR2\_2\_1-LR2\_2\_5 (Type 2).
We vary the number of requests ($n \in \{40, 60\}$), vehicles ($|K| = n$), regions ($z \in \{2,4\}$), and machines ($|H| \in \{5,6\}$).
}

\okC{All instances were initially generated such that the first and second sets contained four and six machines, respectively.} 
However, we ensure feasibility (see Section \ref{subsubsec:inst-feas}) considering only the first three or five machines in the first and second sets, respectively, according to the data set instance.
Thus, the last machine of the data set instance is an optional additional machine.
Furthermore, we define three types of vehicles with variation $\Delta q = 20\%$, i.e., $Q_k = \mathcal{B}_Q + \texttt{random}(-1,2) \times \texttt{round}(0.2 \cdot \mathcal{B}_Q)$.
Therefore, there is at least one vehicle per request, and all vehicles can load the maximum weight request. 
This ensures instance feasibility regarding capacity constraints.
\okC{Besides, the first three vehicles of an instance are fixed to represent the three distinct types.}

For each combination, we adopt the following convention name: $n_rR\_n_vV\_n_iI\_n_mM$ (multi-island set) and $n_rR\_n_vV\_n_fF\_n_mM$ (multi-floor set), where $n_r=n$ is the number of requests, $n_v=|K|$ is the number of vehicles, $n_i=|F|$ is the number of islands, $n_f=|F|$ is the number of floors, $n_m=|H|$ is the number of machines.
For all the instances, we consider that any distance between two nodes in the same region $f \in F$ or between two machine stations is the Euclidean distance.
\okC{In total, we created 480 instances: 320 instances in the first set, and 160 in the second set.}

\subsection{\okC{Baseline MIP results}}
\label{subsec:mipresults_no_vis}

Tables \ref{tab:overview_results_multi_island} and \ref{tab:overview_results_multi_floor} present an overview of the results obtained by the solver \okC{Gurobi} using the proposed \okC{baseline} MIP formulation for the two instance families.
\okC{These experiments include the MIP preprocessing step described in Section \ref{sec:mippreprocessing}, but do not apply the valid inequalities proposed in Section \ref{sec:mip_tightening}.
Besides, we only considered the instances with up to 12 requests, given that the remaining instances are too large for the solver to find at least a feasible solution.}
In these tables, each line represents an instance group with five instances, which is identified in the first column.
The next columns indicate, for each instance type, the number of instances for which the solver: proved solution optimality (Opt.), exceeded the time limit (TLE), was not able to finish the execution because it was killed by the operating system (Killed), and obtained at least one feasible solution (Feas.).
Notice that columns Opt., TLE, and Killed are complementary, i.e., they sum up to five.

%\FloatBarrier
\begin{table}[!ht]
    \centering
    \caption{\okC{Overview of the baseline MIP results for the multi-island family.}}
    \label{tab:overview_results_multi_island}
    \scriptsize
    \okC{
    \begin{tabular}{l r rrrr r rrrr}
        \toprule
        \multicolumn{1}{c}{} &\hspace{0.1cm} &\multicolumn{4}{c}{Type~1} & \hspace{0.2cm} &\multicolumn{4}{c}{Type~2} \\
        \cmidrule(lr){3-6}
        \cmidrule(lr){8-11}
        Group & \hspace{0.1cm}  & Opt. & TLE & Killed & Feas. & \hspace{0.2cm} & Opt. & TLE & Killed & Feas.\\
        \cmidrule(lr){1-1}
        \cmidrule(lr){3-6}
        \cmidrule(lr){8-11}
        06R\_06V\_02I\_03M & \hspace{0.1cm} & 5 & 0 & 0 & 5 & \hspace{0.2cm} & 5 & 0 & 0 & 5\\
        06R\_06V\_02I\_04M & \hspace{0.1cm} & 5 & 0 & 0 & 5 & \hspace{0.2cm} & 5 & 0 & 0 & 5\\
        06R\_06V\_04I\_03M & \hspace{0.1cm} & 5 & 0 & 0 & 5 & \hspace{0.2cm} & 4 & 1 & 0 & 5\\
        06R\_06V\_04I\_04M & \hspace{0.1cm} & 5 & 0 & 0 & 5 & \hspace{0.2cm} & 4 & 1 & 0 & 5\\
        \addlinespace
        08R\_08V\_02I\_03M & \hspace{0.1cm} & 3 & 2 & 0 & 5 & \hspace{0.2cm} & 2 & 3 & 0 & 5\\
        08R\_08V\_02I\_04M & \hspace{0.1cm} & 3 & 2 & 0 & 5 & \hspace{0.2cm} & 2 & 3 & 0 & 5\\
        08R\_08V\_04I\_03M & \hspace{0.1cm} & 3 & 2 & 0 & 5 & \hspace{0.2cm} & 1 & 4 & 0 & 5\\
        08R\_08V\_04I\_04M & \hspace{0.1cm} & 3 & 2 & 0 & 5 & \hspace{0.2cm} & 1 & 4 & 0 & 5\\
        \addlinespace
        10R\_10V\_02I\_03M & \hspace{0.1cm} & 3 & 2 & 0 & 5 & \hspace{0.2cm} & 0 & 5 & 0 & 5\\
        10R\_10V\_02I\_04M & \hspace{0.1cm} & 3 & 2 & 0 & 5 & \hspace{0.2cm} & 0 & 5 & 0 & 5\\
        10R\_10V\_04I\_03M & \hspace{0.1cm} & 2 & 3 & 0 & 5 & \hspace{0.2cm} & 0 & 4 & 1 & 4\\
        10R\_10V\_04I\_04M & \hspace{0.1cm} & 3 & 2 & 0 & 5 & \hspace{0.2cm} & 0 & 3 & 2 & 3\\
        \addlinespace
        12R\_12V\_02I\_03M & \hspace{0.1cm} & 2 & 3 & 0 & 5 & \hspace{0.2cm} & 0 & 4 & 1 & 4\\
        12R\_12V\_02I\_04M & \hspace{0.1cm} & 2 & 3 & 0 & 5 & \hspace{0.2cm} & 0 & 5 & 0 & 4\\
        12R\_12V\_04I\_03M & \hspace{0.1cm} & 1 & 4 & 0 & 5 & \hspace{0.2cm} & 0 & 3 & 2 & 2\\
        12R\_12V\_04I\_04M & \hspace{0.1cm} & 1 & 4 & 0 & 5 & \hspace{0.2cm} & 0 & 4 & 1 & 4\\
        \midrule
        Sum & \hspace{0.1cm} & 49 & 31 & 0 & 80 & \hspace{0.2cm} & 24 & 49 & 7 & 71\\
        \bottomrule
    \end{tabular}
    \begin{tablenotes}
      \scriptsize
      \item Instances killed:  lr204 from 10R\_10V\_04I\_03M; lr203 and lr204 from 10R\_10V\_04I\_04M; lr204 from 12R\_12V\_02I\_03M; lr202 and lr203 from 12R\_12V\_04I\_03M; lr203 from 12R\_12V\_04I\_04M
    \end{tablenotes}
    }
\end{table}
%\FloatBarrier

\begin{table}[!ht]
    \centering
    \caption{\okC{Overview of the baseline MIP results for the multi-floor family.}}
    \label{tab:overview_results_multi_floor}
    \scriptsize
    \okC{
    \begin{tabular}{l r rrrr r rrrr}
        \toprule
        \multicolumn{1}{c}{} &\hspace{0.1cm} &\multicolumn{4}{c}{Type~1} &\hspace{0.2cm} &\multicolumn{4}{c}{Type~2} \\
        \cmidrule(lr){3-6}
        \cmidrule(lr){8-11}
        Group &\hspace{0.1cm} & Opt. & TLE & Killed & Feas. &\hspace{0.2cm} & Opt. & TLE & Killed & Feas.\\
        \cmidrule(lr){1-1}
        \cmidrule(lr){3-6}
        \cmidrule(lr){8-11}
        06R\_06V\_02F\_03M & \hspace{0.1cm} & 5 & 0 & 0 & 5 & \hspace{0.2cm} & 5 & 0 & 0 & 5\\
        06R\_06V\_02F\_04M & \hspace{0.1cm} & 5 & 0 & 0 & 5 & \hspace{0.2cm} & 4 & 1 & 0 & 5\\
        06R\_06V\_04F\_03M & \hspace{0.1cm} & 5 & 0 & 0 & 5 & \hspace{0.2cm} & 4 & 1 & 0 & 5\\
        06R\_06V\_04F\_04M & \hspace{0.1cm} & 5 & 0 & 0 & 5 & \hspace{0.2cm} & 4 & 1 & 0 & 5\\
        \addlinespace
        08R\_08V\_02F\_03M & \hspace{0.1cm} & 5 & 0 & 0 & 5 & \hspace{0.2cm} & 1 & 4 & 0 & 5\\
        08R\_08V\_02F\_04M & \hspace{0.1cm} & 5 & 0 & 0 & 5 & \hspace{0.2cm} & 1 & 4 & 0 & 5\\
        08R\_08V\_04F\_03M & \hspace{0.1cm} & 5 & 0 & 0 & 5 & \hspace{0.2cm} & 1 & 4 & 0 & 5\\
        08R\_08V\_04F\_04M & \hspace{0.1cm} & 5 & 0 & 0 & 5 & \hspace{0.2cm} & 1 & 4 & 0 & 5\\
        \addlinespace
        10R\_10V\_02F\_03M & \hspace{0.1cm} & 2 & 3 & 0 & 5 & \hspace{0.2cm} & 0 & 5 & 0 & 5\\
        10R\_10V\_02F\_04M & \hspace{0.1cm} & 3 & 2 & 0 & 5 & \hspace{0.2cm} & 1 & 4 & 0 & 5\\
        10R\_10V\_04F\_03M & \hspace{0.1cm} & 3 & 2 & 0 & 5 & \hspace{0.2cm} & 0 & 4 & 1 & 3\\
        10R\_10V\_04F\_04M & \hspace{0.1cm} & 3 & 2 & 0 & 5 & \hspace{0.2cm} & 0 & 3 & 2 & 3\\
        \addlinespace
        12R\_12V\_02F\_03M & \hspace{0.1cm} & 1 & 4 & 0 & 5 & \hspace{0.2cm} & 0 & 5 & 0 & 5\\
        12R\_12V\_02F\_04M & \hspace{0.1cm} & 2 & 3 & 0 & 5 & \hspace{0.2cm} & 0 & 4 & 1 & 4\\
        12R\_12V\_04F\_03M & \hspace{0.1cm} & 2 & 3 & 0 & 5 & \hspace{0.2cm} & 0 & 4 & 1 & 2\\
        12R\_12V\_04F\_04M & \hspace{0.1cm} & 2 & 3 & 0 & 5 & \hspace{0.2cm} & 0 & 5 & 0 & 3\\
        \midrule
        Sum & \hspace{0.1cm} & 58 & 22 & 0 & 80 & \hspace{0.2cm} & 22 & 53 & 5 & 70\\
        \bottomrule
    \end{tabular}
    \begin{tablenotes}
        \scriptsize
        \item Instances killed: lr204 from 10R\_10V\_04F\_03M; lr203 and lr204 from 10R\_10V\_04F\_04M; lr204 from 12R\_12V\_02F\_04M; lr202 from 12R\_12V\_04F\_03M
    \end{tablenotes}
    }
\end{table}

\okC{
The results in Table \ref{tab:overview_results_multi_island} show that the solver obtained, for the multi-island instances, feasible solutions for all Type~1 instances and 71 out of 80 (88.8\%) Type~2 instances. 
Table \ref{tab:overview_results_multi_floor} shows that, for the multi-floor instances, the solver obtained feasible solutions for all Type~1 instances and 70 out of 80 (87.5\%) Type~2 instances.
\okC{Together, there are} 301 out of 320 (94.1\%) instances with known feasible solutions.
For the 19 remaining instances without feasible solutions, the solver was killed by the operating system in 12 of them (63.2\%), seven multi-island, and five multi-floor instances.
Besides, all the killed processes were running Type~2 instances with at least ten requests.
The solver proved solution optimality for 49 and 58 out of 80 (61.3\% and 72.5\%) Type~1 instances for the multi-island and multi-floor families, respectively.
On the other hand, 24 and 22 out of 80 (30.0\% and 27.5\%) Type~2 instances were solved to optimality on the multi-island and multi-floor families, respectively.
Remarkably, the solver proved solution optimality for all Type~1 multi-floor instances with six and eight requests, in addition to all Type~1 multi-island instances with six requests.
Furthermore, the solver still proved solution optimality for 29 out of 60 (48.3\%) Type~1 multi-island instances with at least eight requests, and 18 out of 40 (45.0\%) Type~1 multi-floor instances with at least ten requests.
For the Type~2 instances, the solver proved solution optimality for 18 and 17 out of 20 (90.0\% and 85.0\%) instances with six requests for the multi-island and the multi-floor families, respectively.
However, the optimal solution rate dropped significantly for instances with at least eight requests. 
Only 6 out of 60 (10.0\%) instances for the multi-island family and 5 out of 60 (8.3\%) for the multi-floor family were solved to optimality in such scenarios.
\okC{In total, the number of instances in which the solver proved solution optimality was 153 out of 320 (47.8\%).}
}

Tables \ref{tab:mean_summary_results_multi_island} and \ref{tab:mean_summary_results_multi_floor} summarize the results of the MIP formulation for the two instance families.
Each line in these tables represents the five instances of the corresponding instance group, indicated in the first column.
The next columns show, for each instance type, the mean of: the best incumbent solution value (Sol.), the percentage gap between the best incumbent solution value and the best bound (Gap), and the running time in seconds (Time).
The value “N/A” is placed in the table whenever there is at least one instance within each type of group for which the solver did not find a feasible solution within the time limit or was killed by the operating system.
The row “Mean” presents the mean of each column, excluding the groups with an “N/A” value, regardless of whether the “N/A” occurs in Type~1 or Type~2 of the group.

\begin{table}[!ht]
    \centering
    \caption{\okC{Mean summary of the baseline MIP results for the multi-island family}}
    \label{tab:mean_summary_results_multi_island}
    \scriptsize
    \okC{
    \begin{tabular}{lrrrrrrrrr}\toprule
        &\hspace{0.1cm} &\multicolumn{3}{c}{Type~1} &\hspace{0.2cm} &\multicolumn{3}{c}{Type~2} \\
        \cmidrule(lr){3-5} \cmidrule(lr){7-9}
        Group &\hspace{0.1cm} &Sol. &Gap (\%) &Time(s) &\hspace{0.2cm} &Sol. &Gap (\%) &Time(s) \\
        \cmidrule(lr){1-1} \cmidrule(lr){3-5} \cmidrule(lr){7-9}
        06R\_06V\_02I\_03M &\hspace{0.1cm} &534.04 &0.0 &7.1 &\hspace{0.2cm} &597.86 &0.0 &590.0 \\
        06R\_06V\_02I\_04M &\hspace{0.1cm} &528.45 &0.0 &8.2 &\hspace{0.2cm} &594.99 &0.0 &709.0 \\
        06R\_06V\_04I\_03M &\hspace{0.1cm} &566.87 &0.0 &22.3 &\hspace{0.2cm} &603.80 &20.0 &838.8 \\
        06R\_06V\_04I\_04M &\hspace{0.1cm} &557.76 &0.0 &20.2 &\hspace{0.2cm} &598.21 &20.0 &866.8 \\
        \addlinespace
        08R\_08V\_02I\_03M &\hspace{0.1cm} &674.08 &18.5 &1474.7 &\hspace{0.2cm} &721.50 &57.3 &2416.6 \\
        08R\_08V\_02I\_04M &\hspace{0.1cm} &674.00 &17.3 &1472.5 &\hspace{0.2cm} &711.60 &57.2 &2535.1 \\
        08R\_08V\_04I\_03M &\hspace{0.1cm} &779.23 &11.5 &1458.3 &\hspace{0.2cm} &776.39 &62.1 &2982.1 \\
        08R\_08V\_04I\_04M &\hspace{0.1cm} &758.48 &6.5 &1505.3 &\hspace{0.2cm} &767.44 &72.4 &2964.8 \\
        \addlinespace
        10R\_10V\_02I\_03M &\hspace{0.1cm} &854.23 &28.8 &1472.3 &\hspace{0.2cm} &774.67 &81.7 &3600.2 \\
        10R\_10V\_02I\_04M &\hspace{0.1cm} &841.11 &32.5 &1485.1 &\hspace{0.2cm} &793.43 &80.9 &3600.3 \\
        10R\_10V\_04I\_03M &\hspace{0.1cm} &840.74 &34.4 &2171.3 &\hspace{0.2cm} &N/A &N/A &N/A \\
        10R\_10V\_04I\_04M &\hspace{0.1cm} &827.11 &35.9 &1840.5 &\hspace{0.2cm} &N/A &N/A &N/A \\
        \addlinespace
        12R\_12V\_02I\_03M &\hspace{0.1cm} &1048.66 &40.4 &2224.6 &\hspace{0.2cm} &N/A &N/A &N/A \\
        12R\_12V\_02I\_04M &\hspace{0.1cm} &1024.37 &41.0 &2352.6 &\hspace{0.2cm} &N/A &N/A &N/A \\
        12R\_12V\_04I\_03M &\hspace{0.1cm} &955.11 &40.6 &3059.0 &\hspace{0.2cm} &N/A &N/A &N/A \\
        12R\_12V\_04I\_04M &\hspace{0.1cm} &946.13 &40.7 &3190.7 &\hspace{0.2cm} &N/A &N/A &N/A \\
        \midrule
        Mean &\hspace{0.1cm} &676.82 &11.5 &892.6 &\hspace{0.2cm} &693.99 &45.2 &2110.4 \\
        \bottomrule
    \end{tabular}
    }
\end{table}

\begin{table}[!ht]
    \centering
    \caption{\okC{Mean summary of the baseline MIP results for the multi-floor family}}
    \label{tab:mean_summary_results_multi_floor}
    \scriptsize
    \okC{
    \begin{tabular}{lrrrrrrrrr}\toprule
        &\hspace{0.1cm} &\multicolumn{3}{c}{Type~1} &\hspace{0.2cm} &\multicolumn{3}{c}{Type~2} \\
        \cmidrule(lr){3-5} \cmidrule(lr){7-9}
        Group &\hspace{0.1cm} &Sol. &Gap (\%) &Time(s) &\hspace{0.2cm} &Sol. &Gap (\%) &Time(s) \\
        \cmidrule(lr){1-1} \cmidrule(lr){3-5} \cmidrule(lr){7-9} 
        06R\_06V\_02F\_03M &\hspace{0.1cm} &611.30 &0.0 &1.5 &\hspace{0.2cm} &601.77 &0.0 &710.7 \\
        06R\_06V\_02F\_04M &\hspace{0.1cm} &606.03 &0.0 &1.9 &\hspace{0.2cm} &601.52 &3.7 &785.8 \\
        06R\_06V\_04F\_03M &\hspace{0.1cm} &778.45 &0.0 &0.5 &\hspace{0.2cm} &713.27 &14.5 &839.4 \\
        06R\_06V\_04F\_04M &\hspace{0.1cm} &770.46 &0.0 &0.9 &\hspace{0.2cm} &697.11 &11.3 &900.2 \\
        \addlinespace
        08R\_08V\_02F\_03M &\hspace{0.1cm} &855.86 &0.0 &524.4 &\hspace{0.2cm} &756.58 &59.8 &2895.5 \\
        08R\_08V\_02F\_04M &\hspace{0.1cm} &848.18 &0.0 &824.3 &\hspace{0.2cm} &756.00 &70.2 &2906.6 \\
        08R\_08V\_04F\_03M &\hspace{0.1cm} &1084.83 &0.0 &224.1 &\hspace{0.2cm} &944.07 &74.4 &2986.0 \\
        08R\_08V\_04F\_04M &\hspace{0.1cm} &1068.07 &0.0 &852.2 &\hspace{0.2cm} &948.73 &74.0 &2946.6 \\
        \addlinespace
        10R\_10V\_02F\_03M &\hspace{0.1cm} &1064.17 &18.9 &2160.1 &\hspace{0.2cm} &893.29 &86.2 &3600.1 \\
        10R\_10V\_02F\_04M &\hspace{0.1cm} &1058.77 &14.0 &1825.2 &\hspace{0.2cm} &886.11 &77.6 &3599.7 \\
        10R\_10V\_04F\_03M &\hspace{0.1cm} &1342.23 &1.1 &1441.1 &\hspace{0.2cm} &N/A &N/A &N/A \\
        10R\_10V\_04F\_04M &\hspace{0.1cm} &1325.67 &3.7 &1445.6 &\hspace{0.2cm} &N/A &N/A &N/A \\
        \addlinespace
        12R\_12V\_02F\_03M &\hspace{0.1cm} &1302.23 &37.1 &2880.5 &\hspace{0.2cm} &1089.34 &97.9 &3600.5 \\
        12R\_12V\_02F\_04M &\hspace{0.1cm} &1291.96 &39.7 &2261.8 &\hspace{0.2cm} &N/A &N/A &N/A \\
        12R\_12V\_04F\_03M &\hspace{0.1cm} &1584.13 &22.1 &2161.9 &\hspace{0.2cm} &N/A &N/A &N/A \\
        12R\_12V\_04F\_04M &\hspace{0.1cm} &1558.61 &24.4 &2162.2 &\hspace{0.2cm} &N/A &N/A &N/A \\
        \midrule
        Mean &\hspace{0.1cm} &913.49 &6.4 &845.1 &\hspace{0.2cm} &807.98 &51.8 &2342.8 \\
        \bottomrule
    \end{tabular}
    }
\end{table}

\okC{
We can notice in Tables \ref{tab:mean_summary_results_multi_island} and \ref{tab:mean_summary_results_multi_floor} that in both multi-island and Type~1 multi-floor instances with six requests, the solver proved solution optimality in less than 23 seconds on average for each group.
In contrast, the mean times and mean gaps significantly increased in both Type~2 multi-island and multi-floor instances with six requests, particularly when the number of regions increased.
In Type~1 instances, when the number of requests increased from six to eight, the mean gaps only increased for multi-island instances, with a minimum mean gap of 6.5\%.
The mean times also increased, reaching a maximum variation of 0.9 seconds in group 06R\_06V\_04F\_04M to 852.2 seconds in group 08R\_08V\_04F\_04M.
Besides that, as seen in Tables \ref{tab:overview_results_multi_island} and \ref{tab:overview_results_multi_floor}, the solver exceeded the time limit in 30 out of the 40 (75.0\%) Type~2 instances with eight requests \okC{combining both families}.
Thus, the mean times were above 2416.6 seconds in each of these groups.
Also, when looking at the mean gaps of these groups, it is noticeable that the solver had difficulty increasing the best bound: the mean gaps ranged from 57.2\% to 74.4\% in general for the Type~2 instances with eight requests.
This scenario is even worse for Type~2 multi-island instances with at least ten requests, especially in groups with 12 requests and two floors.
On the other hand, the mean gaps ranged from 1.1\% to 41.0\% for Type~1 instances with at least 10 requests.
In general, considering only groups without “N/A” values, the mean gap and time are similar between multi-island and multi-floor instances. 
For Type~1 instances, the mean gaps were 11.5\% (multi-island) and 6.4\% (multi-floor), with mean times of 892.6 seconds and 845.1 seconds, respectively. 
Conversely, for Type~2 instances, the mean gaps were 45.2\% (multi-island) and 51.8\% (multi-floor), with mean times of 2110.4 seconds and 2342.8 seconds, respectively. 
Besides that, both the mean gap and time increase significantly from Type~1 to Type~2 instances.
The maximum variation on the mean gap ranged from 6.4\% to 51.8\%, with a mean time ranging from 845.1 seconds to 2342.8 seconds.
}

\okC{
In \ref{app:grb_vs_hx}, a comparison is provided between the solution values obtained by the solvers Gurobi and Hexaly on the baseline MIP for the instances with up to 12 requests.
Briefly, we show that all Gurobi's solution values are at least as good as those obtained by Hexaly. 
Hence, we decided not to conduct further experiments using Hexaly.
}

\subsection{\okC{The impact of valid inequalities on the MIP formulation}}
\label{subsec:impact_valid_inequalities}

\okC{
We conducted preliminary experiments to evaluate the impact of the valid inequalities proposed in Section \ref{sec:mip_tightening} on the  MIP formulation.
In total, there are 12 families of valid inequalities proposed, described by the constraints \eqref{mip:vi:00}-\eqref{mip:vi:55}.
The experiments consider 26 different valid inequality configurations: one containing no valid inequalities for reference, 12 encompassing only one valid inequality family, 12 comprising all valid inequality families except one, and one involving all valid inequalities.
It is worth mentioning that we only considered paths of size two for the valid inequalities \eqref{mip:vi:20}-\eqref{mip:vi:25}.
For these preliminary experiments, we imposed a time limit of 600 seconds (10 minutes) on each solver execution.
The remaining solver settings and computational environment are the same as those from the previous MIP experiments (see Section \ref{sec:experiments}).
Finally, these preliminary tests were conducted on part of the additional instances generated from the sixth to the tenth instance of each type for the PDPTW: lr106-lr110 (Type~1) and lr206-lr210 (Type~2).
These instances were generated in the same way as described in Section \ref{subsec:inst_gen}, varying the number of requests ($n \in \{6, 8, 10, 12\}$), vehicles ($|K| = n$), regions ($z \in \{2,4\}$) and machines ($|H| \in \{3,4\}$).
We randomly selected 16 out of those instances while maintaining the distribution of different instance size characteristics.
}

\okC{
Table \ref{tab:ranking_vi_configs} presents the ranking of the valid inequality configurations considered.
The criterion adopted for ranking the configurations is based on the feasibility rate, the Average Relative Percentage Deviation (ARPD), and the percentage gap between the best incumbent solution and the best bound (Gap).
The feasibility rate of a configuration is calculated as $\frac{n^{feas}}{n^{inst}}$, where $n^{feas}$ is the number of instances in which the solver found a feasible solution, and $n^{inst} = 16$ is the total number of instances selected.
The ARPD of a configuration is the average of the Relative Percentage Deviations (RPDs) for each configuration, each one calculated as $\frac{sol - best}{best}$, where $sol$ is the value of the solution obtained by the solver for a valid inequality configuration, and $best$ is the value of the best solution achieved by the solver for all valid inequality configurations.
To compare these three metrics, we took the complement of the ARPD and Gap metrics, normalized using the min-max feature scaling.
\okC{
Consequently, the metrics fall within the interval $[0,1]$, and the higher the metric value, the better.}
Then, we took the average of those three metrics and arranged the table in decreasing order by this average.
The first column ``Pos.'' displays the position of the valid inequality configuration in the ranking.
The second column ``Configuration'' presents which valid inequality constraints were applied.
\okC{
Columns ``Feas. (\%)'', ``$\overline{\text{ARPD}}$ (\%)'', and ``$\overline{\text{Gap}}$ (\%)'' exhibit the feasibility rate, ARPD, and GAP-related metrics.
The column ``Avg. (\%)'' provides the average of the metrics.}
The last column ``Opt.'' shows the number of instances in which the solver found an optimal solution for the configuration out of the 16 instances selected.
}

\begin{table}[ht]
\centering
\scriptsize
\caption{\okC{Ranking of the valid inequality configurations \label{tab:ranking_vi_configs}}} 
\okC{
\begin{tabular}{r l c rrr c r c r}
  \toprule
  & & & \multicolumn{3}{c}{Normalized metrics} & & \\
  \cmidrule{4-6}
  Pos. & Configuration & & Feas. (\%) & $\overline{\text{ARPD}}$ (\%) & $\overline{\text{Gap}}$ (\%) & & Avg. (\%) & Opt. \\ 
  \midrule
  1  & (35)-(39),(41)-(46)   & & 100.00  & 91.90     & 81.18     & & 91.02 & 5 \\ 
  2  & (35),(37)-(46)        & & 100.00  & 100.00    & 64.13     & & 88.04 & 4 \\ 
  3  & (35)-(36),(38)-(46)   & & 100.00  & 94.60     & 63.44     & & 86.01 & 4 \\ 
  4  & (35)-(44),(46)        & & 93.75   & 83.32     & 76.74     & & 84.60 & 4 \\ 
  5  & (36)-(46)             & & 87.50   & 56.43     & 100.00    & & 81.31 & 4 \\ 
  6  & (35)-(38),(40)-(46)   & & 93.75   & 78.93     & 67.19     & & 79.96 & 5 \\ 
  7  & (35)-(41),(43)-(46)   & & 93.75   & 77.14     & 63.57     & & 78.15 & 4 \\ 
  8  & (45)                  & & 93.75   & 75.19     & 63.26     & & 77.40 & 4 \\ 
  9  & (35)-(37),(39)-(46)   & & 87.50   & 62.54     & 74.76     & & 74.93 & 4 \\ 
  10 & (35)-(43),(45)-(46)   & & 87.50   & 64.48     & 70.91     & & 74.30 & 4 \\ 
  11 & (40)                  & & 87.50   & 55.41     & 64.41     & & 69.11 & 4 \\ 
  12 & (35)-(42),(44)-(46)   & & 87.50   & 63.62     & 53.05     & & 68.06 & 4 \\ 
  13 & (35)-(46)             & & 81.25   & 44.43     & 60.80     & & 62.16 & 3 \\ 
  14 & (36)                  & & 81.25   & 40.70     & 51.39     & & 57.78 & 3 \\ 
  15 & (35)-(40),(42)-(46)   & & 81.25   & 41.59     & 38.58     & & 53.81 & 4 \\ 
  16 & (35)-(45)             & & 81.25   & 46.12     & 29.23     & & 52.20 & 3 \\ 
  17 & (41)                  & & 75.00   & 21.70     & 57.07     & & 51.26 & 4 \\ 
  18 & (44)                  & & 75.00   & 23.22     & 52.36     & & 50.19 & 4 \\ 
  19 & (42)                  & & 75.00   & 21.81     & 41.66     & & 46.15 & 3 \\ 
  20 & None                  & & 75.00   & 19.17     & 42.71     & & 45.63 & 3 \\ 
  21 & (46)                  & & 68.75   & 3.58      & 55.46     & & 42.60 & 4 \\ 
  22 & (38)                  & & 75.00   & 24.53     & 23.62     & & 41.05 & 3 \\ 
  23 & (39)                  & & 75.00   & 23.97     & 20.78     & & 39.92 & 4 \\ 
  24 & (37)                  & & 68.75   & 3.10      & 25.63     & & 32.49 & 3 \\ 
  25 & (43)                  & & 68.75   & 3.08      & 16.52     & & 29.45 & 3 \\ 
  26 & (35)                  & & 68.75   & 0.00      & 0.00      & & 22.92 & 3 \\ 
   \bottomrule
\end{tabular}
}
\end{table}

\okC{
As can be seen in Table \ref{tab:ranking_vi_configs}, there are 19 out of the 26 (73.1\%) valid inequality configurations that performed better than the solver on the baseline MIP formulation.
Besides, nine out of the top ten (90.0\%) valid inequality configurations consider more than one family of valid inequalities.
The top three valid inequality configurations found feasible solutions for all instances considered.
Finally, the best valid inequality configuration disregards the valid inequalities \eqref{mip:vi:25} and obtained an average score of 91.02\%, reaching an optimal solution for five out of the 16 (31.33\%) instances.
}

\okC{
Then, we selected the top-ranked valid inequality configuration to run additional experiments on the benchmark instances.
We executed these experiments using the same computational environment described at the beginning of Section \ref{sec:experiments}.
Tables~\ref{tab:overview_results_multi_island_vi} and \ref{tab:overview_results_multi_floor_vi} present a summary overview of the impact of the valid inequalities on the solver.
The columns are the same as in Tables~\ref{tab:overview_results_multi_island} and \ref{tab:overview_results_multi_floor}.
Rows “MIP” and “MIP+VI” summarize each metric by instance type for the baseline MIP formulation and the MIP with valid inequalities, respectively.
}

\FloatBarrier

\begin{table}[!ht]
    \centering
    \caption{\okC{Summary overview of the results of the MIP with valid inequalities for the multi-island family.}}
    \label{tab:overview_results_multi_island_vi}
    \scriptsize
    \okC{
    \begin{tabular}{llrrrrlrrrr}
        \toprule
        \multicolumn{1}{c}{} &\hspace{0.1cm} &\multicolumn{4}{c}{Type~1} & \hspace{0.2cm} &\multicolumn{4}{c}{Type~2} \\
        \cmidrule(lr){3-6}
        \cmidrule(lr){8-11}
        Group & \hspace{0.1cm}  & Opt. & TLE & Killed & Feas. & \hspace{0.2cm} & Opt. & TLE & Killed & Feas.\\
        \cmidrule(lr){1-1}
        \cmidrule(lr){3-6}
        \cmidrule(lr){8-11}
        MIP & \hspace{0.1cm} & 49 & 31 & 0 & 80 & \hspace{0.2cm} & 24 & 49 & 7 & 71\\
        MIP+VI & \hspace{0.1cm} & 52 & 28 & 0 & 80 & \hspace{0.2cm} & 23 & 50 & 7 & 71\\
    \bottomrule
    \end{tabular}
    \begin{tablenotes}
      \scriptsize
      \item Instances killed with valid inequalities:  lr204 from 10R\_10V\_04I\_03M; lr203 and lr204 from 10R\_10V\_04I\_04M; lr204 from 12R\_12V\_02I\_03M; lr204 from 12R\_12V\_02I\_04M; lr203 from 12R\_12V\_04I\_03M; lr202 from 12R\_12V\_04I\_04M
    \end{tablenotes}
    }
\end{table}

\begin{table}[!ht]
    \centering
    \caption{\okC{Summary overview of the results of the MIP with valid inequalities for the multi-floor family.}}
    \label{tab:overview_results_multi_floor_vi}
    \scriptsize
    \okC{
    \begin{tabular}{llrrrrlrrrr}
        \toprule
        \multicolumn{1}{c}{} &\hspace{0.1cm} &\multicolumn{4}{c}{Type~1} & \hspace{0.2cm} &\multicolumn{4}{c}{Type~2} \\
        \cmidrule(lr){3-6}
        \cmidrule(lr){8-11}
        Group & \hspace{0.1cm}  & Opt. & TLE & Killed & Feas. & \hspace{0.2cm} & Opt. & TLE & Killed & Feas.\\
        \cmidrule(lr){1-1}
        \cmidrule(lr){3-6}
        \cmidrule(lr){8-11}
        MIP & \hspace{0.1cm} & 58 & 22 & 0 & 80 & \hspace{0.2cm} & 22 & 53 & 5 & 70\\
        MIP+VI & \hspace{0.1cm} & 59 & 21 & 0 & 80 & \hspace{0.2cm} & 22 & 53 & 5 & 73\\
        \bottomrule
    \end{tabular}
    \begin{tablenotes}
      \scriptsize
      \item Instances killed with valid inequalities: lr204 from 10R\_10V\_04F\_03M; lr204 from 10R\_10V\_04F\_04M; lr203 from 12R\_12V\_04F\_03M; lr202 from 12R\_12V\_04F\_04M; lr203 from 12R\_12V\_04F\_04M
    \end{tablenotes}
    }
\end{table}

\okC{As seen in Tables~\ref{tab:overview_results_multi_island_vi} and \ref{tab:overview_results_multi_floor_vi}, the solver achieved the same number of feasible solutions for each instance type and family, except for Type~2 multi-floor instances, which increased from 70 to 73 out of 80 (from 87.5\% to 91.2\%).
Thus, the total number of instances in which the solver found feasible solutions increased from 301 to 304 out of 320 (from 94.1\% to 95.0\%) \okC{with the new strengthened formulation}.
Despite the number of optimal solutions for the Type~2 multi-island instances decreased from 24 to 23 out of 80 (from 30.0\% to 28.8\%), the number of optimal solutions for the Type~1 multi-island instances increased from 49 to 52 out of 80 (from 61.3\% to 65.0\%).
Besides, the solver proved solution optimality for an additional Type~1 multi-floor instance, increasing from 58 to 59 out of 80 (from 72.5\% to 73.8\%).
Hence, the total number of instances in which the solver proved solution optimality increased from 153 to 156 out of 320 (from 47.8\% to 48.8\%) \okC{with the new strengthened formulation}.
Although the number of killed instance executions is equal for each instance family and type, the two killed instance sets are different.}

\begin{table}[h]
    \centering
    \okC{
    \resizebox{\textwidth}{!}{
    \begin{tabular}{l r rrrr c rrrr}\toprule
        & &\multicolumn{4}{c}{Multi-island} & &\multicolumn{4}{c}{Multi-floor} \\
        \cmidrule{3-6}
        \cmidrule{8-11}
         & &MIP or MIP+VI & MIP and MIP+VI &MIP &MIP+VI & &MIP or MIP+VI & MIP and MIP+VI &MIP &MIP+VI \\
         \cmidrule{1-1}
        \cmidrule{3-6}
        \cmidrule{8-11}
        Type~1 & &52/80 &49/80 &0/80 &3/80 & &61/80 &56/80 &2/80 &3/80 \\
        Type~2 & &35/80 &30/80 &3/80 &2/80 & &35/80 &26/80 &6/80 &3/80 \\
        \bottomrule
    \end{tabular}
    }
    }
    \caption{\okC{Missing observations in each gap RPD box plot}}
    \label{tab:missing_obs_gap_dev_box_plots}
\end{table}

\begin{figure}[!ht]
    \centering
    \subfigure[Multi-island instances]{
        \includegraphics[width=0.45\linewidth]{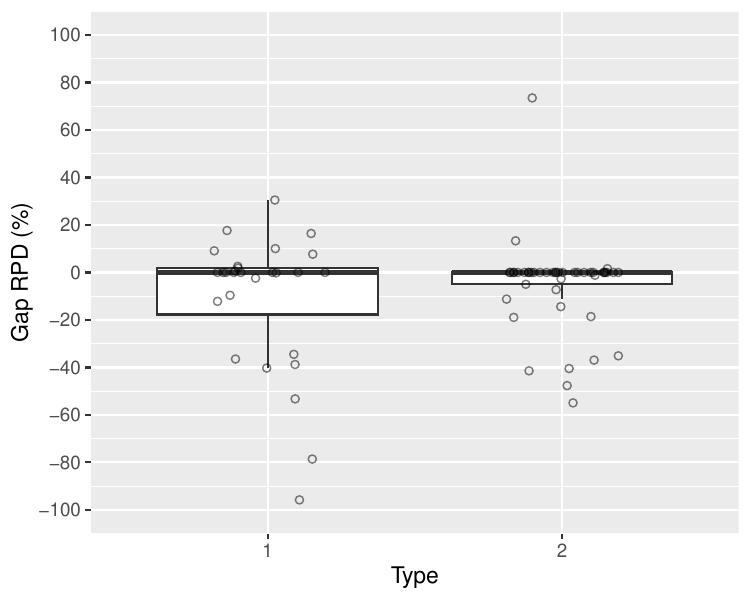}
    }
    \subfigure[Multi-floor instances]{
        \includegraphics[width=0.45\linewidth]{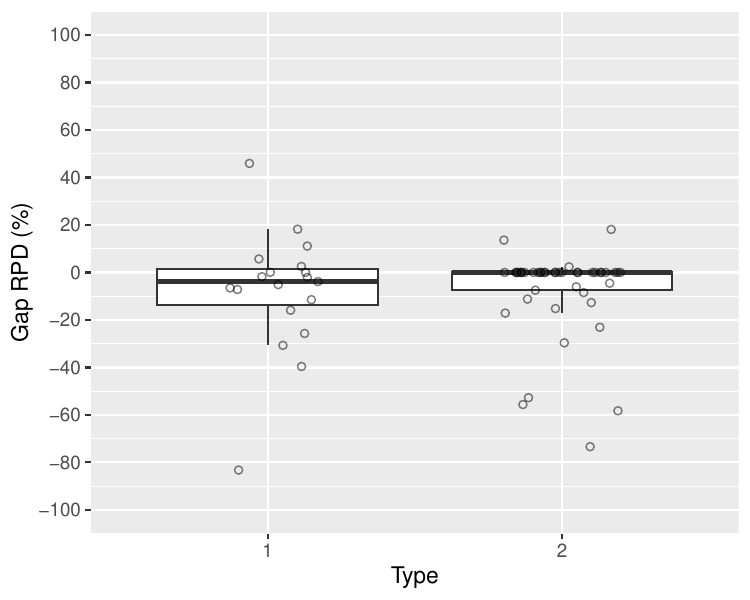}
    }
    \caption{\okC{Box plots of the gap RPDs by instance type for the baseline MIP and the MIP with valid inequalities. Negative RPDs indicate that the version with valid inequalities outperforms the baseline MIP. Observations only for instances where a feasible solution was found, but the time limit was exceeded.}}
    \label{fig:gap_deviation_with_without_vi}
\end{figure}

\okC{
The box plots of Figure \ref{fig:gap_deviation_with_without_vi} depict the solver gap RPDs by instance type for the baseline MIP and the MIP with valid inequalities.
We consider only the gaps for instances where both approaches found a feasible solution, but the time limit was exceeded.
Table \ref{tab:missing_obs_gap_dev_box_plots} summarizes the number of missing observations in each box plot, separated by instance family and type.
Column ``MIP or MIP+VI" shows the total number of missing gaps.
Column ``MIP and MIP+VI" indicates cases where both approaches \okC{miss each one its  respective gap}.
Columns ``MIP" and ``MIP+VI" count cases where only for the corresponding formulation (the baseline MIP and the MIP with valid inequalities, respectively), the gap value is missing.
The gap RPD is calculated as $100 \times\frac{gap^{\text{MIP+VI}} - gap^{\text{MIP}}}{gap^{\text{MIP}}}$, in which $gap^{\text{MIP}}$ and $gap^{\text{MIP+VI}}$ are the gaps of the solutions for an instance obtained by the solver on the baseline MIP and on the MIP with valid inequalities, respectively.
Thereupon, negative values indicate that the gap of the solver on the MIP with valid inequalities outperformed the gap of the baseline MIP.
Otherwise, no improving gap was found.
}

\okC{
As seen in the box plots of Figure~\ref{fig:gap_deviation_with_without_vi}, the valid inequalities contributed to at least match or improve the gaps of at least 75\% of the observations for each instance family and type.
Additionally, the box plots of Type~1 instances have the largest gap RPD ranges, although the improved observations have the largest gap RPDs.
Conversely, there are few positive gap observations in the box plots of Type~2 instances, and most of them are outliers.
}

\begin{table}[h]
    \centering
    \okC{
    \resizebox{\textwidth}{!}{
    \begin{tabular}{l r rrrr c rrrr}\toprule
        & &\multicolumn{4}{c}{Multi-island} & &\multicolumn{4}{c}{Multi-floor} \\
        \cmidrule{3-6}
        \cmidrule{8-11}
         & &MIP or MIP+VI &  MIP and MIP+VI &MIP &MIP+VI & &MIP or MIP+VI & MIP and MIP+VI &MIP &MIP+VI \\
        \cmidrule{1-1}
        \cmidrule{3-6}
        \cmidrule{8-11}
        Type~1 & &0/80  &0/80 &0/80 &0/80 & &0/80  &0/80 &0/80 &0/80 \\
        Type~2 & &11/80 &7/80 &2/80 &2/80 & &12/80 &5/80 &5/80 &2/80 \\
        \bottomrule
    \end{tabular}
    }
    }
    \caption{\okC{Missing observations in each solution deviation box plot}}
    \label{tab:missing_obs_sol_dev_box_plots}
\end{table}

\begin{figure}[!ht]
    \centering
    \subfigure[Multi-island instances]{
        \includegraphics[width=0.45\linewidth]{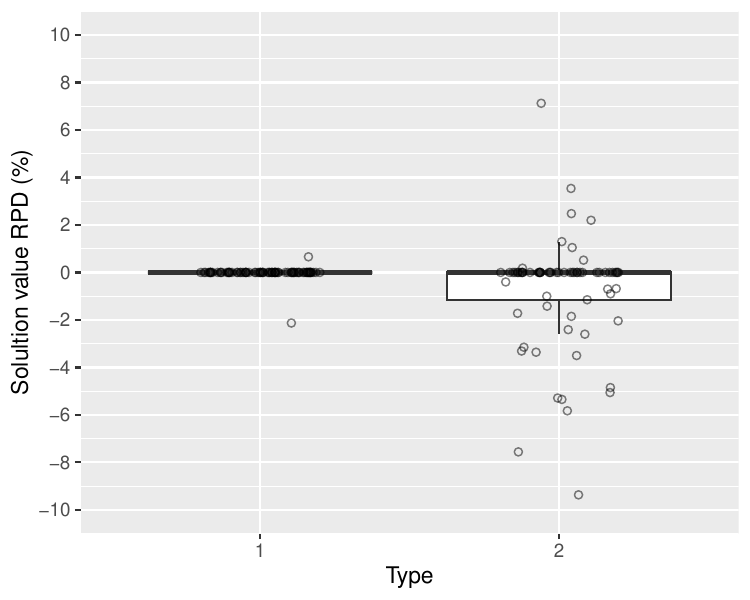}
    }
    \subfigure[Multi-floor instances]{
        \includegraphics[width=0.45\linewidth]{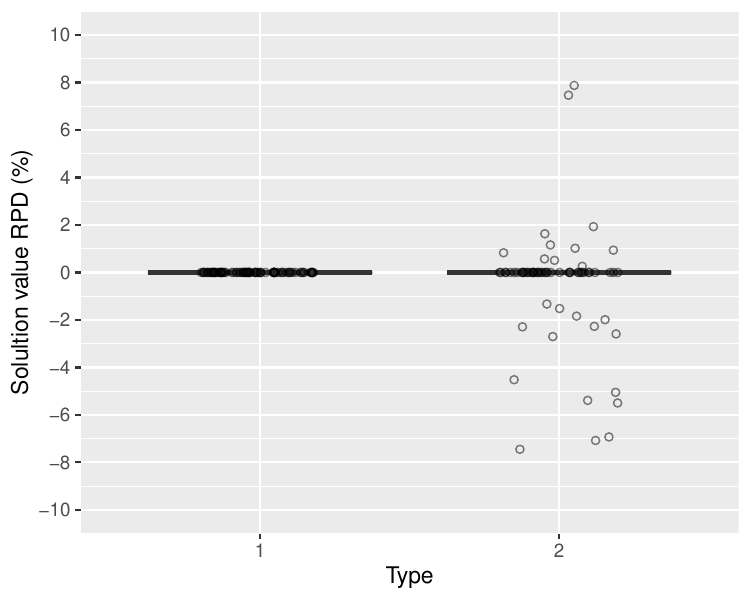}
    }
    \caption{\okC{Box plots of the solution value RPDs by instance type for the baseline MIP and the MIP with valid inequalities. Negative RPD indicate that the version with valid inequalities outperforms the baseline MIP. Observations only for instances where a feasible solution was found (optimal or not).}}
    \label{fig:sol_deviation_with_without_vi}
\end{figure}

\okC{
The box plots in Figure \ref{fig:sol_deviation_with_without_vi} depict the solver solution value RPDs by instance type for the baseline MIP and the MIP with valid inequalities.
We consider only the solution values for instances where both approaches found a feasible solution, whether optimal or not.
Table \ref{tab:missing_obs_sol_dev_box_plots} summarizes the number of missing observations in each box plot, separated by instance family and type.
Column ``MIP or MIP+VI" shows the total number of missing solution values.
Column ``MIP and MIP+VI" indicates cases where both approaches \okC{miss each one its respective solution value}.
Columns ``MIP" and ``MIP+VI" count cases where only the corresponding formulation (the baseline MIP and the MIP with valid inequalities, respectively) is missing a solution.
The solution value RPD is calculated as $100 \times\frac{sol^{\text{MIP+VI}} - sol^{\text{MIP}}}{sol^{\text{MIP}}}$, in which $sol^{\text{MIP}}$ and $sol^{\text{MIP+VI}}$ are the solution values for an instance obtained by the solver on the baseline MIP and on the MIP with valid inequalities, respectively.
Therefore, negative values indicate that the solution value of the MIP with valid inequalities surpassed the solution value of the solver on the baseline MIP.
Otherwise, no improving solution was found.
}

\okC{
We can observe from the box plots of Figure \ref{fig:sol_deviation_with_without_vi} that the valid inequalities mostly benefited Type~2 instances in terms of solution value improvement.
Notice that for at least 75\% of the Type~2 multi-island observations, the same solution or a better one was achieved.
Besides, there are more negative outliers than positive outliers on the Type~2 instances in both families.
Conversely, all solution value RPDs for Type~1 instances in both instance families are equal to zero or irrelevant, disregarding the outliers.
Given that Type~1 instances benefited more from valid inequalities on gap reductions (see Figure \ref{fig:gap_deviation_with_without_vi}), and there are no missing observations on the respective solution value RPD box plots (see Table \ref{tab:missing_obs_sol_dev_box_plots}), the solutions acquired for them are near-optimal solutions.
\okC{Finally, we conclude that the MIP with valid inequalities improves solver's performance not only in terms of optimality gap, but also in terms of solution value.}
}

\okC{
\ref{app:mip_results_details} presents the box plots of time and gap distributions across all instances related to the results of the MIP with valid inequalities.
Concerning the time distributions, although increasing the number of requests results in an exponential increase in the total time spent, many solver times are too sparse for the same-size requests.
For Type~2 instances, the solver reaches its time limit for most of the instances with at least eight requests, while it only occurs more frequently on Type~1 instances with 12 requests.
Regarding the gap distributions, the gaps are at least above 80\% for almost all Type~2 instances with at least ten requests in both families, and for more than 50\% of the instances with eight requests.
On the other hand, most of the gaps are equal to zero for Type~1 instances, as the number of optimal solutions is high (see Table \ref{tab:overview_results_multi_floor_vi}).
}

\subsection{Heuristic results}
\label{subsec:heur_results}

\okC{
In this section, we provide the results of the MSLP heuristic described in Section \ref{sec:msheurlp}.
First, we performed preliminary experiments (see \ref{app:preliminary_tests}) to define the parameters for MSLP.
In short, we selected a diverse set of instances outer of the test data, and conducted  these experiments with the same stop criterion, varying the values of the threshold parameter $\alpha \in \{0.05, 0.10, 0.15, 0.20\}$.
Based on the solution value RPDs for each parameter configuration, we set $\alpha$ to $0.05$.
}

For the test data, ten independent runs were performed for each instance. The results \okC{for the 360 instances with up to 12 requests} are summarized in Tables \ref{tab:mean_summary_mslp_results_multi_island} and \ref{tab:mean_summary_mslp_results_multi_floor} for the two instance families.
Each line in these tables represents the five instances of the corresponding instance group, indicated in the first column.
The next columns show, for each instance type, the mean of: the best incumbent solution value for the MIP with valid inequalities (Sol.), the minimum of the best MSLP solution values (Min. Sol.), the mean of the best MSLP solution values (Mean Sol.), the CPU time in seconds (Time (s)), and the CPU time to best in seconds (TTB (s)).
The value “N/A” is placed whenever there is at least one instance within each type of group for which the solver did not find a feasible solution within the time limit or was killed by the operating system.
The row “Mean --- MSLP” presents the mean of all the groups for the columns related to the multi-start heuristic.
The row “Mean --- Sol.” presents the mean of the columns for solution values, excluding the groups with an “N/A” value, regardless of whether the “N/A” occurs in Type~1 or Type~2 of the group.

% \begin{landscape}
% \begin{table}
\begin{sidewaystable}
    \centering
    \caption{\okC{Mean summary of the MSLP results for the smaller instances of the multi-island family}}
    \label{tab:mean_summary_mslp_results_multi_island}
    \scriptsize
    \okC{
    \begin{tabular}{lrcrrrrrrcrrrrrr}\toprule
        &\hspace{0.05cm} &\multicolumn{6}{c}{Type~1} &\hspace{0.05cm} &\multicolumn{6}{c}{Type~2} \\
        \cmidrule{3-8} \cmidrule{10-15}
        &\hspace{0.05cm} &MIP+VI &\hspace{0.05cm} &\multicolumn{4}{c}{MSLP} &\hspace{0.05cm} &MIP+VI &\hspace{0.05cm} &\multicolumn{4}{c}{MSLP} \\
        \cmidrule{3-3} \cmidrule{5-8} \cmidrule{10-10} \cmidrule{12-15}
        Group &\hspace{0.05cm} &Sol. &\hspace{0.05cm} &Min. Sol. &Mean Sol. &Time (s) &TTB (s) &\hspace{0.05cm} &Sol. &\hspace{0.05cm} &Min. Sol. &Mean Sol. &Time (s) &TTB (s) \\
        \midrule
        06R\_06V\_02I\_03M &\hspace{0.05cm} &534.04 &\hspace{0.05cm} &537.02 &537.02 &88.44 &0.27 &\hspace{0.1cm} &597.86 &\hspace{0.05cm} &647.21 &647.21 &96.04 &4.93 \\
        06R\_06V\_02I\_04M &\hspace{0.05cm} &528.45 &\hspace{0.05cm} &531.41 &531.41 &88.84 &0.33 &\hspace{0.1cm} &594.99 &\hspace{0.05cm} &645.02 &645.02 &98.05 &1.50 \\
        06R\_06V\_04I\_03M &\hspace{0.05cm} &566.87 &\hspace{0.05cm} &568.91 &568.91 &55.94 &0.26 &\hspace{0.1cm} &600.75 &\hspace{0.05cm} &656.47 &656.47 &112.86 &2.62 \\
        06R\_06V\_04I\_04M &\hspace{0.05cm} &557.76 &\hspace{0.05cm} &559.16 &559.16 &90.40 &0.34 &\hspace{0.1cm} &598.21 &\hspace{0.05cm} &652.49 &652.49 &116.14 &1.41 \\
        \addlinespace
        08R\_08V\_02I\_03M &\hspace{0.05cm} &674.08 &\hspace{0.05cm} &675.92 &675.92 &95.56 &0.76 &\hspace{0.1cm} &715.74 &\hspace{0.05cm} &731.14 &732.37 &124.95 &27.77 \\
        08R\_08V\_02I\_04M &\hspace{0.05cm} &674.00 &\hspace{0.05cm} &676.09 &676.09 &104.26 &0.64 &\hspace{0.1cm} &712.54 &\hspace{0.05cm} &727.26 &728.61 &129.04 &23.89 \\
        08R\_08V\_04I\_03M &\hspace{0.05cm} &779.23 &\hspace{0.05cm} &779.30 &779.30 &62.75 &3.54 &\hspace{0.1cm} &766.74 &\hspace{0.05cm} &772.36 &774.73 &148.09 &30.56 \\
        08R\_08V\_04I\_04M &\hspace{0.05cm} &758.48 &\hspace{0.05cm} &760.86 &760.86 &95.21 &3.55 &\hspace{0.1cm} &769.51 &\hspace{0.05cm} &768.85 &769.16 &153.97 &33.49 \\
        \addlinespace
        10R\_10V\_02I\_03M &\hspace{0.05cm} &854.23 &\hspace{0.05cm} &856.61 &856.61 &100.63 &6.91 &\hspace{0.1cm} &782.37 &\hspace{0.05cm} &815.60 &818.46 &155.51 &29.14 \\
        10R\_10V\_02I\_04M &\hspace{0.05cm} &841.11 &\hspace{0.05cm} &845.27 &845.31 &113.04 &7.94 &\hspace{0.1cm} &782.90 &\hspace{0.05cm} &814.55 &816.91 &162.97 &40.27 \\
        10R\_10V\_04I\_03M &\hspace{0.05cm} &840.74 &\hspace{0.05cm} &849.41 &849.67 &75.13 &11.92 &\hspace{0.1cm} &N/A &\hspace{0.05cm} &920.61 &928.14 &185.21 &80.05 \\
        10R\_10V\_04I\_04M &\hspace{0.05cm} &827.11 &\hspace{0.05cm} &833.79 &833.81 &119.35 &13.27 &\hspace{0.1cm} &N/A &\hspace{0.05cm} &911.72 &915.54 &208.80 &96.54 \\
        \addlinespace
        12R\_12V\_02I\_03M &\hspace{0.05cm} &1045.28 &\hspace{0.05cm} &1050.76 &1051.05 &94.15 &26.44 &\hspace{0.1cm} &N/A &\hspace{0.05cm} &1004.23 &1016.12 &192.95 &81.66 \\
        12R\_12V\_02I\_04M &\hspace{0.05cm} &1025.31 &\hspace{0.05cm} &1030.01 &1030.20 &113.83 &30.25 &\hspace{0.1cm} &N/A &\hspace{0.05cm} &1004.69 &1013.51 &203.21 &91.86 \\
        12R\_12V\_04I\_03M &\hspace{0.05cm} &955.11 &\hspace{0.05cm} &965.21 &968.95 &88.29 &24.62 &\hspace{0.1cm} &N/A &\hspace{0.05cm} &1039.11 &1054.80 &243.75 &111.40 \\
        12R\_12V\_04I\_04M &\hspace{0.05cm} &946.13 &\hspace{0.05cm} &953.82 &954.22 &135.50 &34.04 &\hspace{0.1cm} &N/A &\hspace{0.05cm} &1032.90 &1046.00 &267.79 &143.42 \\
        \midrule
        Mean --- MSLP &\hspace{0.05cm} &- &\hspace{0.05cm} &779.60 &779.90 &95.08 &10.32 &\hspace{0.1cm} &- &\hspace{0.05cm} &821.51 &825.97 &162.46 &50.03 \\
        Mean --- Sol. &\hspace{0.05cm} &676.82 &\hspace{0.05cm} &679.05 &679.06 &- &- &\hspace{0.1cm} &692.16 &\hspace{0.05cm} &723.10 &724.14 &- &- \\
        \bottomrule
    \end{tabular}
    }
% \end{sidewaystable}
\vspace{0.5cm}
% \begin{sidewaystable}
    \centering
    \caption{\okC{Mean summary of the MSLP results for the smaller instances of the multi-floor family}}
    \label{tab:mean_summary_mslp_results_multi_floor}
    \scriptsize
    \okC{
    \begin{tabular}{lrcrrrrrrcrrrrrr}\toprule
        &\hspace{0.05cm} &\multicolumn{6}{c}{Type~1} &\hspace{0.1cm} &\multicolumn{6}{c}{Type~2} & \\
        \cmidrule{3-8} \cmidrule{10-15}
        &\hspace{0.05cm} &MIP+VI &\hspace{0.05cm} &\multicolumn{4}{c}{MSLP} &\hspace{0.1cm} &MIP+VI &\hspace{0.05cm} &\multicolumn{4}{c}{MSLP} & \\
        \cmidrule{3-3} \cmidrule{5-8} \cmidrule{10-10} \cmidrule{12-15}
        Group &\hspace{0.05cm} &Sol. &\hspace{0.05cm} &Min. Sol. &Mean Sol. &Time (s) &TTB (s) &\hspace{0.1cm} &Sol. &\hspace{0.05cm} &Min. Sol. &Mean Sol. &Time (s) &TTB (s) \\\midrule
        06R\_06V\_02F\_03M &\hspace{0.05cm} &611.30 &\hspace{0.05cm} &617.09 &617.09 &96.81 &0.31 &\hspace{0.1cm} &601.77 &\hspace{0.05cm} &635.45 &635.45 &113.20 &1.46 \\
        06R\_06V\_02F\_04M &\hspace{0.05cm} &606.03 &\hspace{0.05cm} &611.87 &611.87 &100.56 &0.33 &\hspace{0.1cm} &601.52 &\hspace{0.05cm} &635.34 &635.34 &115.93 &3.04 \\
        06R\_06V\_04F\_03M &\hspace{0.05cm} &778.45 &\hspace{0.05cm} &789.03 &789.03 &89.77 &0.41 &\hspace{0.1cm} &713.27 &\hspace{0.05cm} &741.33 &741.33 &130.94 &8.12 \\
        06R\_06V\_04F\_04M &\hspace{0.05cm} &770.46 &\hspace{0.05cm} &783.34 &783.34 &93.73 &0.42 &\hspace{0.1cm} &697.11 &\hspace{0.05cm} &725.22 &725.22 &133.28 &1.94 \\
        \addlinespace
        08R\_08V\_02F\_03M &\hspace{0.05cm} &855.86 &\hspace{0.05cm} &858.49 &858.49 &133.20 &0.79 &\hspace{0.1cm} &756.58 &\hspace{0.05cm} &764.22 &764.29 &141.62 &11.98 \\
        08R\_08V\_02F\_04M &\hspace{0.05cm} &848.18 &\hspace{0.05cm} &851.49 &851.49 &136.17 &0.82 &\hspace{0.1cm} &754.03 &\hspace{0.05cm} &761.04 &761.10 &146.38 &17.35 \\
        08R\_08V\_04F\_03M &\hspace{0.05cm} &1084.83 &\hspace{0.05cm} &1092.18 &1092.18 &70.31 &1.03 &\hspace{0.1cm} &940.12 &\hspace{0.05cm} &954.79 &957.00 &174.31 &27.42 \\
        08R\_08V\_04F\_04M &\hspace{0.05cm} &1068.07 &\hspace{0.05cm} &1074.15 &1074.15 &114.75 &1.47 &\hspace{0.1cm} &937.41 &\hspace{0.05cm} &950.10 &952.42 &180.08 &27.60 \\
        \addlinespace
        10R\_10V\_02F\_03M &\hspace{0.05cm} &1064.17 &\hspace{0.05cm} &1070.36 &1070.51 &92.24 &6.69 &\hspace{0.1cm} &897.63 &\hspace{0.05cm} &881.37 &884.43 &181.18 &69.63 \\
        10R\_10V\_02F\_04M &\hspace{0.05cm} &1058.77 &\hspace{0.05cm} &1062.59 &1062.61 &134.89 &7.04 &\hspace{0.1cm} &882.61 &\hspace{0.05cm} &878.06 &882.57 &189.18 &65.52 \\
        10R\_10V\_04F\_03M &\hspace{0.05cm} &1342.23 &\hspace{0.05cm} &1361.22 &1361.83 &52.76 &7.22 &\hspace{0.1cm} &N/A &\hspace{0.05cm} &1116.73 &1126.59 &226.82 &116.32 \\
        10R\_10V\_04F\_04M &\hspace{0.05cm} &1325.67 &\hspace{0.05cm} &1331.88 &1331.88 &105.77 &5.50 &\hspace{0.1cm} &N/A &\hspace{0.05cm} &1111.74 &1117.62 &238.74 &124.46 \\
        \addlinespace
        12R\_12V\_02F\_03M &\hspace{0.05cm} &1302.23 &\hspace{0.05cm} &1307.48 &1307.48 &114.11 &3.11 &\hspace{0.1cm} &1090.67 &\hspace{0.05cm} &1063.61 &1073.25 &219.72 &109.27 \\
        12R\_12V\_02F\_04M &\hspace{0.05cm} &1291.96 &\hspace{0.05cm} &1297.40 &1297.40 &143.83 &4.37 &\hspace{0.1cm} &1085.13 &\hspace{0.05cm} &1059.27 &1070.69 &233.56 &111.75 \\
        12R\_12V\_04F\_03M &\hspace{0.05cm} &1584.13 &\hspace{0.05cm} &1591.64 &1597.14 &57.15 &25.47 &\hspace{0.1cm} &N/A &\hspace{0.05cm} &1278.11 &1294.37 &285.48 &119.78 \\
        12R\_12V\_04F\_04M &\hspace{0.05cm} &1558.53 &\hspace{0.05cm} &1572.46 &1573.06 &108.74 &53.34 &\hspace{0.1cm} &N/A &\hspace{0.05cm} &1274.66 &1288.06 &303.64 &147.97 \\
        \midrule
        Mean --- MSLP &\hspace{0.05cm} &- &\hspace{0.05cm} &1079.54 &1079.97 &102.80 &7.40 &\hspace{0.1cm} &- &\hspace{0.05cm} &926.94 &931.86 &188.38 &60.22 \\
        Mean --- Sol. &\hspace{0.05cm} &945.03 &\hspace{0.05cm} &951.29 &951.30 &- &- &\hspace{0.1cm} &829.82 &\hspace{0.05cm} &837.48 &840.26 &- &- \\
        \bottomrule
    \end{tabular}
    }
\end{sidewaystable}
% \end{table}
% \end{landscape}

\begin{figure}[!ht]
    \centering
    \subfigure[Multi-island instances]{
        \includegraphics[width=0.45\linewidth]{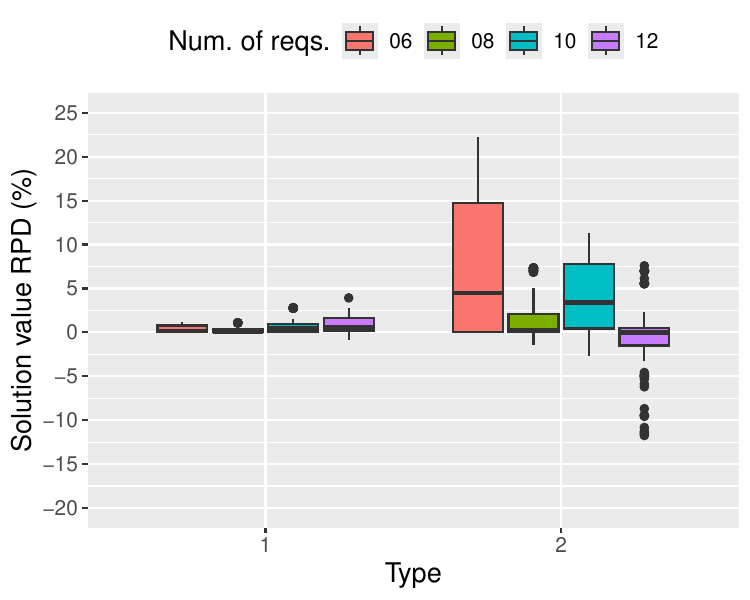}
    }
    \subfigure[Multi-floor instances]{
        \includegraphics[width=0.45\linewidth]{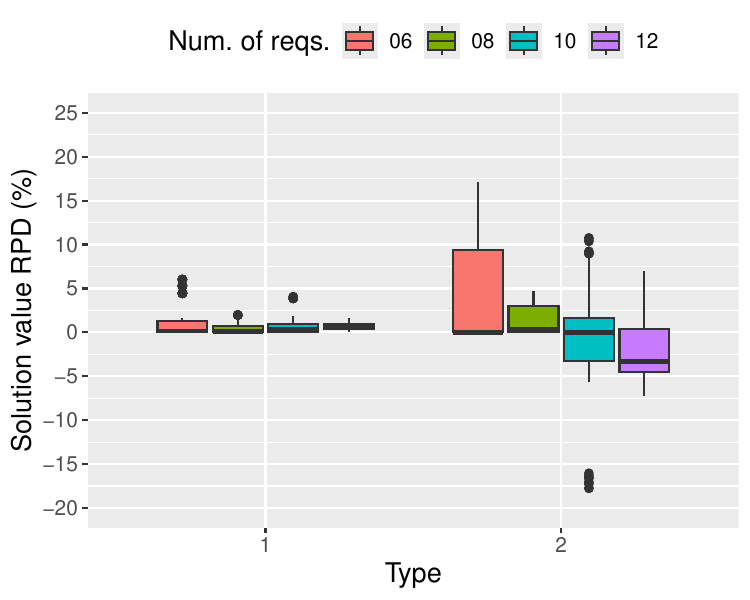}
    }
    \caption{\okC{Box plots of the solution value RPDs by instance type for the MIP with valid inequalities and the MSLP heuristic. Negative values indicate that the MSLP heuristic outperforms the MIP with valid inequalities. Observations only for instances with up to 12 requests where a feasible solution was found (optimal or not).}}
    \label{fig:dev_mip_mslp_boxplots}
\end{figure}

\okC{
We can notice in Tables \ref{tab:mean_summary_mslp_results_multi_island} and \ref{tab:mean_summary_mslp_results_multi_floor} that the heuristic found feasible solutions in all independent executions for all 320 instances with up to 12 requests of the data set.
Furthermore, the mean CPU times of the heuristic are slightly greater on the multi-floor instances for each type compared to the multi-island instances.
Besides, the MSLP took much less CPU time to find its best solution than the total CPU time spent.
For the multi-island and multi-floor Type~1 instances, the mean CPU times across all groups were 95.08 and 102.80 seconds, respectively.
Regarding the mean CPU times to best, it took only 10.32 and 7.40 seconds for the multi-island and multi-floor Type~1 instances, respectively.
For the multi-island and multi-floor Type~2 instances, the mean CPU times across all groups increased, reaching 162.46 and 188.38 seconds, respectively.
Moreover, the mean CPU times to best increased for the multi-island and multi-floor Type~2 instances, taking 50.03 and 60.22 seconds, respectively.
As discussed in more detail later in this section (see Figure \ref{fig:dev_mip_mslp_boxplots}), the MSLP heuristic can find solutions with reasonably low RPDs and short times when compared with the MIP with valid inequalities formulation, matching, or even surpassing solutions with at least eight requests.
It is worth mentioning that the mean CPU times for Type~1 instances do not increase as in Type~2 instances due to solution feasibility.
As described in Section \ref{sec:heur_gen_framework}, the LP improvement procedure is only executed when the (semi-)greedy solution is feasible.
Therefore, iterations without feasible solutions speed up the heuristic, as the LP improvement procedure is skipped.
Later in this section, we summarize the analysis of solution feasibility throughout the iterations of the MSLP presented in \ref{app:mslp_solution_feasibility}.
}

\okC{
The box plots of Figure \ref{fig:dev_mip_mslp_boxplots} depict an instance-by-instance solution value RPD comparison between the MIP with valid inequalities and the MSLP heuristic for the multi-island and multi-floor families, respectively.
For each instance type, the observations are grouped by the number of requests of the instance.
Each Type~1 box plot contains 800 observations.
For the Type~2 box plots, there are 90 and 70 missing observations out of the 800 (11.2\% and 8.8\%) corresponding to nine and seven instances in the multi-island and multi-floor families, respectively.
Among the tested instances, three of the multi-island and two of the multi-floor observations correspond to instances with ten requests, while the remaining correspond to instances with 12 requests.
All missing observations result from the solver's limitation to find a feasible solution for an instance using the MIP with valid inequalities.
The solution value RPDs are calculated as $100\times\frac{sol^{\text{MSLP}} - sol^{\text{MIP+VI}}}{sol^{\text{MIP+VI}}}$, where $sol^{\text{MSLP}}$ is the solution value obtained by the MSLP heuristic in one independent execution, and $sol^{\text{MIP+VI}}$ is the best solution value obtained by the \okC{solver on the MIP with valid inequalities}.
Therefore, negative values indicate that the MSLP heuristic achieved a better solution value compared to the solver on the MIP with valid inequalities.
Otherwise, no improving solution was found.
}

\okC{
Notice from the box plots of Figure \ref{fig:dev_mip_mslp_boxplots} that \okC{despite almost all observations are positive, the solution value RPDs are near zero} for Type~1 observations in both multi-island and multi-floor families.
Besides, the MSLP heuristic improved at least 25\% solution values for Type~2 multi-island instances with eight and ten requests.
Similarly, it also matched or improved solutions for at least 50\% of the Type~2 multi-floor instances with ten requests.
Furthermore, the best heuristic results happened for \okC{Type~2} instances with twelve requests in both instance families, where most of them improve or at least match the best solution values of the MIP with valid inequalities.
\okC{Thereupon, we advocate the usage of the solver on the MIP with valid inequalities on Type~1 instances, and some instances with up to eight requests.
For Type~2 instances with at least ten requests, the MSLP heuristic guarantees both solution feasibility and competitiveness.}
}

\okC{
We also executed tests on the larger instances (with 40 and 60 requests) of the second set.
From the 160 instances in this set, the MSLP heuristic found feasible solutions for all 80 Type~2 instances in the ten independent executions.
However, the heuristic failed to find a feasible solution for some Type~1 instances in the ten independent executions or in at least some of them. 
Table \ref{tab:overview_mslp_type_1_larger_split} presents a feasibility overview of the results for the Type~1 instances in both families.
Column ``All Feas." shows the number of instances in which the heuristic found a feasible solution in all ten independent executions.
Column ``[1,9] Feas." displays the number of instances where the heuristic found a feasible solution in at least one of the ten independent executions, but not in all of them.
Row ``Sum" presents the sum of each metric column.
Notice that for 34 out of the 40 (85.0\%) instances in both families, the MSLP heuristic found a feasible solution in all ten executions.
Besides, for only one of the remaining instances in each family, the MSLP heuristic found feasible solutions in only some executions: eight for the multi-island instance and three for the multi-floor instance.
Therefore, the MSLP heuristic found feasible solutions for 150 out of the 160 (93.8\%) instances with 40 and 60 requests.
Notably, most of the instances without feasible solutions have four regions (islands or floors) and five machines.
\okC{Henceforth, the heuristic indicates that some Type~1 instances might have very tight constraints, despite its feasibility ensurance (see Section \ref{subsubsec:inst-feas}).
In these scenarios, the Table \ref{tab:overview_mslp_type_1_larger_split} displays that providing the additional machine increased the number of independent executions generating feasible solutions.
On the other hand, the heuristic remains reliable regarding solution feasibility for Type~2 instances and for most Type~1 instances.
}}

\begin{table}[!htp]\centering
\caption{\okC{Feasibility overview of the MSLP results for the larger Type~1 instances}}
\label{tab:overview_mslp_type_1_larger_split}
\scriptsize
\begin{minipage}{0.48\textwidth}
\centering
\okC{
\text{Multi-island instances}
\label{tab:overview_mslp_type_1_larger_island}
\begin{tabular}{lrr}
\toprule
Group & All Feas. & [1,9] Feas. \\
\cmidrule{1-1} \cmidrule{2-3}
40R\_40V\_02I\_05M & 5 & 0 \\
40R\_40V\_02I\_06M & 5 & 0 \\
40R\_40V\_04I\_05M & 3 & 0 \\
40R\_40V\_04I\_06M & 5 & 0 \\
\addlinespace
60R\_60V\_02I\_05M & 5 & 0 \\
60R\_60V\_02I\_06M & 5 & 0 \\
60R\_60V\_04I\_05M & 2 & 0 \\
60R\_60V\_04I\_06M & 4 & 1 \\
\midrule
Sum & 34 & 1 \\
\bottomrule
\end{tabular}
}
\end{minipage}
\hfill
\begin{minipage}{0.48\textwidth}
\centering
\okC{
\text{Multi-floor instances}
\label{tab:overview_mslp_type_1_larger_floor}
\begin{tabular}{lrr}
\toprule
Group & All Feas. & [1,9] Feas. \\
\cmidrule{1-1} \cmidrule{2-3}
40R\_40V\_02F\_05M & 5 & 0 \\
40R\_40V\_02F\_06M & 5 & 0 \\
40R\_40V\_04F\_05M & 4 & 0 \\
40R\_40V\_04F\_06M & 5 & 0 \\
\addlinespace
60R\_60V\_02F\_05M & 5 & 0 \\
60R\_60V\_02F\_06M & 5 & 0 \\
60R\_60V\_04F\_05M & 2 & 0 \\
60R\_60V\_04F\_06M & 3 & 1 \\
\midrule
Sum & 34 & 1 \\
\bottomrule
\end{tabular}
}
\end{minipage}
\end{table}

\okC{Tables \ref{tab:mean_summary_mslp_results_larger_instances_multi_island} and \ref{tab:mean_summary_mslp_results_larger_instances_multi_floor} summarize the results for the instances with 40 and 60 requests.
The column names are equal to the column names of the Tables \ref{tab:mean_summary_mslp_results_multi_island} and \ref{tab:mean_summary_results_multi_floor} earlier in this section.
Row \mbox{``Mean --- Feas."} displays the mean of the columns, excluding the groups with an ``N/A" value, i.e., groups with which the heuristic failed to find a feasible solution in the ten independent executions.
Row \mbox{``Mean --- All"} shows the mean of the columns for Type~2 instances.
Observe that the MSLP heuristic spent much more CPU time on Type~2 instances of both families compared to the CPU time spent on Type~1 instances, which supports the results achieved on the smaller instances (see Tables~\ref{tab:mean_summary_mslp_results_multi_island} and \ref{tab:mean_summary_mslp_results_multi_floor}).
Furthermore, the CPU times to best for all instance groups are still significantly lower compared to the total CPU time spent.
Consequently, the number of iterations set for the heuristic execution may be reduced without compromising solution quality.}

\begin{table}[!htp]\centering
\caption{\okC{Mean summary of the MSLP results for the larger instances of the multi-island family}}
    \label{tab:mean_summary_mslp_results_larger_instances_multi_island}
    %\resizebox{\textwidth}{!}{ % use this if the table is too large
    \scriptsize
    \okC{
    \begin{tabular}{lrrrrrrrrrrr}\toprule
        & &\multicolumn{4}{c}{Type~1} & &\multicolumn{4}{c}{Type~2} \\
        \cmidrule{3-6}\cmidrule{8-11}
        Group &\hspace{0.05cm} &Min. Sol. &Mean Sol. &Time (s) &TTB (s) &\hspace{0.05cm} &Min. Sol. &Mean Sol. &Time (s) &TTB (s) \\
        \cmidrule{1-1} \cmidrule{3-6} \cmidrule{8-11} 
        40R\_40V\_02I\_05M &\hspace{0.1cm} &6148.55 &6227.75 &219.60 &91.37 &\hspace{0.2cm} &6172.13 &6363.94 &1721.57 &891.91 \\
        40R\_40V\_02I\_06M &\hspace{0.1cm} &5914.80 &6085.02 &319.31 &154.66 &\hspace{0.2cm} &6156.67 &6354.61 &1846.12 &908.77 \\
        40R\_40V\_04I\_05M &\hspace{0.1cm} &N/A &N/A &N/A &N/A &\hspace{0.2cm} &6330.37 &6619.39 &2236.89 &1117.10 \\
        40R\_40V\_04I\_06M &\hspace{0.1cm} &5238.27 &5423.64 &388.05 &190.86 &\hspace{0.2cm} &6419.31 &6592.43 &2398.56 &1306.78 \\
        \addlinespace
        60R\_60V\_02I\_05M &\hspace{0.1cm} &8197.12 &8252.74 &574.35 &239.30 &\hspace{0.2cm} &9498.10 &9855.85 &3010.72 &1615.88 \\
        60R\_60V\_02I\_06M &\hspace{0.1cm} &8184.79 &8300.70 &673.12 &236.56 &\hspace{0.2cm} &9460.83 &9845.82 &3129.14 &1640.45 \\
        60R\_60V\_04I\_05M &\hspace{0.1cm} &N/A &N/A &N/A &N/A &\hspace{0.2cm} &9878.48 &10245.50 &3343.94 &1786.05 \\
        60R\_60V\_04I\_06M &\hspace{0.1cm} &N/A &N/A &N/A &N/A &\hspace{0.2cm} &9936.92 &10166.36 &3548.44 &1947.87 \\
        \midrule
        Mean --- Feas. &\hspace{0.1cm} &6736.70 &6857.97 &434.89 &182.55 &\hspace{0.2cm} &7541.41 &7802.53 &2421.22 &1272.76 \\
        Mean --- All &\hspace{0.1cm} &N/A &N/A &N/A &N/A &\hspace{0.2cm} &7981.60 &8255.49 &2654.42 &1401.85 \\
        \bottomrule
    \end{tabular}
    }
\end{table}

\begin{table}[!htp]\centering
\caption{\okC{Mean summary of the MSLP results for the larger instances of the multi-floor family}}
    \label{tab:mean_summary_mslp_results_larger_instances_multi_floor}
    %\resizebox{\textwidth}{!}{ % use this if the table is too large
    \scriptsize
    \okC{
    \begin{tabular}{lrrrrrrrrrrr}\toprule
        & &\multicolumn{4}{c}{Type~1} & &\multicolumn{4}{c}{Type~2} \\
        \cmidrule{3-6}\cmidrule{8-11}
        Group &\hspace{0.05cm} &Min. Sol. &Mean Sol. &Time (s) &TTB (s) &\hspace{0.05cm} &Min. Sol. &Mean Sol. &Time (s) &TTB (s) \\
        \cmidrule{1-1} \cmidrule{3-6} \cmidrule{8-11} 
        40R\_40V\_02F\_05M &\hspace{0.05cm} &6690.01 &6864.21 &543.26 &276.43 &\hspace{0.1cm} &6854.00 &7075.67 &1941.55 &973.00 \\
        40R\_40V\_02F\_06M &\hspace{0.05cm} &6698.24 &6764.26 &752.31 &423.67 &\hspace{0.1cm} &6724.06 &7057.67 &2085.85 &948.25 \\
        40R\_40V\_04F\_05M &\hspace{0.05cm} &N/A &N/A &N/A &N/A &\hspace{0.1cm} &8117.65 &8280.66 &2378.25 &1160.92 \\
        40R\_40V\_04F\_06M &\hspace{0.05cm} &8935.57 &9103.35 &418.77 &163.32 &\hspace{0.1cm} &8029.69 &8225.96 &2526.59 &998.50 \\
        \addlinespace
        60R\_60V\_02F\_05M &\hspace{0.05cm} &9983.71 &10268.51 &1013.04 &471.13 &\hspace{0.1cm} &10819.28 &11148.64 &3219.82 &1677.48 \\
        60R\_60V\_02F\_06M &\hspace{0.05cm} &9895.43 &10065.90 &1352.49 &730.80 &\hspace{0.1cm} &10867.11 &11157.69 &3335.26 &1701.31 \\
        60R\_60V\_04F\_05M &\hspace{0.05cm} &N/A &N/A &N/A &N/A &\hspace{0.1cm} &12664.95 &12957.25 &3376.01 &1755.50 \\
        60R\_60V\_04F\_06M &\hspace{0.05cm} &N/A &N/A &N/A &N/A &\hspace{0.1cm} &12689.80 &12880.67 &3525.98 &1706.96 \\
        \midrule
        Mean --- Feas. &\hspace{0.05cm} &8440.59 &8613.25 &815.97 &413.07 &\hspace{0.1cm} &8658.83 &8933.13 &2621.82 &1259.71 \\
        Mean --- All &\hspace{0.05cm} &N/A &N/A &N/A &N/A &\hspace{0.1cm} &9595.82 &9848.03 &2798.66 &1365.24 \\
        \bottomrule
    \end{tabular}
    }
\end{table}

\begin{figure}[!ht]
    \centering
    \subfigure[Multi-island instances]{
        \includegraphics[width=0.45\linewidth]{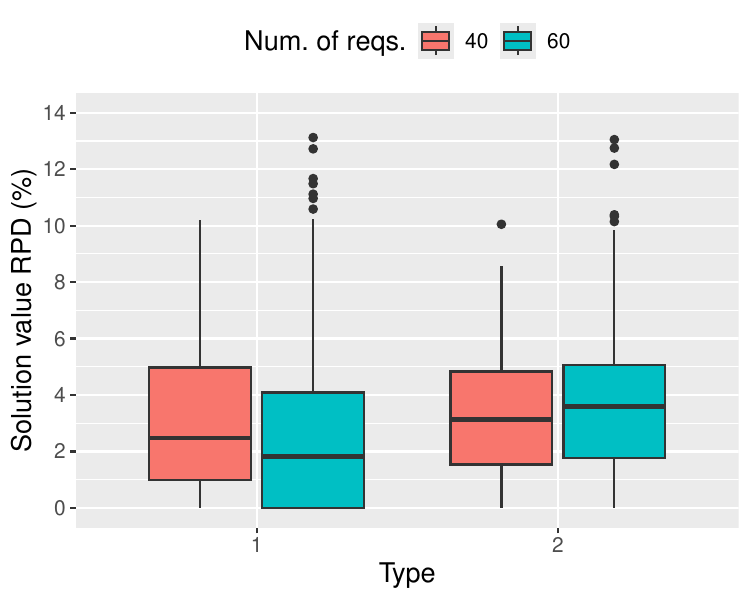}
    }
    \subfigure[Multi-floor instances]{
        \includegraphics[width=0.45\linewidth]{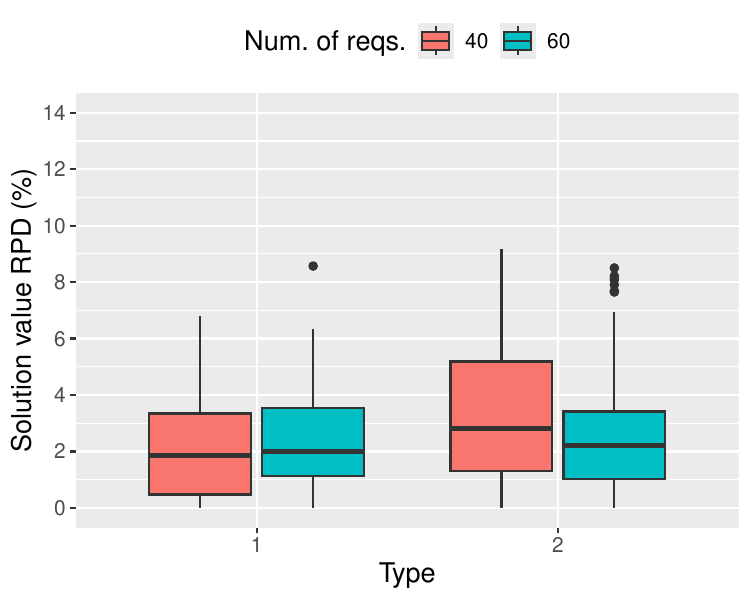}
    }
    \caption{\okC{Box plots of the solution value RPDs by instance type for the solution value obtained in an independent execution compared to the best solution value achieved in the ten independent executions. Observations only for instances with 40 and 60 requests, and independent executions that found a feasible solution.}}
    \label{fig:dev_mslp_larger_boxplots}
\end{figure}

\okC{
Regarding the quality of the heuristic solutions for the instances with 40 and 60 requests, the box plots of Figure \ref{fig:dev_mslp_larger_boxplots} depict the solution value RPDs by instance type for the solution value obtained in an independent execution compared to the best solution value achieved in the ten independent executions.
There are 80 observations in the Type~2 box plots of both families.
However, 52 and 57 out of the 400 (13.0\% and 14.3\%) observations are missing in the Type~1 multi-island and multi-floor box plots, respectively. 
Among them, 20 and 10 correspond to instances with 40 requests, while the remaining ones correspond to 60-request instances.
Disregarding the outliers, the solutions obtained are generally at most 10\% worse than the best solution obtained by the heuristic for a given instance. 
Moreover, about 75\% of the solution value RPD observations are below 5\%.
}

\okC{
In the appendices, we provide some useful insights about the MSLP heuristic.
First,~\ref{app:mslp_solution_feasibility} discusses the solution feasibility throughout the iterations of the MSLP heuristic. 
In summary, although both the number of requests and the number of regions reduced the number of feasible solutions obtained, the MSLP heuristic still found feasible solutions for most of the iterations, except for Type~2 instances with 40 and 60 requests.
Moreover, for these larger instances, it might be useful to provide an additional machine to increase the chances of finding feasible solutions. 
The heuristic may also be used to verify whether the number of machines is sufficient to generate high-quality solutions.
Second, \ref{app:lps_relevance_mslp} displays the \okC{LP improvement procedure's relevance} in the MSLP heuristic to improve the value of the solutions found by the (semi-)greedy heuristic.
Even for the larger instances, the mean LP improvement percentage is generally above 5\%, reaching a maximum improvement of almost 43\% on some 6-request instances.
Besides, we show that Type~2 instances benefited more from this LP procedure compared to Type~1 instances.
}

\subsection{\okC{Summary of the characteristics of the solutions}}
\label{subsec:sol_chars_summary}

\okC{
In this section, we provide a summary of \ref{app:additional_machine_impact} and \ref{app:solution_chars}. 
\ref{app:additional_machine_impact} analyzes the impact of the additional machine on the solution value of the instances.
In essence, the results indicate that providing an additional machine contributed to the solution value reduction in most of the instances analyzed.
Besides that, the Type~1 instances benefited more from the additional machine, especially the larger ones, reaching a maximum improvement of over 10\% in a 40-request multi-floor instance.
}

\okC{
\ref{app:solution_chars} shares some additional details regarding the characteristics of the solutions obtained by the MIP with valid inequalities and the MSLP heuristic.
In brief, the appendix shows that the Type~1 instance solutions made more use of the vehicles without using much of their capacity.
Therefore, the requests in Type~1 instance solutions are more evenly distributed among the vehicles due to the instances' short-scheduling horizons.
Additionally, the percentage of vehicles used does not increase proportionally as the number of requests increases for both Type~1 and Type~2 instances.
The appendix also notes that providing an additional machine decreases the average machine usage time, even though the percentage of machines used generally didn't increase.
Clearly, machine usage times are generally greater on Type~1 instances compared to Type~2 instances.
Besides, almost all larger instance solutions used all the machines available, although the machine usage time was not greater.
Finally, the appendix shows that increasing the number of regions increases the machine usage time on average.
}

\section{Concluding remarks}
\label{sec:finalremarks}

This paper introduced and formalized the pickup and delivery problem with time windows and scheduling on the edges (PDPTW-SE).
The problem extends the PDPTW by requiring that some edges be traversed by special machines that must be scheduled.
The PDPTW-SE enables a variety of potential applications.
Firstly, it offers avenues for developing novel system designs, such as logistical planning for pickups and deliveries within multiple islands, as well as in areas affected by large-scale disasters. 
Secondly, it can enhance the tactical and operational decision-making processes within existing system designs, such as in the logistical planning of pickups and deliveries within multi-floor hospitals using automated guided vehicles (AGVs). 

We devised a mixed integer programming (MIP) formulation\okC{, a MIP preprocessing step, valid inequalities,} and a multi-start heuristic with an LP improvement procedure (MSLP).
Furthermore, we proposed a benchmark set adapted from widely used PDPTW instances.
The proposed benchmark set comprises two families of instances representing the previously mentioned applications.
More specifically, the first family addresses the multi-island scenario, where requests are spread across different islands, and a limited number of cargo ships can be scheduled to transport the vehicles between the islands.
The second family addresses the multi-floor scenario, in which the requests are spread across different building floors, such as hospitals or hotels. In the latter, a limited number of dedicated elevators can be scheduled to transport AGVs between the floors. In hospitals, AGVs can replace the labor-intensive, repetitive, low-value-added, and mundane rule-based tasks such as delivery and replacement of linen, food services to the wards during mealtimes, delivery of consumables or hazardous waste, among others.

\okC{We conducted computational experiments using the set of smaller instances with up to 12 requests to analyze the performance of the solver on the MIP formulation.
The results indicated that the solver on the formulation could solve to optimality instances with up to 12 requests.
Indeed, the solver was much more effective at proving optimality for instances with short-scheduling horizons than for those with long ones.
For short-scheduling horizon instances, the solver proved optimality in 61.2\% and 72.5\% of the instances with up to 12 requests in the multi-island and multi-floor families, respectively.
For long-scheduling horizon instances, these percentages decreased to 30.0\% and 27.5\%, respectively.
Besides, the solver achieved significantly lower gaps on short-scheduling horizon instances compared to those for long-scheduling horizon instances.
It is noteworthy, however, that the solver on the MIP formulation found feasible solutions for all \okC{short-horizon} instances, and for 94.1\% of the instances in the whole benchmark set.
We also slightly improved the MIP formulation by adding several valid inequalities.
With this new formulation, the solver achieved three more optimal solutions and feasible solutions for three more instances.
Consequently, the MIP with valid inequalities solved 48.8\% of the 320 instances with up to 12 requests to optimality, while feasible solutions were found for 95.0\% of them, including all \okC{short-horizon} instances as in the baseline MIP.
Moreover, among the instances where the solver found a feasible solution but reached the time limit, the solver on the MIP with valid inequalities achieved equal or smaller optimality gaps in at least 75\% of them compared to the baseline MIP, especially \okC{short-horizon} instances. Additionally, for instances where a feasible solution was found, the strengthened formulation also improved several \okC{long-horizon} instance solution values.}

\okC{We also carried out experiments to analyze the performance of the MSLP heuristic using the whole benchmark set.
This approach found feasible solutions in all independent executions for the 320 instances with up to 12 requests of the benchmark set.
Additionally, the heuristic found solutions with reasonably low deviations and short times when compared to the solver on the MIP with valid inequalities.
The heuristic could also find several solutions that at least match or surpass the solutions of the solver on the MIP with valid inequalities for instances in both families, particularly in \okC{long-horizon} instances with 12 requests.
For the larger instances with 40 and 60 requests, the MSLP heuristic found feasible solutions for 150 out of the 160 (93.8\%) instances with 40 and 60 requests, whereas the solver on MIP formulations couldn't find a feasible solution for any of them.
It is clear from the results that providing an additional machine increased the chances of finding a feasible solution for these instances, especially for instances with four regions.
Besides, this additional machine contributed to the solution value reduction in most of the instances analyzed, in particular \okC{short-horizon} instances with a larger number of requests.
Regarding solution quality, about 75\% of the solution values obtained by the heuristic for these instances are at most 5\% worse than the best solution value obtained.
Furthermore, the results indicated that the MSLP heuristic may also be used to verify whether the number of machines is sufficient to generate high-quality solutions.
It is also shown that the LP improvement procedure played an important role in improving the solution values, even for the larger instances.
These improvements were generally above 5.0\% on average, reaching a maximum improvement of 43.0\% on average, and \okC{long-horizon} instances benefited more from this procedure compared to \okC{short-horizon} instances.}

\okC{When comparing the results of multi-island and multi-floor instances, there's no significant difference regarding mean gap, mean time, and other solution characteristics for both methods.
It was observed that long scheduling horizons allow more vehicle capacity utilization and fewer vehicles in both families.
Additionally, offering a new machine option can reduce machine usage on average, even though the number of machines used does not necessarily increase.
Finally, as expected, increasing the number of regions and shortening the scheduling horizon implies more machine usage time on average.}

\vspace{0.5cm}

{
\footnotesize
% Acknowledgments here
%\section*{Acknowledgements}

\noindent \textbf{Acknowledgements:} This study was financed in part by the Coordenação de Aperfeiçoamento de Pessoal de Nível Superior - Brasil (CAPES) - Finance Code 001.
This work was partially supported by the Brazilian National Council for Scientific and Technological Development (CNPq), grants 314662/2020-0, 314718/2023-0, \okC{and 445324/2024-4. This work was partially} supported by the FAPESB INCITE PIE0002/2022 grant.
} 
 
% Leave this (end of acknowledgment)

\bibliography{main}

\newpage
\appendix

\makeatletter
\renewcommand{\p@subfigure}{\thefigure} % important to add a prefix for figure reference in text.
\makeatother

% \input{appendix_01}
% \newpage
\section{Summary of the notation}
\label{app:notations}

\begin{table}[!ht]
\caption{Parameters}
\label{tab:parameters}
\renewcommand{\arraystretch}{1.3}
\scriptsize
\centering
\begin{tabular}{p{1.5cm} p{9cm}}
    \hline
    \multicolumn{1}{c}{Notation}  & \multicolumn{1}{c}{Definition} \\
    \hline 
    $V$                                     & Set of nodes in the undirected graph $G$, $V=\{0, \ldots, 2n\}$ \\
    $V_p$                                   & Set of pickup nodes, $V_p=\{1,\ldots,n\}$ \\
    $V_d$                                   & Set of delivery nodes, $V_d=\{n+1,\ldots,2n\}$ \\
    $E$                                     & Set of edges in the undirected graph $G$, $E = E^s \cup E^m$ \\
    $E^s$                                   & Set of edges traversed directly by vehicles \\
    $E^m$                                   & Set of edges whose traversal involves using a machine \\
    $R$                                     & Set of unsplittable requests, $R=\{r_1, \ldots, r_n\}$ \\
    $K$                                     & Set of heterogeneous vehicles \\
    $F$                                     & Set of vehicle accessible regions \\
    $H$                                     & Set of machines \\
    $H_{ij}$                                & Set of machines that can be used to traverse edge $e \in E^m, \ e=\{i,j\}$ \\
    $F_h$                                   & Set of vehicle accessible regions that machine $h$ has an instance \\
    $V'$                                    & Set of nodes $V' = V \cup \{2n+1\}$ \\
    $A$                                     & Set of possible arcs \\
    $A^s$                                   & Set of possible arcs, corresponding to $E^s$ \\
    $A^m$                                   & Set of possible arcs, corresponding to $E^m$ \\
    $n$                                     & Number of requests \\
    $r_i$                                   & Request $i$, $r_i~=~\langle v_i, v_{n+i}, q_i \rangle$ \\
    $q_i$                                   & Weight of request $i$  \\
    $e_i$                                   & Time window's earliest time at node $i$ \\
    $l_i$                                   & Time window's latest time at node $i$ \\
    $s_i$                                   & Service time at node $i$ \\
    $Q_k$                                   & Capacity of vehicle $k$ \\
    $\hat{d}^k_{ij}$                        & Time to traverse the edge $ij \in E$\\
    $z$                                     & Number of regions \\
    $f^h_i$                                 & Node of machine $h$ corresponding to node $i$ \\
    $\bar{d}^k_{ih}$                    & Time to travel from node $i$ to station $f^h_i$ of machine $h$ with vehicle $k$ \\
    $O^h_{f^h_i f^h_j}$                     & Time to travel from node $f^h_i$ to node $f^h_j$ with machine $h$ \\
    $d^k_{ij}$                              & Time to traverse arc $(i,j)$ with vehicle $k$ \vspace{2pt}\\
    \hline 
\end{tabular}
\end{table}

\section{Possible big-M values for the MIP formulation}
\label{app:big_m_bounds}

As seen in Section \ref{sec:formulation}, there are some big-M values that must be defined.
Possible values are presented as follows:

\begin{itemize}
    \item $M_1$: max\{$Q_k \, | \,  k \in K\} + \text{max}\{q_i \, | \, i \in V_p \}$ + 1.
    \item $M_2$: $l_0$ + max\{$s_i \ | \ i \in V$\} + max\{$d^k_{ij} \ | \ i \in V', j \in V', k \in K$\} + 1.
    \item $M_3$: $l_0$ + max\{$d^k_{ij} \ | \ i \in V', j \in V', k \in K$\} + 1.
    \item $M_4$: $l_0$ + max\{$s_i \ |\ i \in V$\} + max\{$\bar{d}^k_{i h} \ | \ i \in V', h \in H, k \in K$\} + 1.
    \item $M_5$: $l_0$ + max\{$\bar{d}^k_{i h}\ |\ i \in V', h \in H, k \in K$\} + 1.
    \item $M_6$: $l_0$ + max\{$O^h_{f^h_i f^h_j}\ |\ h \in H, i \in V', j \in V'$\} + max\{$\bar{d}^k_{i h}\ |\ i \in V', h \in H, k \in K$\} + 1.
    \item $M_7$: $l_0$ + $2\times$max\{$O^h_{f^h_i f^h_j}\ |\ h \in H, i \in V', j \in V'$\} + 1.
    \item $M_8$: max\{$O^h_{f^h_i f^h_j}\ | \ h \in H, i \in V', j \in V'$\} + 1.
\end{itemize}
% \newpage
\section{Multi-start heuristic auxiliary algorithms and details}
\label{app:heur_details}

In this section, we describe the remaining algorithms invoked in the MSLP heuristic and other details.
\ref{app:cand_list_and_choice} describes the GET\_INSERTION\_CANDIDATE\_LIST called in line \ref{alg:semi_ins_heur:050} of Algorithm \ref{alg:semi_ins_heur} and the CHOOSE\_CANDIDATE\_BY\_QUALITY invoked in line \ref{alg:semi_ins_heur:060} of Algorithm \ref{alg:semi_ins_heur}.
\ref{app:feasibility_analysis} states the assumptions and introduces the concepts needed to analyze insertion feasibility.
\ref{app:vehicle_travel_time} describes how we calculate vehicle travel times during the analysis of insertion feasibility.
\ref{app:update_solution} presents how the solution is updated in both greedy and semi-greedy insertion heuristics.

\subsection{Definition of the candidate list and element selection}
\label{app:cand_list_and_choice}

Algorithm \ref{alg:ins_cand_list} implements the GET\_INSERTION\_CANDIDATE\_LIST called in line \ref{alg:semi_ins_heur:050} of Algorithm~\ref{alg:semi_ins_heur} by adapting the GET\_BEST\_INSERTION\_CANDIDATE described in Algorithm \ref{alg:best_ins_cand} to get all the feasible candidates in a list $CL$ (line \ref{alg:ins_cand_list:050}), instead of getting only the best one.
After analyzing all candidates, we return $CL$ (line~\ref{alg:ins_cand_list:060}).

\begin{algorithm}[H]
    \scriptsize
    \caption{\footnotesize GET\_INSERTION\_CANDIDATE\_LIST($p$, $d$, $S$, $I$)}
    \label{alg:ins_cand_list}
    $CL \gets \emptyset$\; \label{alg:ins_cand_list:010}
    \For{$k \in K$ \textbf{and} $p_{pos} \in 2:S.vehicles[k].length$ \textbf{and} $d_{pos} \in p_{pos}:S.vehicles[k].length$}{ \label{alg:ins_cand_list:020}
        $\mathcal{C} \gets$ ANALYZE\_INSERTION\_FEASIBILITY($k$, $p_{pos}$, $d_{pos}$, $p$, $d$, $S$, $I$)\; \label{alg:ins_cand_list:030}
        \If{$\mathcal{C}$ is feasible}{ \label{alg:ins_cand_list:040}
            Add $\mathcal{C}$ to CL\; \label{alg:ins_cand_list:050}
        }
    }
    \Return $CL$\;\label{alg:ins_cand_list:060}
\end{algorithm}

Algorithm \ref{alg:cand_by_qlty} implements the CHOOSE\_CANDIDATE\_BY\_QUALITY invoked in line \ref{alg:semi_ins_heur:060} of Algorithm \ref{alg:semi_ins_heur} by describing how a candidate is chosen from the candidate list $CL$ using the threshold parameter $\alpha$.
Firstly, in lines \ref{alg:cand_by_qlty:010} and \ref{alg:cand_by_qlty:020}, if the candidate list $CL$ is empty, we return an infeasible candidate $\mathcal{C}$.
Then, in lines \ref{alg:cand_by_qlty:030} and \ref{alg:cand_by_qlty:040}, we get the minimum cost $c_{min}$ and maximum cost $c_{max}$ of the candidates in $CL$.
In line \ref{alg:cand_by_qlty:050}, a restricted candidate list $RCL$ is created from $CL$ by inserting candidates whose costs are at most $\alpha$ worse than $c_{min}$.
In line \ref{alg:cand_by_qlty:060}, a random candidate $\mathcal{C}'$ is chosen from the restricted candidate list $RCL$.
In line \ref{alg:cand_by_qlty:070}, the chosen candidate $\mathcal{C}'$ is returned.

\begin{algorithm}[H]
\scriptsize
\caption{\footnotesize CHOOSE\_CANDIDATE\_BY\_QUALITY($CL$, $\alpha$)}\label{alg:cand_by_qlty}
    \If{$|CL|\text{ is empty}$}{ \label{alg:cand_by_qlty:010}
        \Return Infeasible $\mathcal{C}$\; \label{alg:cand_by_qlty:020}
    }
    $c_{min} \gets min(\mathcal{C}.cost : \mathcal{C} \in CL)$\; \label{alg:cand_by_qlty:030}
    $c_{max} \gets max(\mathcal{C}.cost : \mathcal{C} \in CL)$\; \label{alg:cand_by_qlty:040}
    $RCL \gets \{\mathcal{C} \in CL \ | \ \mathcal{C}.cost \leq c_{min} + \alpha (c_{max} - c_{min})\}$\; \label{alg:cand_by_qlty:050}
    $\mathcal{C}' \gets$ $RCL[\ rand() \ \textbf{mod} \ |RCL| \ ]$\; \label{alg:cand_by_qlty:060}
    \Return $\mathcal{C}'$\; \label{alg:cand_by_qlty:070}
\end{algorithm}

\subsection{Insertion feasibility analysis}
\label{app:feasibility_analysis}

To analyze the feasibility of inserting nodes $p$ and $d$ at positions $p_{pos}$ and $d_{pos}$ of the route $S.vehicles[k], \ k \in K$, it is necessary to check both time and capacity constraints.
Let $(i_0, i_1, i_2, \ldots, i_{\eta})$, $i_0 = 0, \ i_{\eta} = 2n+1$, be a partially constructed feasible route for a vehicle $k \in K$, for which the service start times and the vehicle loads are known.
Besides that, consider that the pickup node is inserted between $i_{g-1}$ and $i_{g}$, and the delivery node is inserted between $i_{h-1}$ and $i_{h}$, $1 \leq g \leq h \leq \eta$.
Note that when $g = h$, the delivery node $d$ is inserted after the pickup node $p$.
For vehicle's capacity feasibility, the insertion of the pickup node $p$ could exceed vehicle's capacity in the node $p$ or in some node $i_r, \ g \leq r \leq h-1$. 
If $\text{max}_{g \leq r \leq h-1}\{z_{i_r}\} + q_p \leq Q_k$, this insertion is feasible.
For time feasibility, we adapt the concepts described for the VRPTW in \cite{solomon87} to the PDPTW-SE.

Denote by $t^{new}_{i_g}$ and $t^{new}_{i_h}$ the new service start times at nodes $i_{g}$ and $i_{h}$, respectively, given the insertion of nodes $p$ and $d$ in the route of the vehicle $k$.
Also, let $w_{i_r}$ be the waiting time at $i_r$ for $p \leq r \leq \eta$.
As we assumed that the traversal times respect triangular inequalities (see Section~\ref{sec:formalization}), each insertion defines a \textit{push forward} in the schedule at $i_g$ and $i_h$:
$\text{PF}_{i_g} = t^{new}_{i_g} - t_{i_g} \geq 0$, and $\text{PF}_{i_h} = t^{new}_{i_h} - t_{i_h} \geq 0$.
Supposing no machine travel is modified later on, we conclude that $\text{PF}_{i_{r+1}} = \text{max}\{0, \text{PF}_{i_r} - w_{i_{r+1}}\} \text{ for } \ g \leq r \leq h-2 \text{ and } h \leq r \leq \eta-1$.
Therefore, if $\text{PF}_{i_g} > 0$ or $\text{PF}_{i_h} > 0$, serving some nodes $i_r, \ g \leq r \leq \eta,$ could become infeasible.
Besides that, for some node $i_r, \ h \leq r < \eta,$ the service start time can remain unchanged, i.e., $\text{PF}_{i_r} = 0$.
Therefore, the vehicle will \okC{serve} each node in the interval $(i_r, i_{\eta}]$ at the same previous scheduling time.
However, each $\text{PF}_{i_r}, \ g \leq r \leq \eta$ also defines a \textit{push forward} for the next machine travel.
If the machine travel is pushed forward, an ejection chain procedure would be necessary because it could also push forward other machine travels and, consequently, other vehicle service start times.
To avoid an ejection chain, we disregard all machine travels held with the vehicle $k$ from $p_{pos}$ onward, i.e., machine travels are rescheduled for each analysis.
Furthermore, we ensure that machine travels held with other vehicles remains unchanged.
This simplification can lead to a different scenario from \cite{solomon87}, described in the following.

Let $k_1$ be a vehicle with a partially constructed route with some machine travels scheduled, which we denote as $mtrvs_{k_1}$.
Let $k_2$ be another vehicle in which a request insertion is being analyzed.
Suppose a machine travel $mtrv_2$ is necessary at some point of the route of $k_2$, and the best moment to insert $mtrv_2$ competes with a machine travel $mtrv_1 \in mtrvs_{k_1}$.
As $mtrv_1$ is not allowed to be modified, $mtrv_2$ is allocated at a different moment or on another machine.
However, the next request insertion is placed in the route of vehicle $k_1$ and $mtrv_1$ is rescheduled, opening  space for $mtrv_2$ to be rescheduled in a better position.
If the next request insertion occurs in the route of $k_2$ and $mtrv_2$ is disregarded, then $mtrv_2$ can be rescheduled to a better position, reducing $mtrv_2.\Delta t$.
Thus, the new service start time at node $i_r, \ h < r \leq \eta$, can be pushed backward, i.e., $\text{PF}_{i_r} < 0$.
If $\text{PF}_{i_{\eta}} < 0$, then the completion time of vehicle $k_2$ is reduced.

Therefore, the ANALYZE\_INSERTION\_FEASIBILITY, called in line \ref{alg:best_ins_cand:040} of Algorithm~\ref{alg:best_ins_cand} and line \ref{alg:ins_cand_list:030} of Algorithm \ref{alg:ins_cand_list}, should examine these nodes sequentially for time and capacity feasibility. The analysis continues until it finds: some node $i_r, \ g \leq r < h$, with infeasible capacity; or some node $i_r, \ g \leq r < \eta$, node $p$, or node $d$ with infeasible time.
In the worst case, all the nodes $i_r, \ p \leq r \leq \eta$ are examined, and the insertion candidate is feasible.
Besides that, it is also necessary to save all machine travel data associated with this insertion candidate, which is done during the insertion feasibility analysis.

\subsection{Vehicle travel time}
\label{app:vehicle_travel_time}

Note that, in PDPTW-SE, there are two possibilities to determine the start time for the service at a node $i \in V_p \cup V_d$.
For intra-region nodes, the start time at node $j \in V' \setminus \{0\}$ can be calculated directly, i.e., $t_j = \text{max}\{e_j, t_i + s_i + d^k_{ij}\}, \ \forall \ (i,j) \in A^s, \ k \in K$.
Conversely, for inter-region nodes, the start time at node $j \in V' \setminus \{0\}$ involves the usage of a machine to traverse the arc $(i,j) \in A^m$.
It can be calculated as $t_j = \text{max}\{e_j, t_i + s_i + \bar{d}^k_{ih} + w^h_{ij} + O^h_{f^h_i f^h_j} + \bar{d}^k_{jh}\}, \ (i,j) \in A^m, \ k \in K, \ h \in H_{ij}$, where $w^h_{ij}$ is the waiting time at machine station $f^h_i$ when traversing arc $(i,j) \in A^m$.
The waiting time is calculated as follows: $w^h_{ij} = \text{max}(0, \alpha^h_{i'j'} + O^h_{f^h_{i'} f^h_{j'}} - k_{arr}), \ (i',j') \in A^m$, where $k_{arr} = t_i + s_i + \bar{d}^k_{ih}$ and~$\gamma^h_{i'j'ij} = 1$.

Suppose that during the feasibility analysis we examine the node $j = i_r, \ g-1 < r \leq m$, which succeeds node $i, \ (i,j) \in A^m$.
By construction, the node $i$ can be $i_{r-1}$, $p$, or $d$, and it was already examined, i.e., we know its new service start time.
To get the best inter-route arc traversal time, we iteratively search the best position to insert a machine travel in $S.machines[h], \ h \in H_{ij}$, without modifying any other machine travel's start time.
Denote by $best\_mtrv$ this best machine travel.
Given that we know the new machine travels previously calculated in this analysis, which we denote as $mtrvs$, this search can start after the last $mtrv \in mtrvs[h]$.
Besides that, we can stop the analysis at a machine sequence of travels $S.machines[h], \ h \in H_{ij}$, at the first fit position, or if the new machine travel insertion leads to violating the latest start time $l_j$.
We can do this because inserting a new machine travel after one of these positions will delay even more the arrival time at the next visit.
If a feasible machine travel was found, we calculate $mtrv.\Delta t = s_i + \bar{d}^k_{ih} + w^h_{ij} + O^h_{f^h_i f^h_j} + \bar{d}^k_{jh}$ and update $best\_mtrv$ if $mtrv.\Delta t < best\_mtrv.\Delta t$.
After analyzing each machine $h \in H_{ij}$, $best\_mtrv$ is added to $mtrvs$.
If no feasible machine travel was found, the request insertion is infeasible.
Thus, we set $best\_mtrv.\Delta t = \infty$ and return $best\_mtrv.\Delta t$ as the best inter-route arc traversal time.

\subsection{Solution update}
\label{app:update_solution}

The UPDATE\_SOLUTION algorithm invoked in Algorithms \ref{alg:ins_heur} and \ref{alg:semi_ins_heur} is similar to the ANALYZE\_{\linebreak}INSERTION\_FEASIBILITY algorithm described in \ref{app:feasibility_analysis}.
Given the chosen candidate $\mathcal{C}$ and the nodes $p$ and $d$ to insert in solution $S$, we iteratively update $visit.st$ and $visit.load$ for every $visit \in S.vehicles[\mathcal{C}.k]$, starting from $\mathcal{C}.p_{pos}$.
Then, we insert the nodes $p$ and $d$ at $S.vehicles[\mathcal{C}.k]$ in positions $\mathcal{C}.p_{pos}$ and $\mathcal{C}.d_{pos}$, respectively.
We also insert each $mtrv \in \mathcal{C}.mtrvs$ in its respective position $h_{pos}$ in $S.machines[mtrv.h], \ h \in H$.
Finally, we remove the machine travels disregarded (see \ref{app:feasibility_analysis}) and update $S.cost$ by adding $\mathcal{C}.cost$.

% \newpage
%\newpage
\section{Procedure to ensure instance feasibility}
\label{app:alg_to_ensure_feas}

In this section, we adapt the greedy constructive heuristic described in Section \ref{sec:greedy_insertion_heuristic} to ensure feasibility by allowing modifications in the input instance characteristics.
In \ref{app:gen_structure_inst_feas}, we present the general structure of the algorithm to ensure feasibility.
In \ref{app:ins_feas_cost}, we describe the algorithm for the feasibility analysis and the calculus of the request insertion cost.
In \ref{app:vehicle_travel_time_mod}, we present how we calculate vehicle travel times during the insertion analysis.
In \ref{app:instance_modification}, we present how the input instance characteristics are modified.

\subsection{General structure}
\label{app:gen_structure_inst_feas}

Algorithm~\ref{alg:ins_heur_mod} details the procedure to ensure feasibility. Its main lines are highlighted in what follows.
In line~\ref{alg:ins_heur_mod:050}, the algorithm searches the best feasible insertion candidate $\mathcal{C}$ and the best infeasible insertion candidate $I\mathcal{C}$ for the request nodes $p$ and $d$ in the partially constructed solution $S$.
If the insertion candidate $\mathcal{C}$ is not feasible, we update $S$ by inserting the request nodes $p$ and $d$ in solution $S$ following the best infeasible insertion candidate $I\mathcal{C}$ data (line~\ref{alg:ins_heur_mod:090}).
Besides that, we modify the input instance $I$, which makes $S$ infeasible (line~\ref{alg:ins_heur_mod:100}).
After inserting all requests, if $S$ is not feasible (line~\ref{alg:ins_heur_mod:105}), then the instance $I$ was modified.
Thus, we save the modified instance $I$ as a new instance (line~\ref{alg:ins_heur_mod:107}).

\begin{algorithm}[H]
\scriptsize
\caption{\footnotesize GREEDY\_INSERTION\_HEURISTIC\_TO\_ENSURE\_FEASIBILITY($I$)}
\label{alg:ins_heur_mod}
    $S \gets $ Empty solution $S$ with $S.feasible=\textbf{true}$ and no request inserted\; \label{alg:ins_heur_mod:010}
    $\mathcal{L} \gets$ List of pickup nodes sorted in nondecreasing order of time window width\; \label{alg:ins_heur_mod:020}
    \For{$p \in \mathcal{L}$\label{alg:ins_heur_mod:030}}{ 
        $d \gets p+n$\; \label{alg:ins_heur_mod:040}
        $\mathcal{C}, \ I\mathcal{C} \gets$ GET\_BEST\_FEASIBLE\_AND\_INFEASIBLE\_INSERTION\_CANDIDATE($p$, $d$, $S$, $I$)\; \label{alg:ins_heur_mod:050}
        \eIf{$\mathcal{C}.feasible$}{ \label{alg:ins_heur_mod:060}
            $S \gets$ UPDATE\_SOLUTION($S$, $p$, $d$, $\mathcal{C}$, $I$)\; \label{alg:ins_heur_mod:070}
        }{\label{alg:ins_heur_mod:080}
            $S \gets$ UPDATE\_SOLUTION\_WITH\_INSTANCE\_MODIFICATIONS($S$, $p$, $d$, $I\mathcal{C}$, $I$)\; \label{alg:ins_heur_mod:090}
            $S.feasible \gets \textbf{false}$\; \label{alg:ins_heur_mod:100}
        }
    }
    \If{!$S.feasible$}{\label{alg:ins_heur_mod:105}
        Save the modified instance $I$ as a new instance\; \label{alg:ins_heur_mod:107}
    }
    \Return $S$\; \label{alg:ins_heur_mod:110}
\end{algorithm}

Algorithm \ref{alg:best_feas_infeas_ins_cand} describes the pseudocode for the GET\_BEST\_FEASIBLE\_AND\_INFEASIBLE\_{\linebreak}INSERTION\_CANDIDATE called in line~\ref{alg:ins_heur_mod:050} of Algorithm \ref{alg:ins_heur_mod}, which searches for the best feasible and infeasible insertion candidates for the nodes $p$ and $d$ in $S$.
In line~\ref{alg:best_feas_infeas_ins_cand:010}, we initialize the best candidate $B\mathcal{C}$ with an infeasible dummy insertion candidate $\mathcal{C}$ with $\mathcal{C}.cost = \infty$ to serve as a sentinel.
Then, in lines \ref{alg:best_feas_infeas_ins_cand:030}-\ref{alg:best_feas_infeas_ins_cand:080}, we analyze the insertion of $p$ and $d$ in positions $p_{pos}$ and $d_{pos}$ of the route $S.vehicles[k]$ for vehicle $k \in K$, respectively (line~\ref{alg:best_feas_infeas_ins_cand:040}).
If the insertion candidate $\mathcal{C}$ is feasible and its cost $\mathcal{C}.cost$ is smaller than $B\mathcal{C}.cost$ (line~\ref{alg:best_feas_infeas_ins_cand:050}), we update $BI\mathcal{C}$ setting candidate $\mathcal{C}$ (line~\ref{alg:best_feas_infeas_ins_cand:060}).
Otherwise, if $\mathcal{C}$ is infeasible and its cost is is smaller than $BI\mathcal{C}.cost$ (line~\ref{alg:best_feas_infeas_ins_cand:070}), we update $BI\mathcal{C}$ setting candidate $\mathcal{C}$ (line~\ref{alg:best_feas_infeas_ins_cand:080}).
Note that, the evaluation order determines $B\mathcal{C}$ to be the first $\mathcal{C}$ evaluated with $\mathcal{C}.cost = B\mathcal{C}.cost$.
The evaluation order follows the given vehicle order defined in the data set, in addition to the values of $p_{pos}$ and $d_{pos}$ indicated in line~\ref{alg:best_feas_infeas_ins_cand:030}.
The same applies to the best infeasible insertion candidate $BI\mathcal{C}$.
After analyzing all candidates, we return both $B\mathcal{C}$ and $BI\mathcal{C}$ (line~\ref{alg:best_feas_infeas_ins_cand:090}).

\begin{algorithm}[H]
    \scriptsize
    \caption{\footnotesize GET\_BEST\_FEASIBLE\_AND\_INFEASIBLE\_INSERTION\_CANDIDATE($p$, $d$, $S$, $I$)}
    \label{alg:best_feas_infeas_ins_cand}
    $B\mathcal{C}, \ BI\mathcal{C} \gets$ Dummy insertion candidates $\mathcal{C}$ with $\mathcal{C}.cost = \infty$ and $\mathcal{C}.feasible = \textbf{false}$\; \label{alg:best_feas_infeas_ins_cand:010}
    \For{$k \in K$ \textbf{and} $p_{pos} \in 2:S.vehicles[k].length$ \textbf{and} $d_{pos} \in p_{pos}:S.vehicles[k].length$}{ \label{alg:best_feas_infeas_ins_cand:030}
        $\mathcal{C} \gets$ ANALYZE\_INSERTION($k,p_{pos}, d_{pos}, p, d$, $I$)\; \label{alg:best_feas_infeas_ins_cand:040}
        \uIf{$\mathcal{C}.feasible$ \textbf{and} $\mathcal{C}.cost < B\mathcal{C}.cost$}{ \label{alg:best_feas_infeas_ins_cand:050}
            $B\mathcal{C} \gets \mathcal{C}$\; \label{alg:best_feas_infeas_ins_cand:060}
        }
        \ElseIf{\textbf{not} $\mathcal{C}.feasible$ \textbf{and} $\mathcal{C}.cost < BI\mathcal{C}.cost$}{        \label{alg:best_feas_infeas_ins_cand:070}
                $BI\mathcal{C} \gets \mathcal{C}$\; \label{alg:best_feas_infeas_ins_cand:080}
        }
    }
    \Return $B\mathcal{C}, \ BI\mathcal{C}$\; \label{alg:best_feas_infeas_ins_cand:090}
\end{algorithm}

\subsection{Insertion analysis}
\label{app:ins_feas_cost}

The ANALYZE\_INSERTION invoked in line~\ref{alg:best_feas_infeas_ins_cand:040} of Algorithm \ref{alg:best_feas_infeas_ins_cand} is similar to the ANALYZE\_{\linebreak}INSERTION\_FEASIBILITY algorithm described in \ref{app:feasibility_analysis}.
It will also examine the nodes sequentially from the pickup position $p_{pos}$ for time and capacity feasibility, but it will not stop until all nodes are examined.
If no constraints were violated, it will return that the request insertion is feasible and its cost, which is vehicle arrival time at the depot.
Otherwise, it will return that the insertion is infeasible and its cost of modifications.
This cost is defined to be $tw\_shifts+exc\_cap$, where $tw\_shifts$ is the sum of time window shifts, and $exc\_cap$ is the amount of vehicle capacity exceeded, weighted by the vehicle capacity.
We weight the amount exceeded by the vehicle capacity to discourage the usage of a larger vehicle.
If there is no vehicle type larger than the current route's total load, we reject the insertion of the request in this position.
We can reject this insertion because, as mentioned in Section \ref{subsubsec:inst-chars}, we also guarantee, by construction, that the algorithm will find another vehicle to load this request.
It is worth remembering that, from $p_{pos}-1$ onward, all machine travels transporting the vehicle are disregarded in this analysis.
Thus, whenever new machine travel is necessary, it will also be recalculated.

\subsection{Vehicle travel time}
\label{app:vehicle_travel_time_mod}

The algorithm to get the best vehicle travel time in arc $(i,j) \in A^m$ is similar to the algorithm described in \ref{app:vehicle_travel_time}.
The only difference is that we do not stop the analysis at a machine sequence of travels $S.machines[h], \ h \in H_{ij}$ if the new machine travel insertion leads to violate $l_j$.
All machine travels that do not violate other machine travels are considered feasible and $best\_mtrv$ can be added to $mtrvs$, which is the list of new machine travels calculated in the ANALYZE\_INSERTION algorithm.

\subsection{Instance modification}
\label{app:instance_modification}

The UPDATE\_SOLUTION\_WITH\_INSTANCE\_MODIFICATIONS called in line~\ref{alg:ins_heur_mod:090} of the Algorithm \ref{alg:ins_heur_mod} also modifies the instance to \okC{serve} some request.
The modifications are performed during the solution update (see \ref{app:update_solution}).
When the time window of a node $i \in V$ is violated, we adjust it by setting the latest time ($l_i$) to the smallest integer greater than or equal to the arrival time of the vehicle in node $i$, while adjusting the earliest time ($e_i$) to maintain the original window width.
When the capacity limit is exceeded, we change the vehicle type to the one with the smallest capacity that is still larger than the current route's total load.
An overview of the instance changes is presented in Appendix~\ref{app:overview_inst_changes}.

% \newpage
\section{{Overview of the instance changes}}
\label{app:overview_inst_changes}

\okC{
In this section, we provide an overview of how many instances were modified and how much each instance was changed by the procedure to ensure instance feasibility (see \ref{app:alg_to_ensure_feas}).}
Table~\ref{tab:inst_change_overview} presents an overview of the instance changes for each family and instance type.
\okC{It is also subdivided by instance size: the first set with 6 to 12 requests, and the second set with 40 and 60 requests.
Given that half of the instances are copies without the last machine (see Section \ref{subsubsec:inst-chars}), the total number of instances for each type and family analyzed in the first and second sets is 40 and 20, respectively.}
The row ``Instances modified (TW/Cap)" shows the percentage of instances with a modification, whether in a time window or in a vehicle capacity.
The row ``Instances with TW shift" shows the percentage of instances with shifted time windows for at least one node.
The row ``Instances with Cap increase" shows the percentage of instances with increased capacity for at least one vehicle.
\okC{The table indicates a prevalence of instance changes on short scheduling horizon instances (Type 1), especially on instances with 40 and 60 requests.
For the first instance set containing 6 to 12 requests, fewer than 8\% of the instances have increased their vehicle capacities.
Conversely, 70\% of Type 1 instances with 40 and 60 requests had their vehicle capacities increased, although only 10\% of the multi-floor Type 2 instances and no multi-island Type 2 instances were modified.
It is also clear that most of the instances modified are from the multi-floor family.}

\begin{table}[!ht]
    \centering
    \caption{Overview of the instance changes}
    \label{tab:inst_change_overview}
    \scriptsize
    \okC{
    \begin{tabular}{lr c rrrr c rrrr}\toprule
        \multirow{3}{*}{} & &\multicolumn{4}{c}{6 to 12 requests} & &\multicolumn{4}{c}{40 and 60 requests} \\
        \cmidrule{3-6} \cmidrule{8-11}
        & &\multicolumn{2}{c}{Multi-island} &\multicolumn{2}{c}{Multi-floor} & &\multicolumn{2}{c}{Multi-island} &\multicolumn{2}{c}{Multi-floor} \\
        \cmidrule{3-6} \cmidrule{8-11}
        & &Type 1 &Type 2 &Type 1 &Type 2 & &Type 1 &Type 2 &Type 1 &Type 2 \\
        % \cmidrule{1-1} \cmidrule{3-6} \cmidrule{8-11}
        \midrule
        Instances modified (TW/Cap) & &40.00\% & 0.00\% & 92.50\% & 22.50\% & & 100.00\% & 20.00\% & 100.00\% & 50.00\% \\
        Instances with TW shift & &40.00\% & 0.00\% & 92.50\% & 22.50\% & &  100.00\% & 20.00\% & 100.00\% & 50.00\% \\
        Instances with Cap increase & &2.50\% & 0.00\% & 2.50\% & 7.50\% & & 70.00\% & 0.00\% & 70.00\% & 10.00\% \\
        \bottomrule
\end{tabular}
}
\end{table}

\begin{figure}[!ht]
    \centering
    \subfigure[Multi-island instances]{
        \centering
        \includegraphics[scale=0.42]{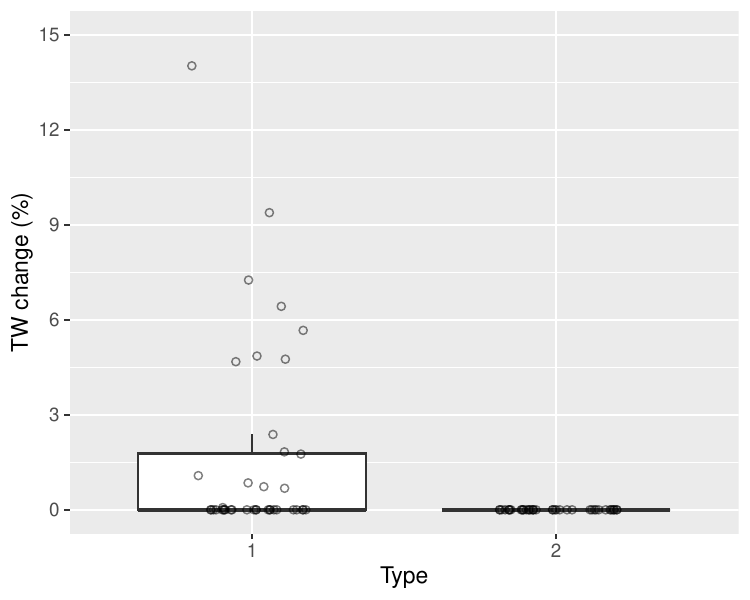}
    }
    \subfigure[Multi-floor instances]{
        \centering
        \includegraphics[scale=0.42]{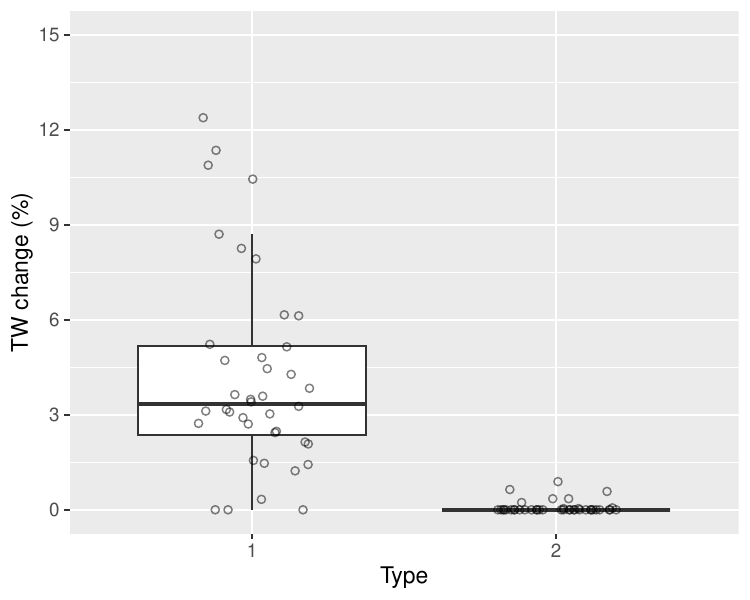}
    }
    \caption{Box plot for Time Window change by type for smaller instances.}
    \label{fig:tw-change-by-type}
\end{figure}

\begin{figure}[!ht]
    \centering
    \subfigure[Multi-island instances]{
        \centering
        \includegraphics[scale=0.42]{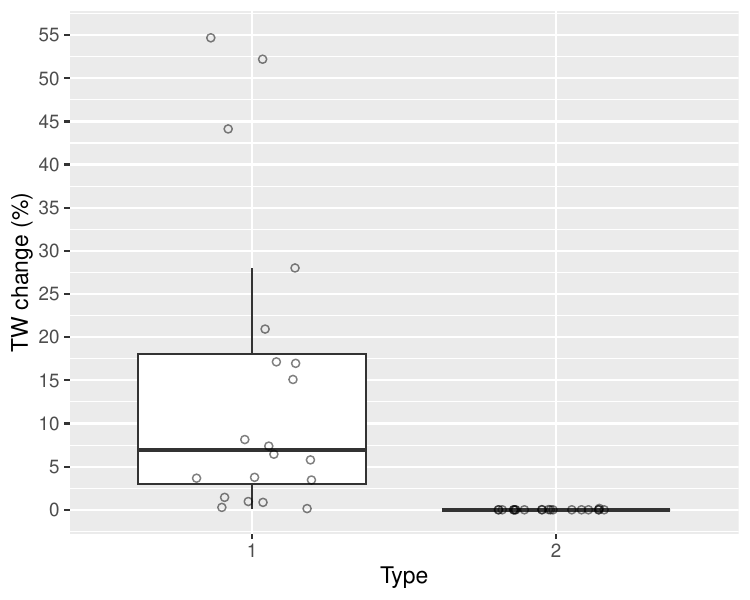}
    }
    \subfigure[Multi-floor instances]{
        \centering
        \includegraphics[scale=0.42]{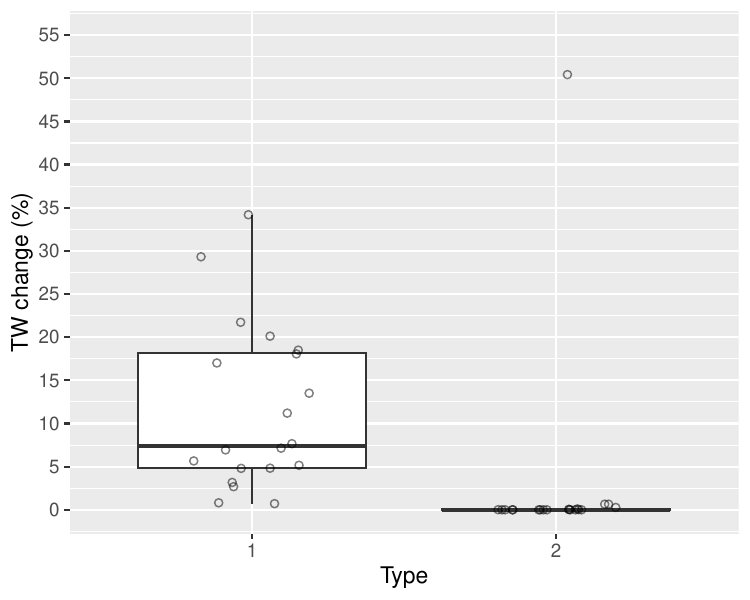}
    }
    \caption{Box plot for Time Window change by type for bigger instances.}
    \label{fig:tw-change-by-type-big}
\end{figure}

\okC{The box plots in Figures~\ref{fig:tw-change-by-type} and \ref{fig:tw-change-by-type-big} depict the mean percentage distribution of time window changes.
In Figure \ref{fig:tw-change-by-type}, each box contains 40 observations representing the 80 instances with 6-12 requests, three and four machines, and the corresponding type.
In Figure \ref{fig:tw-change-by-type-big}, each box contains 20 observations representing the 40 instances with 40 and 60 requests, five and six machines, and the corresponding type.}
A node time window change is calculated as the amount of time shifted divided by its original time window's latest time.
The time window modification for an instance is the mean of those shifts.

Machine travels increase the distance between nodes in different regions, making it difficult to serve requests with close time windows.
\okC{That is clear in the box plots of Figures~\ref{fig:tw-change-by-type} and \ref{fig:tw-change-by-type-big}, given that changes are prevalent in Type 1 instances.}
\okC{Regarding the smaller instances, it's} also remarkable that more instances in the multi-floor \okC{family} were changed\okC{, especially those of Type~1}. \okC{The reason is that the request nodes are distributed uniformly in the multi-floor instances.
Consequently, especially for smaller instances, the requests with originally closer time windows have higher chances of being separated into different regions, increasing the traversal time between them.
However, for larger instances with the same number of regions, the chances of a request being in the same region also increase.
Thereupon, instance changes became more balanced between multi-island and multi-floor families.}
\okC{Furthermore, there are significantly more changes on the larger Type 1 instances, but more than 75\% of the observations have a mean time window change of less than 20\%.}

\begin{figure}[!htp]
    \centering
    \subfigure[Multi-island instances]{
        \centering
        \includegraphics[scale=0.42]{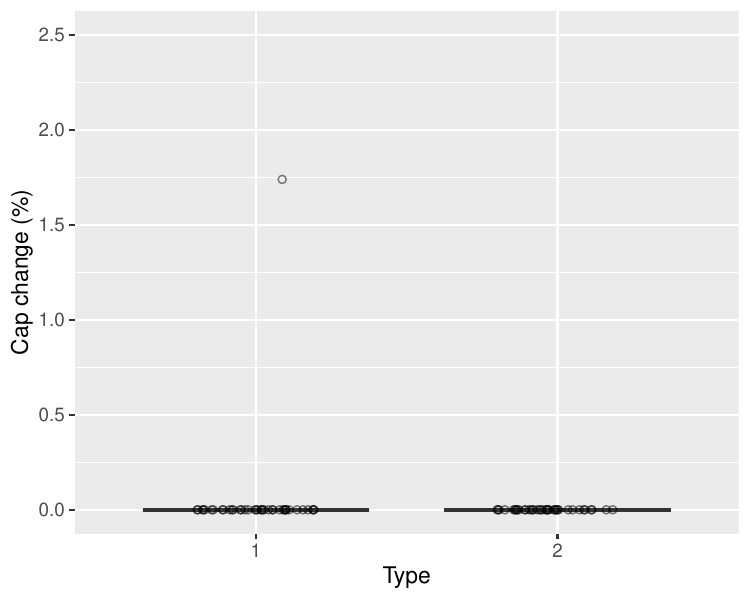}
    }
    \subfigure[Multi-floor instances]{
        \centering
        \includegraphics[scale=0.42]{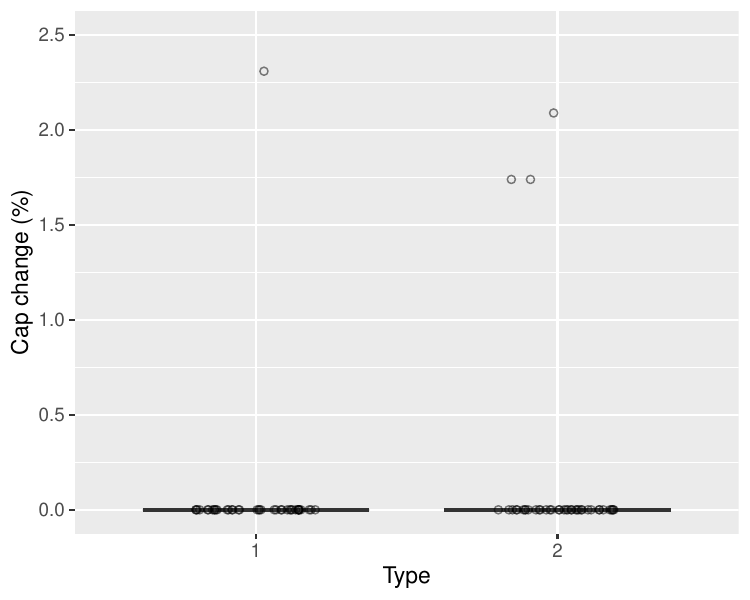}
    }
    \caption{Box plot for Capacity change by type.}
    \label{fig:cap-change-by-type}
\end{figure}

\begin{figure}[!ht]
    \centering
    \subfigure[Multi-island instances]{
        \centering
        \includegraphics[scale=0.42]{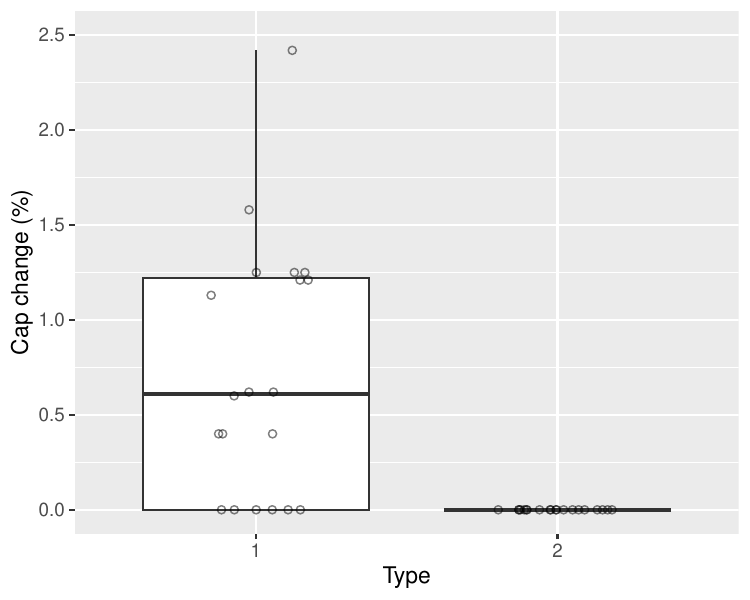}
    }
    \subfigure[Multi-floor instances]{
        \centering
        \includegraphics[scale=0.42]{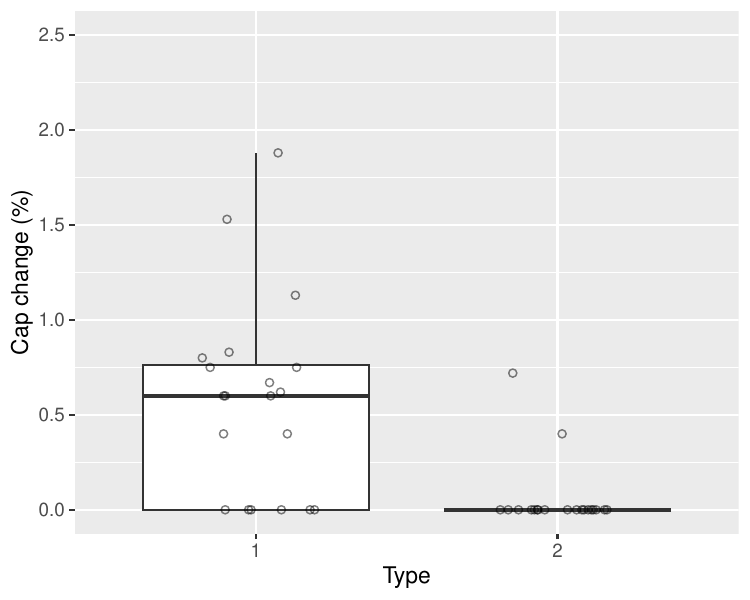}
    }
    \caption{Box plot for Capacity change by type for bigger instances.}
    \label{fig:cap-change-by-type-big}
\end{figure}

A vehicle capacity change is the amount of capacity increased in comparison with its original capacity.
The capacity modification for an instance is the mean of these increments.
\okC{
The box plots in Figure~\ref{fig:cap-change-by-type} and \ref{fig:cap-change-by-type-big} present the mean percentage changes in the vehicle capacities. 
The number of observations is analogous to the number of observations in Figures \ref{fig:tw-change-by-type} and \ref{fig:tw-change-by-type-big}.
}
By construction, it is not necessary to increase a vehicle's capacity, as mentioned in section \ref{subsubsec:inst-chars}.
However, the algorithm described in  \ref{app:alg_to_ensure_feas} can find a lower-cost alternative to \okC{serving} a request by combining both vehicle capacity increases and time window shifts.
\okC{
For this reason, some vehicle capacities increased on the smaller instances, and on several Type 1 larger instances.
}

\FloatBarrier

% \newpage
\section{\okC{Solver comparison: Gurobi vs. Hexaly}}
\label{app:grb_vs_hx}

\begin{figure}[!ht]
    \centering
    \subfigure[Multi-island instances]{
        \centering
        \includegraphics[scale=0.42]{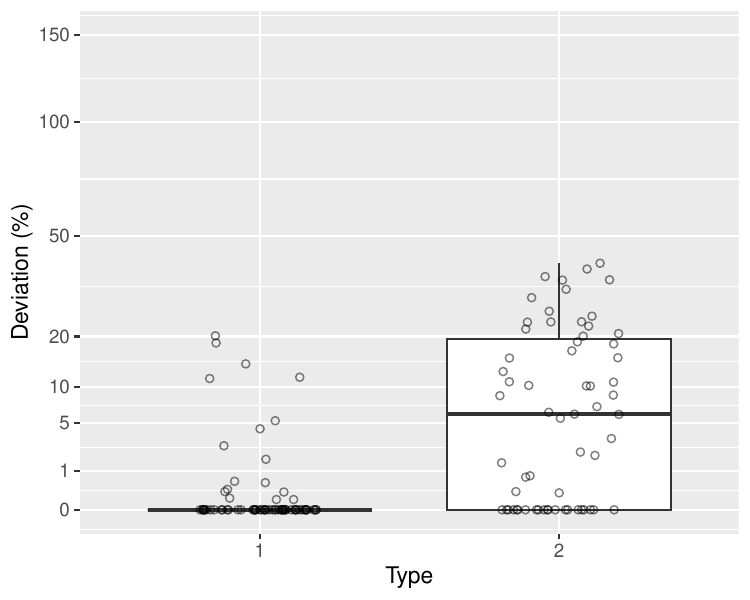}
    }
    \subfigure[Multi-floor instances]{
        \centering
        \includegraphics[scale=0.42]{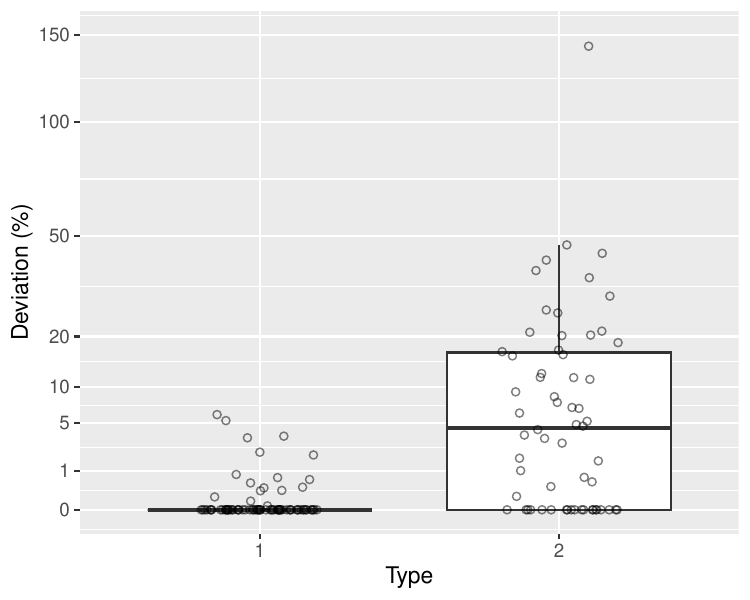}
    }
    \caption{Box plots for Gurobi-Hexaly solution value deviation by type. Positive deviations indicate Gurobi's outperformance.}
    \label{fig:grb-vs-hx-by-type}
\end{figure}

\okC{In this section, we provide a comparison between the solution values obtained by the solvers Gurobi and Hexaly on the baseline MIP for the benchmark instances.
The box plots in Figure~\ref{fig:grb-vs-hx-by-type} depict the Gurobi-Hexaly solution value deviation for each instance family, separated by instance type.
Gurobi's results are the same as presented in Section \ref{subsec:mipresults_no_vis}, while Hexaly's results were obtained by running the solver for each instance with the same time limit of 3600 seconds.
For Hexaly, the baseline MIP was implemented in Python 3.10.12, using Hexaly v14.0.
Each box contains at most 80 observations, representing the 80 instances with the corresponding type.
Missing observations occur because at least one of the solvers was unable to find a feasible solution for an instance.
Table \ref{tab:missing_obs_box_plots} summarizes the number of missing observations in each box plot, separated by instance family and type.
Column ``GRB or HX" shows the total number of missing solutions.
Column ``GRB and HX" indicates cases where both solvers failed.
Columns ``GRB" and ``HX" count cases where only the corresponding solver (Gurobi or Hexaly) is missing a solution.
It is worth mentioning that for only three out of the 35 (8.6\%) Type 2 instances with missing observations, only Hexaly found a feasible solution: two multi-island and one multi-floor.
The deviation is calculated as $100 \cdot\frac{sol^{\text{HX}} -sol^{\text{GRB}}}{sol^{\text{GRB}}}$, \okC{in which $sol^{\text{HX}}$ and $sol^{\text{GRB}}$ are the solution values of Hexaly's and Gurobi's executions on an instance.}
Therefore, negative values indicate that Hexaly's solution value surpasses Gurobi's solution value.
Otherwise, no improving solution was found.
}

\begin{table}[!ht]
    \centering
    \scriptsize
    \okC{
    \begin{tabular}{l r rrrr c rrrr}\toprule
        & &\multicolumn{4}{c}{Multi-island} & &\multicolumn{4}{c}{Multi-floor} \\
        \cmidrule{3-6}
        \cmidrule{8-11}
         & &GRB or HX &GRB and HX &GRB &HX & &GRB or HX &GRB and HX &GRB &HX \\
        \cmidrule{1-1}
        \cmidrule{3-6}
        \cmidrule{8-11}
        Type 1 & &5/80  &0/80 &0/80 &5/80 & &0/80 &0/80 &0 &0 \\
        Type 2 & &17/80 &7/80 &2/80 &8/80 & &18/80 &9/80 &1/80 &8/80 \\
        \bottomrule
    \end{tabular}
    }
    \caption{\okC{Missing observations in each box plot}}
    \label{tab:missing_obs_box_plots}
\end{table}

\okC{
As seen in Figure~\ref{fig:grb-vs-hx-by-type}, all deviations are at least greater than or equal to zero.
Thus, Gurobi obtained better results when compared to Hexaly.
This difference is even greater in Type 2 instances for both families, where more than 50\% of the observations have a deviation greater than 7.5\%.
}

\FloatBarrier

% \newpage
\section{Analysis of time and optimality gap for the MIP formulation}
\label{app:mip_results_details}

\begin{figure}[H]
    \centering
    \subfigure[Multi-island instances\label{fig:time_n_reqs_1}]{
        \centering
        \includegraphics[scale=0.42]{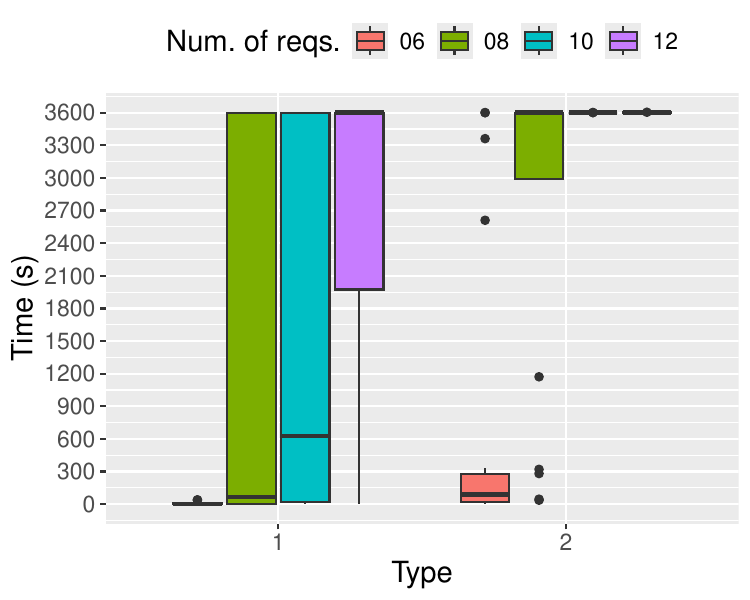}
    }
    \subfigure[Multi-floor instances\label{fig:time_n_reqs_2}]{
        \centering
        \includegraphics[scale=0.42]{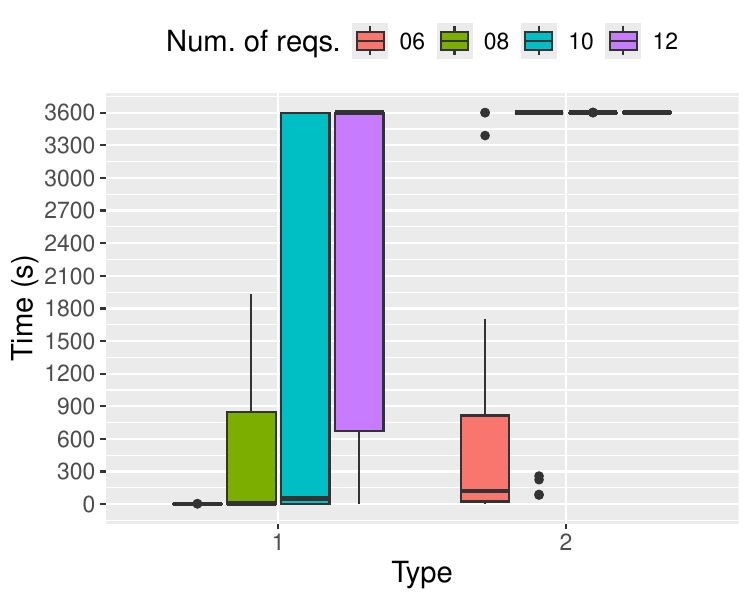}
    }
    \caption{Time x Num. of Requests for the MIP results}
    \label{fig:time_vs_num_reqs}
\end{figure}

\begin{figure}[H]
    \centering
    \subfigure[Multi-island instances]{
        \centering
        \includegraphics[scale=0.42]{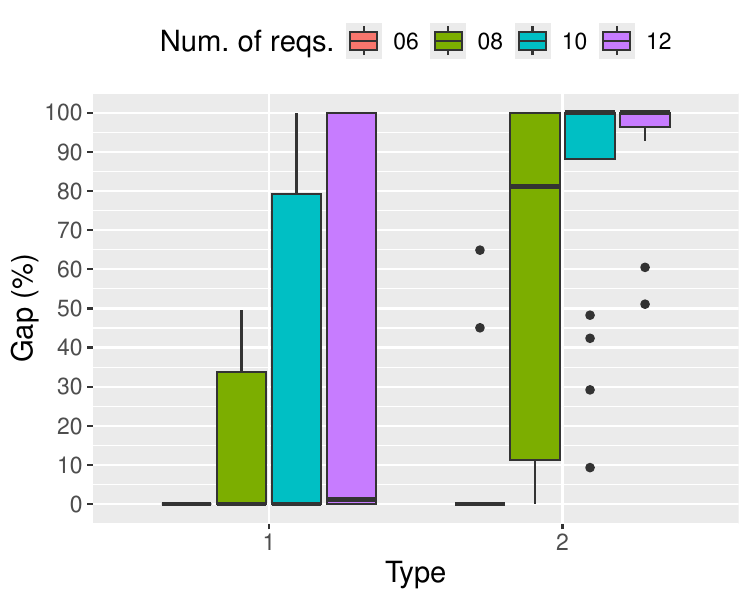}
    }
    \subfigure[Multi-floor instances]{
        \centering
        \includegraphics[scale=0.42]{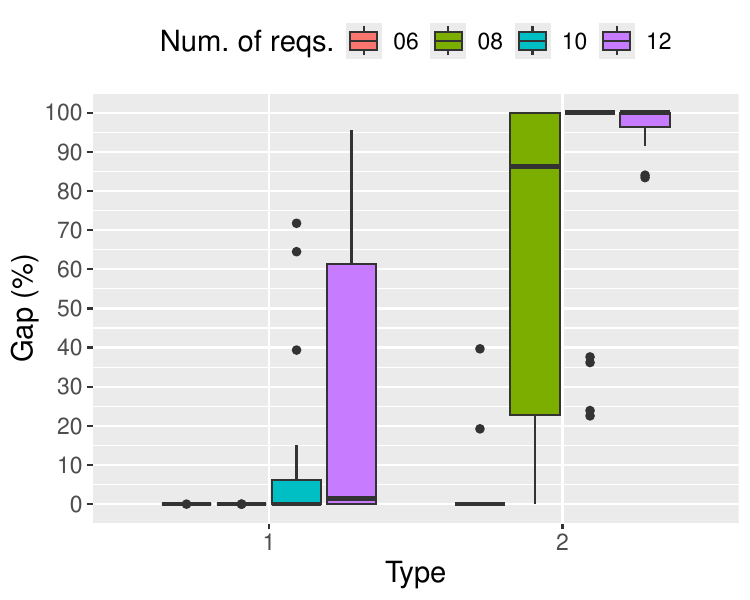}
    }
    \caption{Gap x Num. of Requests for the MIP results}
    \label{fig:gap_vs_num_reqs}
\end{figure}

\okC{
The box plots in Figures \ref{fig:time_vs_num_reqs} and \ref{fig:gap_vs_num_reqs} depict the distributions of, respectively, running time and gap of the MIP with valid inequalities for the instances with the same number of requests in each family.
For each box plot representing the Type~1 instances and the Type~2 instances with up to eight requests, there are 20 observations, given that the solver found a feasible solution for all of them.
Nevertheless, the number of observations for Type~2 instances with at least ten requests is presented in Table \ref{tab:num_obs_box_plots_time_gap_mip_vi}.
}

\begin{table}[h]
    \centering
    \scriptsize
    \okC{
    \begin{tabular}{lrrrr}
        \toprule
           Size     & & Multi-island & & Multi-floor \\
        \cmidrule{1-1} \cmidrule{3-3} \cmidrule{5-5}
         10 reqs.   & & 17 & & 18 \\
         12 reqs.   & & 14 & & 15 \\
         \bottomrule
    \end{tabular}
    }
    \caption{\okC{Number of observations in each box plot of Figures \ref{fig:time_vs_num_reqs} and \ref{fig:gap_vs_num_reqs} for Type~2 instances with at least ten requests.}}
    \label{tab:num_obs_box_plots_time_gap_mip_vi}
\end{table}

\okC{Regarding running time distributions in Figure \ref{fig:time_vs_num_reqs}, it can be seen that, for both multi-island and multi-floor Type 1 instances, the median increases exponentially as the number of requests increases.
However, the time observations are too sparsed in each box plot for the Type~1 multi-island instances with at least eight requests, and for the Type~1 multi-floor instances with at least ten requests.
For Type~2 instances, the times increase even earlier, with some six-request instances and more than 50\% of the eight-request instances already reaching the time limit.
Conversely, the box plots of Figure \ref{fig:gap_vs_num_reqs} indicate that, for both multi-island and multi-floor Type 1 instances, the median gaps approach zero percent.
Similarly as in the time distributions, the gap observations are also too sparsed in each box plot for the Type~1 multi-island instances with at least ten requests, and for the Type~1 multi-floor instances with twelve requests.
For Type~2 instances, the solver achieved gaps above 80\% for almost all instances with at least ten requests, and for more than 50\% of the instances with eight request in both families.}

%\newpage
\section{MSLP preliminary tests} 
\label{app:preliminary_tests}

\okC{
To evaluate the performance of the MSLP heuristic, we first need to configure the stop criterion of the multistart framework (see Section \ref{sec:heur_gen_framework}) and the threshold parameter $\alpha$ for the semi-greedy heuristic (see Section \ref{sec:semi-greedy-ins-heur}).
To set the threshold parameter $\alpha$, preliminary tests were conducted on additional instances generated from the sixth to the tenth instance of each type for the PDPTW:
lr106 - lr110 (Type 1) and lr206 - lr210 (Type 2), or LR1\_2\_6 - LR1\_2\_10 (Type 1) and LR2\_2\_6 - LR2\_2\_10 (Type 2).
These instances were generated in the same way as described in Section \ref{subsec:inst_gen}, varying the number of requests ($n \in \{8,12,40,60\}$), vehicles ($|K| = n$), regions ($z \in \{2,4\}$) and machines ($|H| \in \{3, 4, 5, 6\}$).
From each one of the 32 possible group variations, one instance of each type for each possible variation was selected.
In total, we selected 64 instances, which corresponds to 13.3\% of the total number of instances (480) created for the benchmark.
We set the stop criterion to be 60000 iterations or 3600 seconds, the same as in the benchmark tests.
}

\okC{
The preliminary tests were conducted for $\alpha \in \{0.05, 0.10, 0.15, 0.20\}$.
Given the stochastic nature of the algorithm, 20 independent runs were performed for each instance with each configuration.
We summarize the results in Table \ref{tab:mslp_preliminar_grouped_type_alpha} for each $\alpha$ value by calculating the mean \okC{($\pm$ standard deviation)} of the minimum, mean, and maximum RPDs to the best solution value found for any $\alpha$ value or independent execution, separated by instance type.
Notice that the RPD values for the parameter $\alpha = 0.05$ beat the values for other parameter $\alpha$ values, except the mean of the minimal RPDs for Type 1 instances.
However, its difference from the best value is not significantly better.
Therefore, we selected $\alpha=0.05$ as the parameter value to perform the experiments using the test data.
}

\begin{table}[ht] 
\centering 
\scriptsize 
\caption{\okC{Solution value RPD summary for each parameter configuration on the MSLP preliminary tests} \label{tab:mslp_preliminar_grouped_type_alpha}} 
\okC{
\begin{tabular}{rrlllrlll} 
\toprule & &\multicolumn{3}{c}{Type 1} & & \multicolumn{3}{c}{Type 2} \\
\cmidrule{3-5} \cmidrule{7-9}
$\alpha$ & & Min RPD & Mean RPD & Max RPD & & Min RPD & Mean RPD & Max RPD \\
\midrule 
0.05 & & 0.41 $\pm$ 0.82 & 1.54 $\pm$ 1.68 & 2.57 $\pm$ 2.75 & & 0.31 $\pm$ 0.85 & 2.58 $\pm$ 2.24 & 4.15 $\pm$ 3.16 \\
0.10 & & 0.39 $\pm$ 0.71 & 1.82 $\pm$ 1.95 & 2.84 $\pm$ 3.00 & & 1.54 $\pm$ 1.83 & 3.79 $\pm$ 3.08 & 5.32 $\pm$ 3.92 \\
0.15 & & 0.87 $\pm$ 1.38 & 2.14 $\pm$ 2.33 & 3.12 $\pm$ 3.14 & & 2.52 $\pm$ 2.19 & 4.56 $\pm$ 3.60 & 5.99 $\pm$ 4.44 \\
0.20 & & 0.90 $\pm$ 1.24 & 2.44 $\pm$ 2.60 & 3.35 $\pm$ 3.47 & & 2.59 $\pm$ 2.44 & 4.66 $\pm$ 3.60 & 6.30 $\pm$ 4.61 \\ 
\bottomrule
\end{tabular} 
}
\end{table}

%\newpage
\section{\okC{Solution feasibility throughout the MSLP iterations}}
\label{app:mslp_solution_feasibility}

\begin{figure}[!ht]
        \centering
        \begin{minipage}{.65\textwidth}
            \subfigure[Multi-island instances vs. \# reqs.]{
            \label{fig:feasible_solutions_percentage_n_reqs_1}
                \centering
                \includegraphics[scale=0.42]{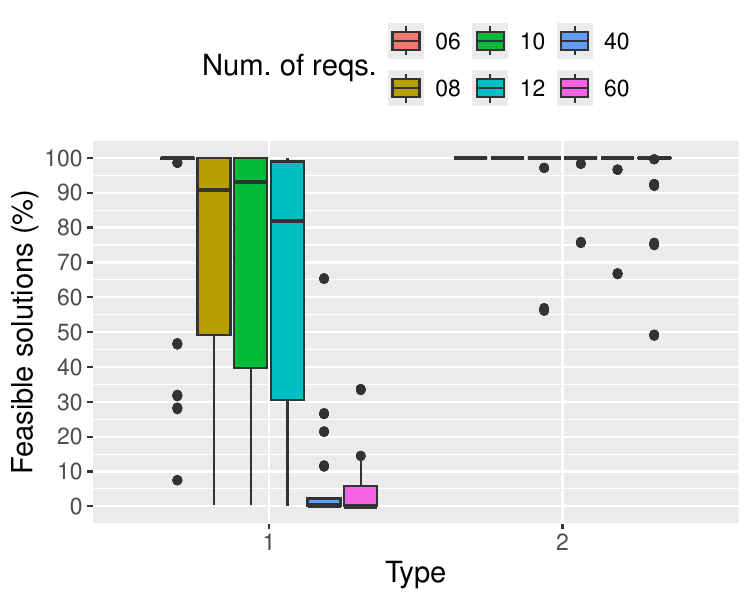}
            }
            \subfigure[Multi-floor instances vs. \# reqs.]{
            \label{fig:feasible_solutions_percentage_n_reqs_2}
                \centering
                \includegraphics[scale=0.42]{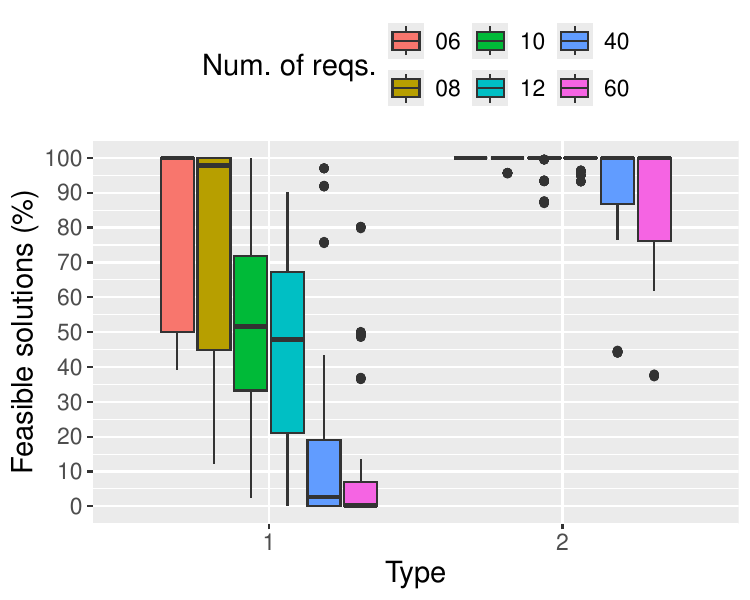}
            }
            \\
            \subfigure[Multi-island instances vs. \# regs.]{
            \label{fig:feasible_solutions_percentage_n_regions_1}
                \centering
                \includegraphics[scale=0.42]{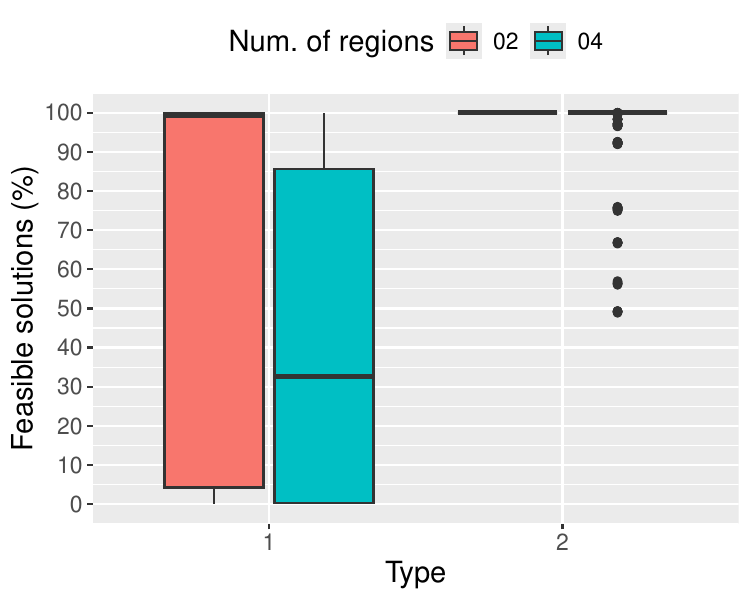}
            }
            \subfigure[Multi-floor instances vs. \# regs.]{
            \label{fig:feasible_solutions_percentage_n_regions_2}
                \centering
                \includegraphics[scale=0.42]{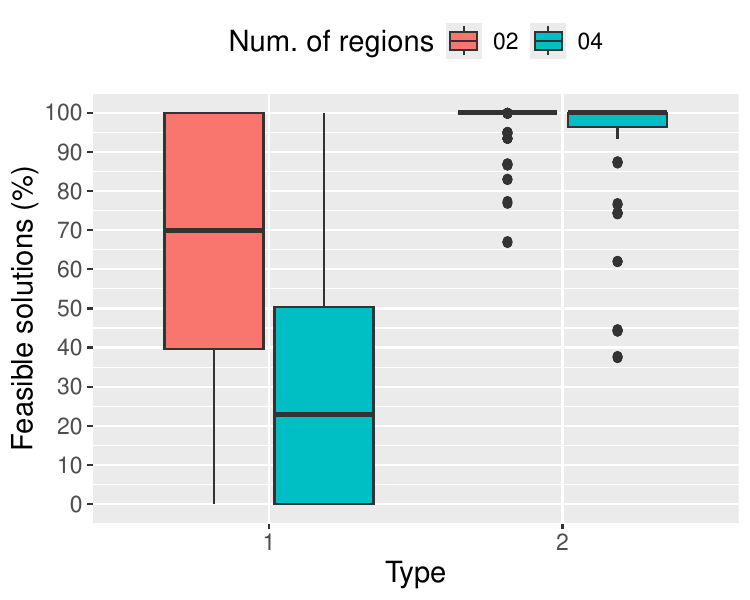}
            }
            \\
            \subfigure[Multi-island instances vs. \# machs.]{
            \label{fig:feasible_solutions_percentage_n_machs_1}
                \centering
                \includegraphics[scale=0.42]{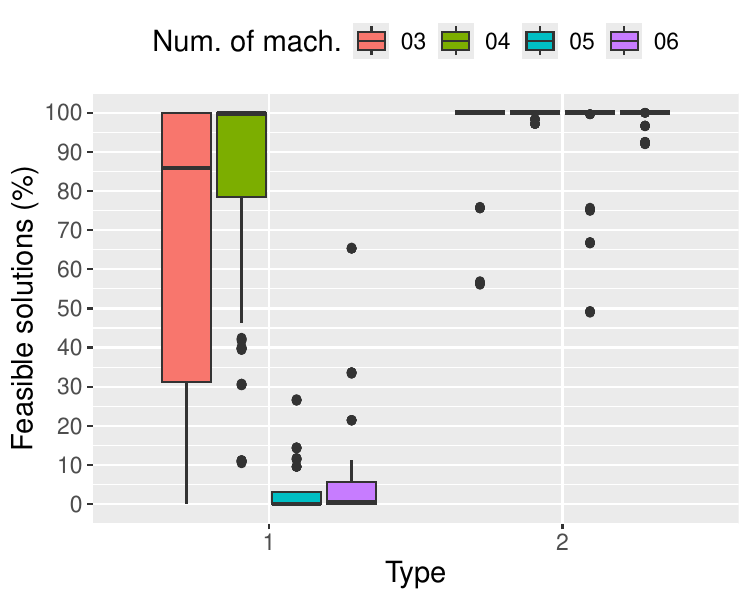}
            }
            \subfigure[Multi-floor instances vs. \# machs.]{
            \label{fig:feasible_solutions_percentage_n_machs_2}
                \centering
                \includegraphics[scale=0.42]{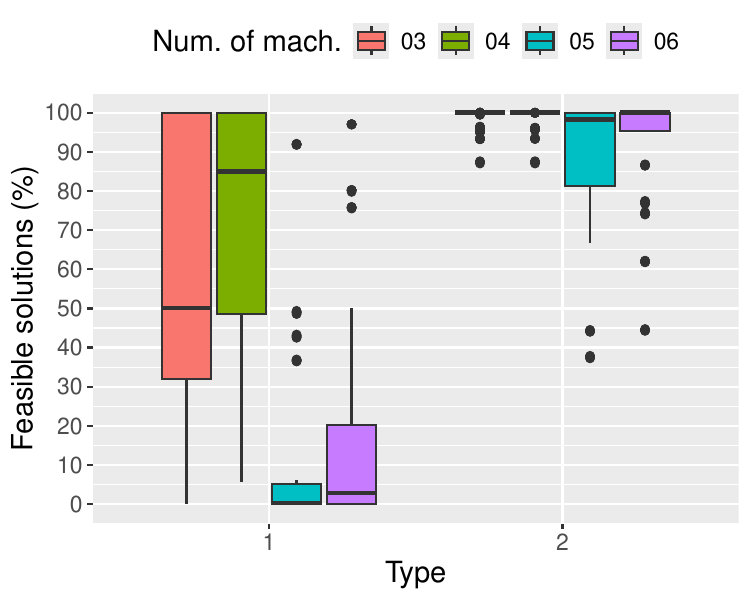}
            }
        \end{minipage}
    \caption{Feasible solutions (\%) for the MSLP heuristic, varying: \# requests, \# regions and \# machines.}
    \label{fig:feasible_solutions_percentage}
\end{figure}

\okC{
The box plots of Figure \ref{fig:feasible_solutions_percentage} depict the percentage distribution of the solution feasibility over the iterations of the ten independent executions of the MSLP heuristic when varying the number of: requests, regions, and machines.
Each box plot for request variation contains 200 observations. 
For region variation, there are 600 observations in each box plot.
Concerning the box plots for machine variation, remember that there are 320 instances with three or four machines, while there are 160 instances with five or six machines.
Therefore, the box plots related to three or four machine instances have each 400 observations, while the box plots related to five or six machine instances have each 200 observations.
Observations for instances with no feasible solution are included.
}

\okC{
Observe in Figures \ref{fig:feasible_solutions_percentage_n_reqs_1} and \ref{fig:feasible_solutions_percentage_n_reqs_2} that although the percentage of feasible solutions decreases as the instance size increases, the MSLP heuristic can consistently find feasible solutions for most of the Type~2 instances and Type~1 instances with up to twelve requests.
In Figures \ref{fig:feasible_solutions_percentage_n_regions_1} and \ref{fig:feasible_solutions_percentage_n_regions_2}, it is also clear that increasing the number of regions also decreases the percentage of feasible solutions.
Nonetheless, the percentage of feasible solutions is still above 20\% for at least 50\% of the Type~1 instances in both families.
Finally, the box plots in Figures \ref{fig:feasible_solutions_percentage_n_machs_1} and \ref{fig:feasible_solutions_percentage_n_machs_2} notably show that the greater the number of machines, the greater the percentage of feasible solutions.
Consequently, as the number of requests increases, the MSLP heuristic might indicate whether the number of machines is an instance's bottleneck to achieving high-quality solutions.
Regarding Type~2 instances, almost all the instances obtained feasible solutions in all iterations of the MSLP heuristic, except for some outliers and multi-floor instances with 40 and 60 requests.
Therefore, it is significantly easier to construct a feasible solution using the MSLP heuristic for long-scheduling horizon instances (Type~2) compared with short-scheduling horizon instances (Type~1).
}

\FloatBarrier

\section{\okC{LP-based improvement procedure's relevance in the MSLP heuristic}}
\label{app:lps_relevance_mslp}

\begin{figure}[!ht]
    \centering
    \subfigure[Multi-island instances]{
        \centering
        \includegraphics[scale=0.42]{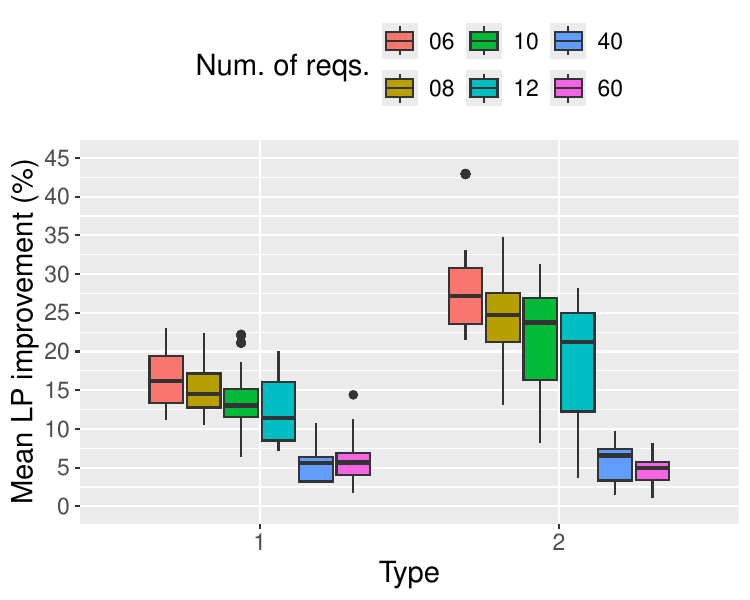}
    }
    \subfigure[Multi-floor instances]{
        \centering
        \includegraphics[scale=0.42]{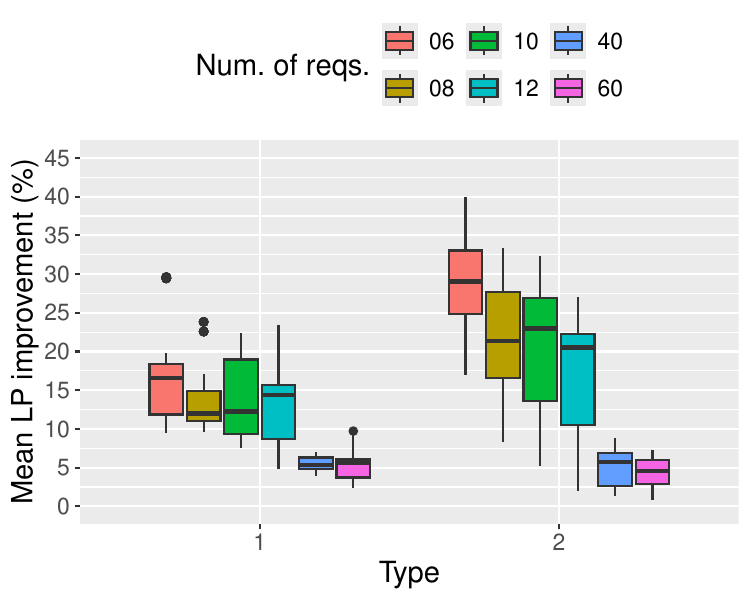}
    }
    \caption{Mean LP improvement percentage for the MSLP heuristic.}
    \label{fig:mean_lp_impr_percentage}
\end{figure}

\okC{
The box plots in Figure \ref{fig:mean_lp_impr_percentage} depict the percentage distribution of the mean LP improvements in the MSLP heuristic for both families, varying the number of requests.
Notice that this metric is only applicable for executions that found a feasible solution.
Thus, only the box plots for instances with up to 12 requests or Type~2 instances with 40 and 60 requests contain 200 observations each.
However, the box plots for 40-request instances contain 180 and 190 observations in the multi-island and multi-floor families, respectively, while the box plots for 60-request instances contain 168 and 153 observations in the multi-island and multi-floor families, respectively.
The LP improvement is calculated as: $100\times\frac{sol^{\text{greedy}} - sol^{\text{LP}}}{sol^{\text{greedy}}}$, where $sol^{\text{greedy}}$ is the (semi-)greedy solution value, and $sol^{\text{LP}}$ is the LP solution value after improving the (semi-)greedy solution.
The mean LP improvement is calculated over the number of executions of the LP procedure.
We can see from these box plots that the LP improvement procedure plays an important role in reducing the solution value of the heuristic for both families.
This improvement is even more significant for Type~2 instances, reaching a mean improvement of almost 43\% in some multi-island Type~2 instances.
}

\FloatBarrier

\section{\okC{The impact of an additional machine}}
\label{app:additional_machine_impact}

\begin{figure}[!ht]
    \centering
    \subfigure[Multi-island instances]{
        \centering
        \includegraphics[scale=0.42]{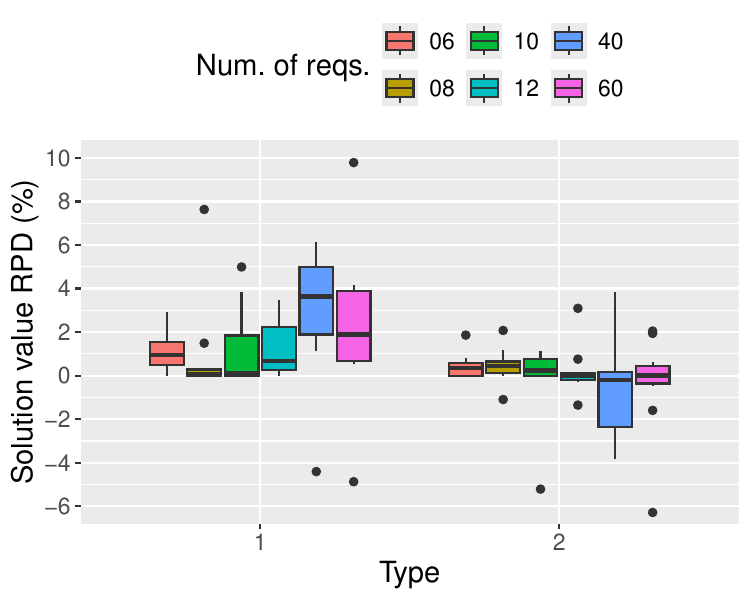}
    }
    \subfigure[Multi-floor instances]{
        \centering
        \includegraphics[scale=0.42]{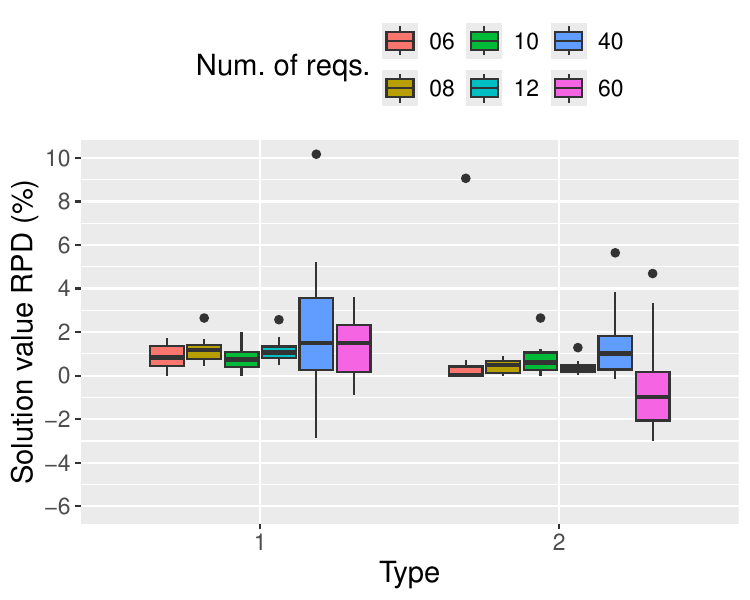}
    }
    \caption{\okC{RPD instance by instance varying \# machines - Best between MIP and MSLP.}}
    \label{fig:dev_inst_to_inst_machs_best}
\end{figure}

\okC{The box plots of Figure \ref{fig:dev_inst_to_inst_machs_best} depict the solution value RPD distribution, instance by instance, when varying the number of machines.
The deviation is calculated as: $100 \times \frac{sol^{\text{M}} - sol^{\text{M+1}}}{sol^{\text{M}}}$, where $sol^{\text{M}}$ and $sol^{\text{M+1}}$ are the best solutions obtained from the baseline MIP, the MIP with valid inequalities, and the MSLP methods for the instance with the minimum and maximum number of machines, respectively.
For the instances with up to 12 requests, the minimum is three machines, and the maximum is four machines.
Conversely, for the instances with 40 and 60 requests, the minimum is five machines, and the maximum is six machines.
Thus, a positive deviation means that providing an additional machine improved the solution value.
On the other hand, a negative deviation indicates that no method was able to improve the solution within the time limit, for MIP, or within the number of iterations of the MSLP heuristic, even though a solution for the 3-machine or 5-machine instance is also feasible for the 4-machine or 6-machine instance, respectively.
Besides that, for the box plots representing the instances with up to 12 requests and the Type~2 instances with 40 and 60 requests, each box plot contains ten observations, one for each pair of instances with the same number of requests, as the MSLP heuristic obtained feasible solutions for these instances.
However, for the box plots representing the Type~1 instances with 40 and 60 requests, we summarize the number of observations in Table \ref{tab:num_obs_box_plots_sol_rpd_varying_machs}.
We can notice that providing an additional machine to the instance can improve solution value, reaching a maximum improvement of over 10\% in a 40-request multi-floor instance.
Furthermore, we can notice that providing an additional machine was more useful for Type~1 instances on average, particularly on multi-floor instances.}

\begin{table}[h]
    \centering
    \scriptsize
    \okC{
    \begin{tabular}{lrrrr}
        \toprule
         & &  Multi-island & & Multi-floor\\
         \cmidrule{1-1} \cmidrule{3-3} \cmidrule{5-5}
        40 reqs. & & 8 & & 9\\
        60 reqs. & & 7 & & 7\\
        \bottomrule
    \end{tabular}
    }
    \caption{\okC{Number of observations in the box plots of Figure \ref{fig:dev_inst_to_inst_machs_best} for Type~1 instances with 40 and 60 requests.}}
    \label{tab:num_obs_box_plots_sol_rpd_varying_machs}
\end{table}

\FloatBarrier

%\newpage
\FloatBarrier
\section{Characteristics of the solutions}
\label{app:solution_chars}

\okC{
In this section, we present some box plots representing the characteristics of each solution of an independent execution of the MIP with valid inequalities (MIP+VI) formulation and the MSLP heuristic.
These box plots group the observations of the instances with the same number of requests, machines, or regions.
There are 20 and 200 observations per request variation box plot for the MIP+VI and the MSLP, respectively, except the cases listed in Table \ref{tab:num_obs_req_var_box_plots}.
For the MIP+VI machine variation box plot, there are 40 observations related to three- or four-machine instances, except the scenarios listed in Table \ref{tab:num_obs_mach_var_box_plots} a).
Conversely, for the MSLP machine variation box plot, there are 400 observations related to three- or four-machine instances, and 200 observations related to five- or six-machine instances, except the scenarios listed in Table \ref{tab:num_obs_mach_var_box_plots} b).
Finally, there are 40 observations per MIP+VI region variation box plot, and 600 observations per MSLP region variation box plot, except the cases listed in Table \ref{tab:num_obs_reg_var_box_plots}.
}

\begin{table}[!htp]\centering
\caption{\okC{Number of observations per request variation box plot}}
\label{tab:num_obs_req_var_box_plots}
\begin{minipage}{0.45\textwidth}
\centering
\scriptsize
\vspace{0.3em}
\okC{
\text{a) MIP+VI}\\[0.3em]
\begin{tabular}{lrr}
\toprule
 & Multi-island & Multi-floor \\
\midrule
Type~2 10-reqs. & 17/20 & 18/20 \\
Type~2 12-reqs. & 14/20 & 15/20 \\
\bottomrule
\end{tabular}
}
\end{minipage}
\hfill
\begin{minipage}{0.45\textwidth}
\centering
\scriptsize
\vspace{0.3em}
\okC{
\text{b) MSLP}\\[0.3em]
\begin{tabular}{lrr}
\toprule
 & Multi-island & Multi-floor \\
\midrule
Type~1 40-reqs. & 180/200 & 190/200 \\
Type~1 60-reqs. & 168/200 & 153/200 \\
\bottomrule
\end{tabular}
}
\end{minipage}

\end{table}

\begin{table}[!htp]\centering
\caption{\okC{Number of observations per machine variation box plot}}
\label{tab:num_obs_mach_var_box_plots}
\begin{minipage}{0.45\textwidth}
\centering
\scriptsize
\vspace{0.3em}
\okC{
\text{a) MIP+VI}\\[0.3em]
\begin{tabular}{lrr}
\toprule
 & Multi-island & Multi-floor \\
\midrule
Type~2 03-mach. & 36/40 & 37/40 \\
Type~2 04-mach. & 35/40 & 36/40 \\
\bottomrule
\end{tabular}
}
\end{minipage}
\hfill
\begin{minipage}{0.45\textwidth}
\centering
\scriptsize
\vspace{0.3em}
\okC{
\text{b) MSLP}\\[0.3em]
\begin{tabular}{lrr}
\toprule
 & Multi-island & Multi-floor \\
\midrule
Type~1 05-mach. & 150/200 & 160/200 \\
Type~1 06-mach. & 198/200 & 183/200 \\
\bottomrule
\end{tabular}
}
\end{minipage}

\end{table}

\begin{table}[!htp]\centering
\caption{\okC{Number of observations per region variation box plot}}
\label{tab:num_obs_reg_var_box_plots}
\scriptsize
\begin{minipage}{0.45\textwidth}
\centering
\vspace{0.3em}
\okC{
\text{a) MIP+VI}\\[0.3em]
\begin{tabular}{lrr}
\toprule
 & Multi-island & Multi-floor \\
\midrule
Type~2 02-regs. & 38/40 & 40/40 \\
Type~2 04-regs. & 33/40 & 33/40 \\
\bottomrule
\end{tabular}
}
\end{minipage}
\hfill
\begin{minipage}{0.45\textwidth}
\centering
\vspace{0.3em}
\okC{
\text{b) MSLP}\\[0.3em]
\begin{tabular}{lrr}
\toprule
 & Multi-island & Multi-floor \\
\midrule
Type~1 04-regs. & 548/600 & 543/600 \\
\bottomrule
\end{tabular}
}
\end{minipage}
\end{table}

\begin{figure}[!ht]
    \centering
    \begin{minipage}{.67\textwidth}
        \centering
        \subfigure[Multi-island instances (MIP+VI)\label{fig:mean_max_load_mip_1}]{
            \centering
            \includegraphics[scale=0.42]{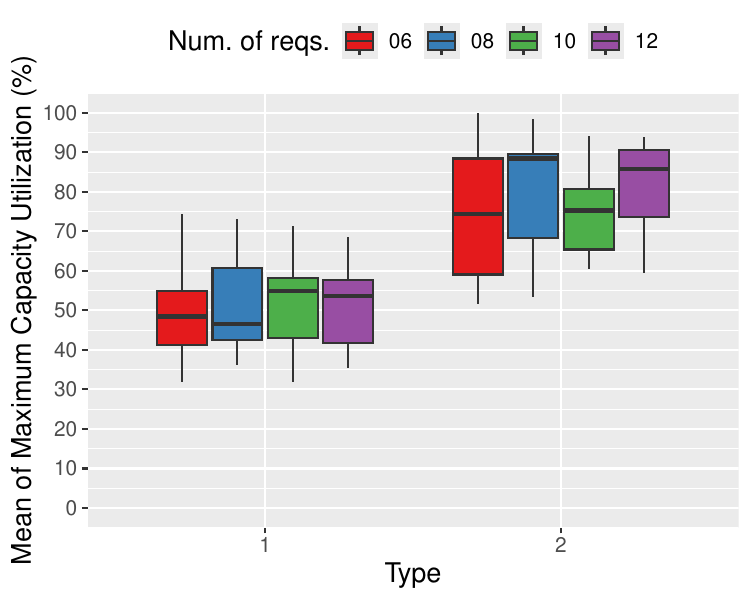}
        }
        \subfigure[Multi-floor instances (MIP+VI)\label{fig:mean_max_load_mip_2}]{
            \centering
            \includegraphics[scale=0.42]{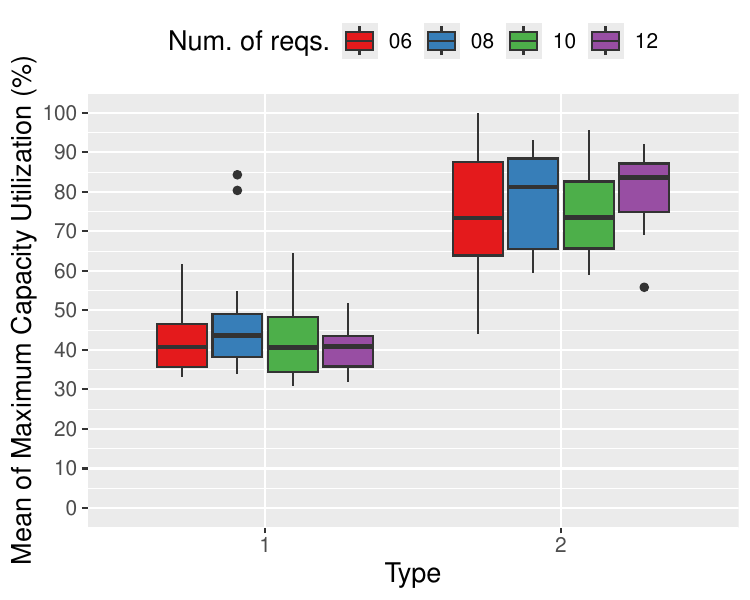}
        }
        \\
        \subfigure[Multi-island instances (MSLP)\label{fig:mean_max_load_mslp_1}]{
            \centering
            \includegraphics[scale=0.42]{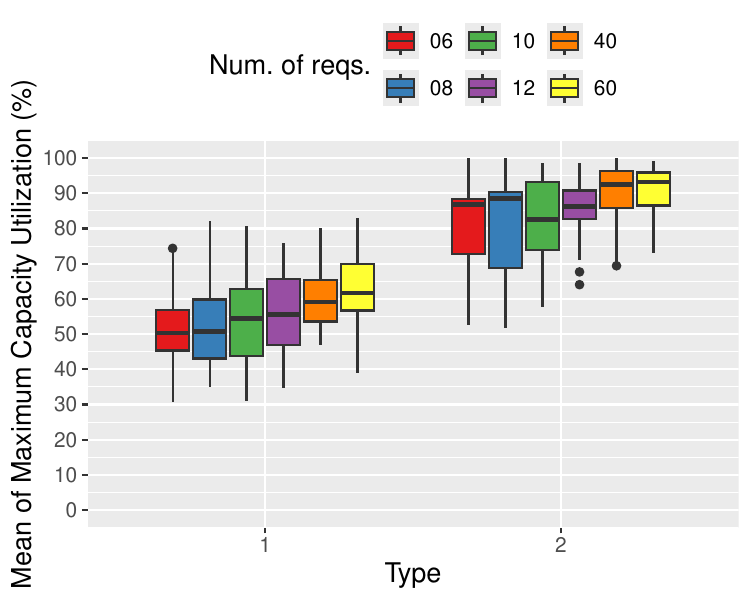}
        }
        \subfigure[Multi-floor instances (MSLP)\label{fig:mean_max_load_mslp_2}]{
            \centering
            \includegraphics[scale=0.42]{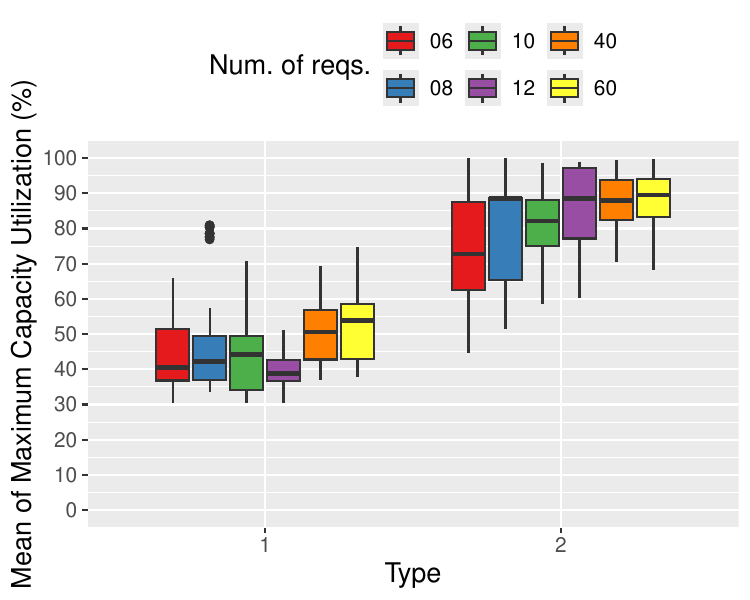}
        }
    \end{minipage}
    \caption{Mean of Maximum Capacity Utilization (\%) x Num. of Requests}
    \label{fig:mean_max_cap_utilization_vs_num_reqs}
\end{figure}

\begin{figure}[!ht]
    \centering
    \begin{minipage}{.67\textwidth}
        \subfigure[Multi-island instances (MIP+VI)\label{fig:n_vehicles_mip_1}]{
            \centering
            \includegraphics[scale=0.42]{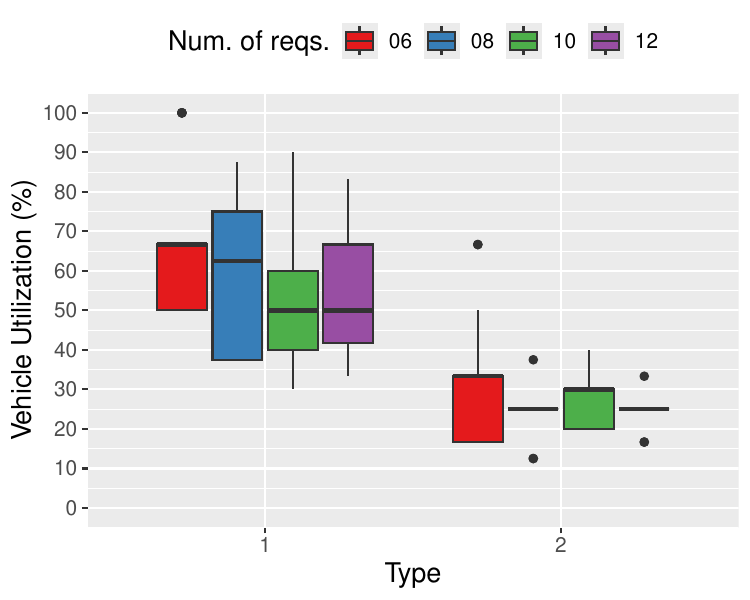}
        }
        \subfigure[Multi-floor instances (MIP+VI)\label{fig:n_vehicles_mip_2}]{
            \centering
            \includegraphics[scale=0.42]{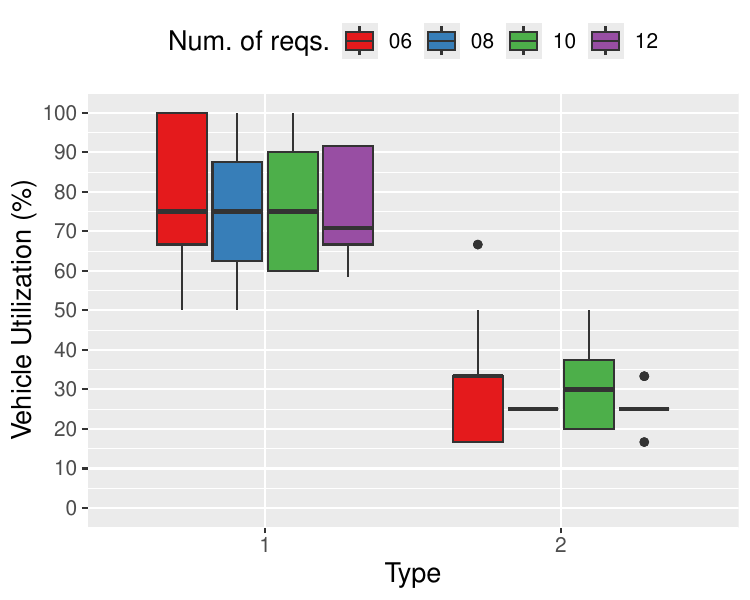}
        }
        \\
        \subfigure[Multi-island instances (MSLP)\label{fig:n_vehicles_mslp_1}]{
            \centering
            \includegraphics[scale=0.42]{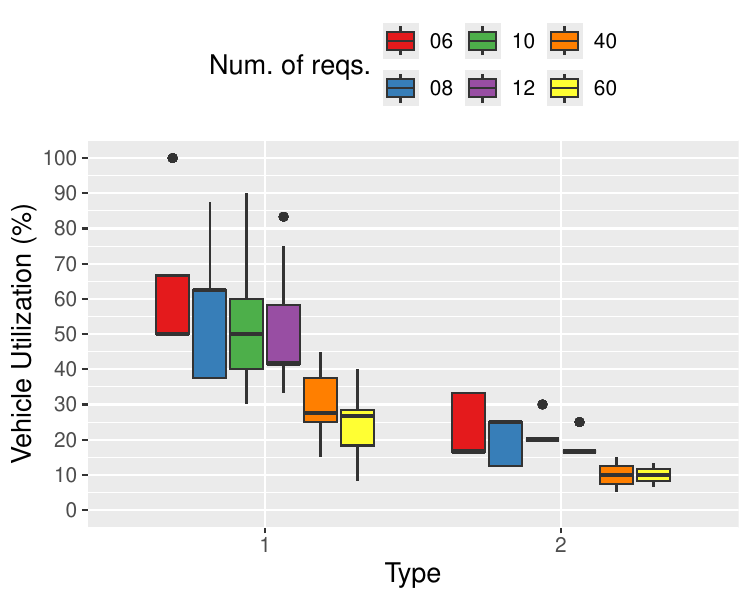}
        }
        \subfigure[Multi-floor instances (MSLP)\label{fig:n_vehicles_mslp_2}]{
            \centering
            \includegraphics[scale=0.42]{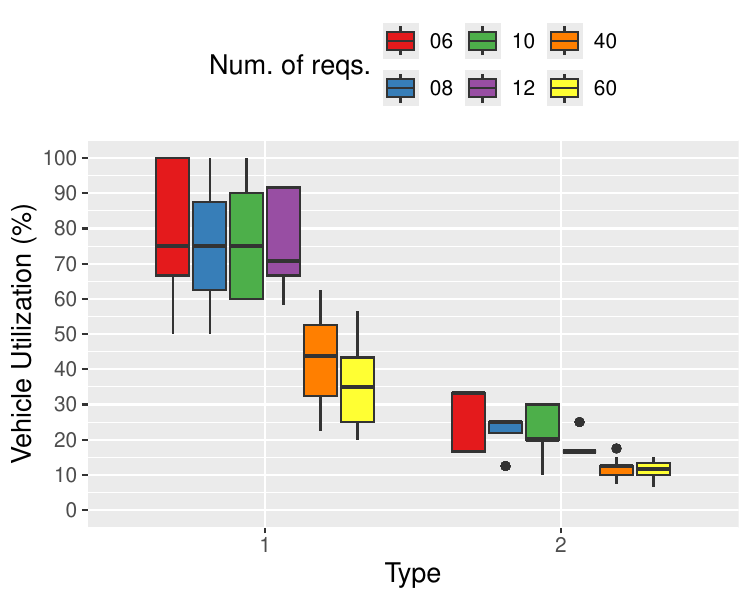}
        }
    \end{minipage}
    \caption{Vehicle Utilization (\%) x Num. of Requests}
    \label{fig:vehicle_utilization_vs_num_reqs}
\end{figure}

The box plots of Figure \ref{fig:mean_max_cap_utilization_vs_num_reqs} depict the distribution of the mean of maximum capacity utilization for each instance family, separated by type.
The maximum capacity utilization measures how much load of a vehicle's capacity is being carried at its peak, expressed as a percentage.
Also, this mean only considers vehicles that \okC{serve} at least one request. It is noticeable that in Type 1 instances, the vehicles are less used up to their maximum capacity throughout the routes.
The box plots of Figure~\ref{fig:vehicle_utilization_vs_num_reqs} explain this scenario by presenting the vehicle utilization distribution for each instance family separated by type.
The vehicle utilization defines how many vehicles were used to \okC{serve} at least one request, expressed as a percentage of the total number of vehicles.
\okC{Notably, the solutions for Type~1 instances made more use of vehicles than the solutions for Type~2 instances.
However, the mean of maximum capacity utilization is lower on Type~1 instance solutions, i.e., the capacity of the vehicles was not well utilized.
Therefore, the requests in Type~1 instance solutions are more evenly distributed among the vehicles due to the instances' short-scheduling horizons.
Additionally, the percentage of vehicles used does not increase proportionally as the number of requests increase.}
\FloatBarrier

\FloatBarrier

\begin{figure}[!htp]
    \centering
    \begin{minipage}{.67\textwidth}
        \subfigure[Multi-island instances (MIP+VI)\label{fig:mean_mach_act_time_mip_1}]{
            \centering
            \includegraphics[scale=0.42]{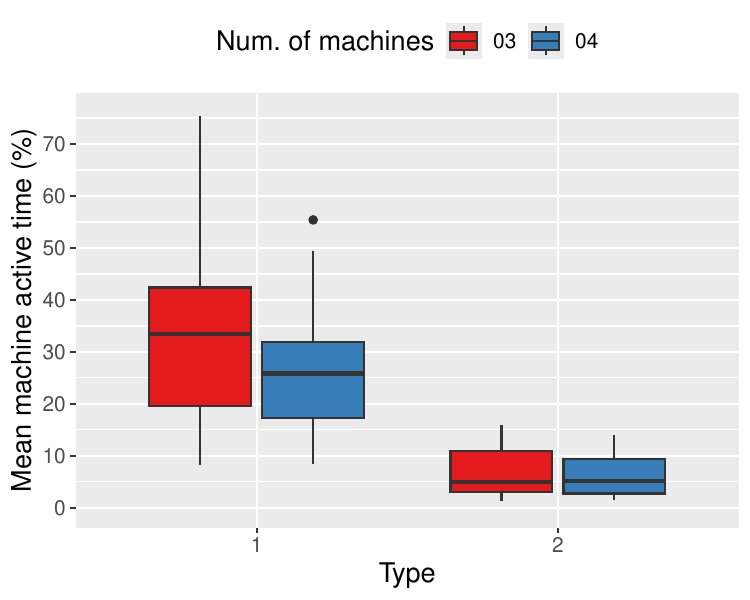}
        }
        \subfigure[Multi-floor instances (MIP+VI)\label{fig:mean_mach_act_time_mip_2}]{
            \centering
            \includegraphics[scale=0.42]{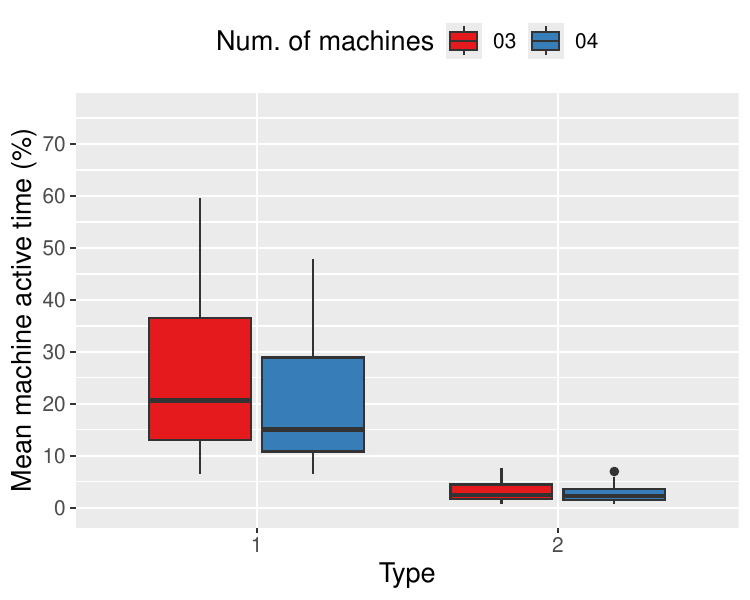}
        }
        \\
        \subfigure[Multi-island instances (MSLP) \label{fig:mean_mach_act_time_mslp_1}]{
            \centering
            \includegraphics[scale=0.42]{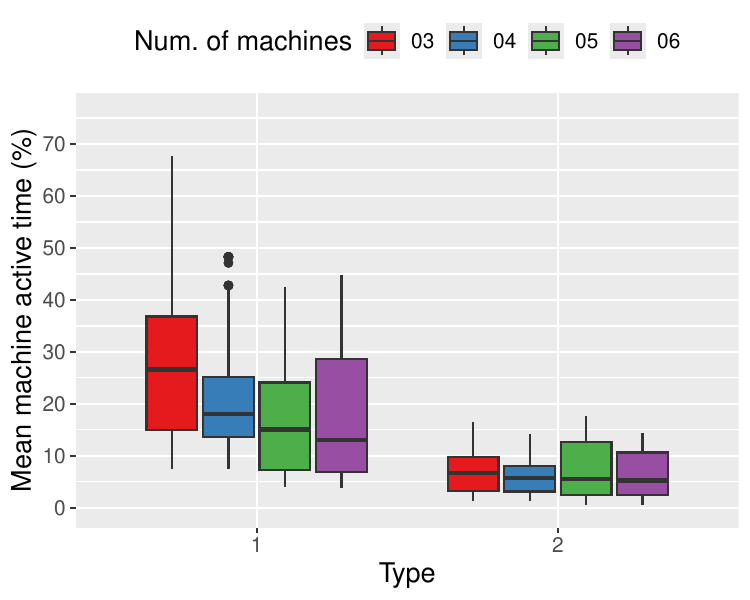}
        }
        \subfigure[Multi-floor instances (MSLP)\label{fig:mean_mach_act_time_mslp_2}]{
            \centering
            \includegraphics[scale=0.42]{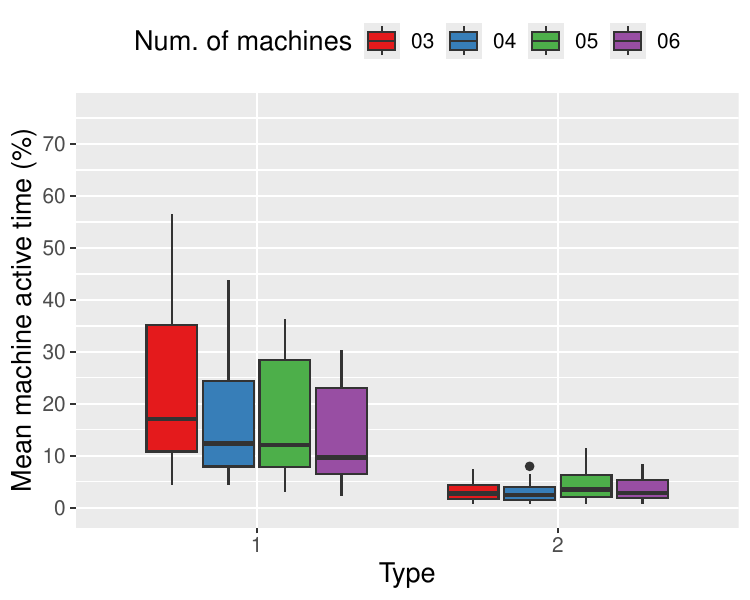}
        }
    \end{minipage}
    \caption{Mean Machine Active Time (\%) x Num. of Machines}
    \label{fig:mean_mach_act_time_vs_num_machs}
\end{figure}

\begin{figure}[!htp]
    \centering
        \begin{minipage}{.67\textwidth}
            \subfigure[Multi-island instances (MIP+VI)\label{fig:n_machines_mip_1}]{
                \centering
                \includegraphics[scale=0.42]{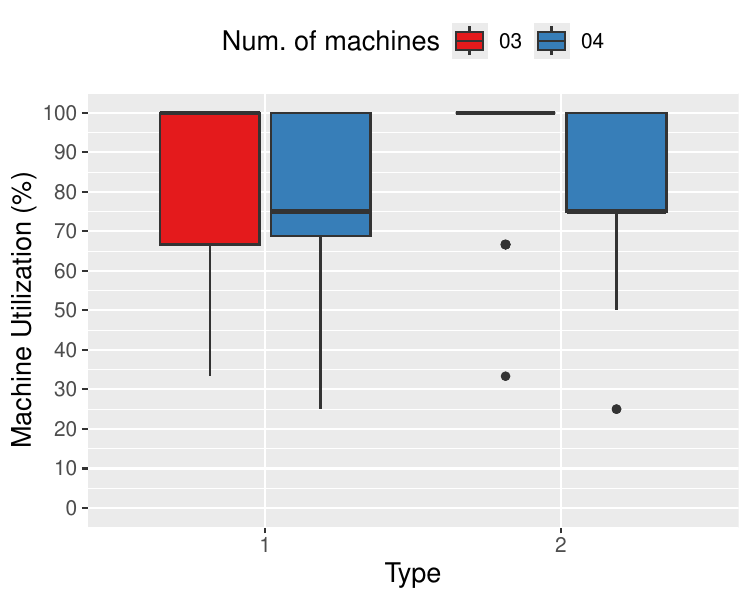}
            }
            \subfigure[Multi-floor instances (MIP+VI)\label{fig:n_machines_mip_2}]{
                \centering
                \includegraphics[scale=0.42]{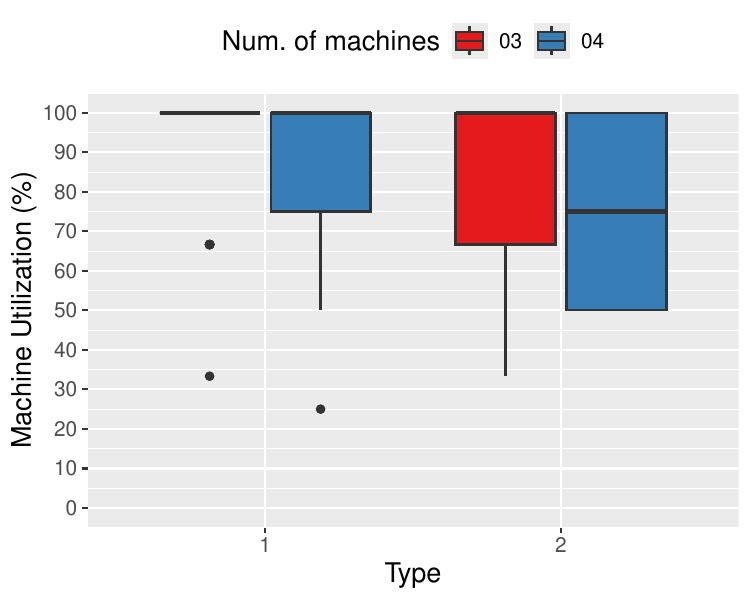}
            }
            \\
            \subfigure[Multi-island instances (MSLP)\label{fig:n_machines_mslp_1}]{
                \centering
                \includegraphics[scale=0.42]{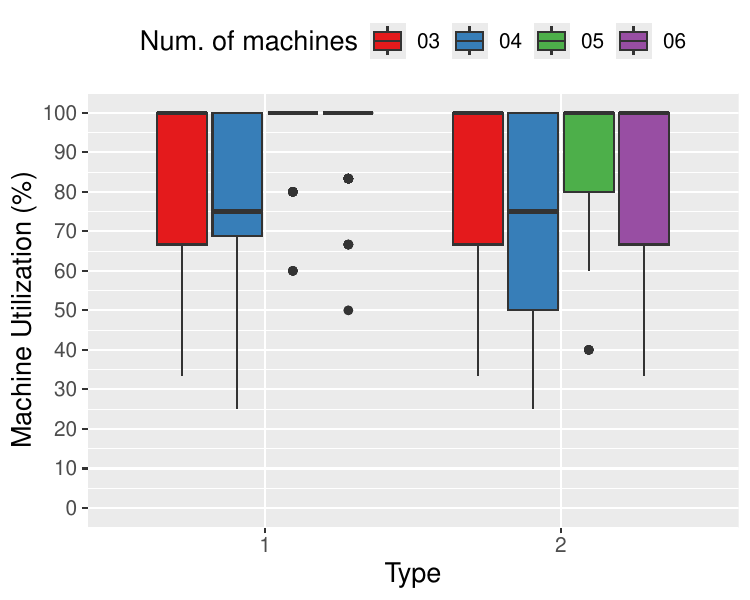}
            }
            \subfigure[Multi-floor instances (MSLP)\label{fig:n_machines_mslp_2}]{
                \centering
                \includegraphics[scale=0.42]{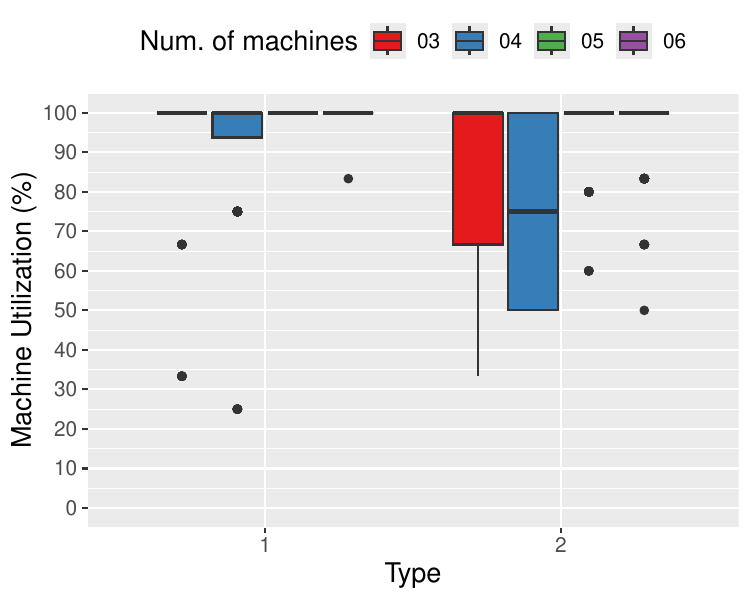}
            }
        \end{minipage}
    \caption{Machine Utilization (\%) x Num. of Machines}
    \label{fig:mach_utilization_vs_num_machs}
\end{figure}

The box plots of Figure \ref{fig:mean_mach_act_time_vs_num_machs} depict the mean machine active time distribution for all instances, varying the number of machines.
A machine is said to be active whenever it moves between stations, carrying or not a vehicle.
Also, this active time is expressed as a percentage of the total planning time horizon and the mean is over all machines used.
The box plots of Figure \ref{fig:mach_utilization_vs_num_machs} depict the machine utilization distribution for all instances, varying the number of machines.
\okC{Note that the box plots related to three- or four-machine instances comprise the smaller instances with up to 12 requests, while those related to five- or six-machine instances encompass the larger instances with 40 and 60 requests.}
It can be seen from the box plots of Figures \ref{fig:mean_mach_act_time_vs_num_machs} and \ref{fig:mach_utilization_vs_num_machs} that the median mean machine active times decrease when the number of machines available increases, even though the machine utilization decreases.
\okC{
Besides, it is notable that almost all solutions for the larger instances used all the machines available (except some Type~2 multi-island instances), but their mean active time was not greater.} \okC{Finally, the mean machine active times are generally greater on Type~1 instances when compared to Type~2 instances.}
\FloatBarrier

\FloatBarrier
% \newpage

\begin{figure}[!ht]
    \centering
    \begin{minipage}{.67\textwidth}
        \subfigure[Multi-island instances (MIP+VI)]{
            \centering
            \includegraphics[scale=0.42]{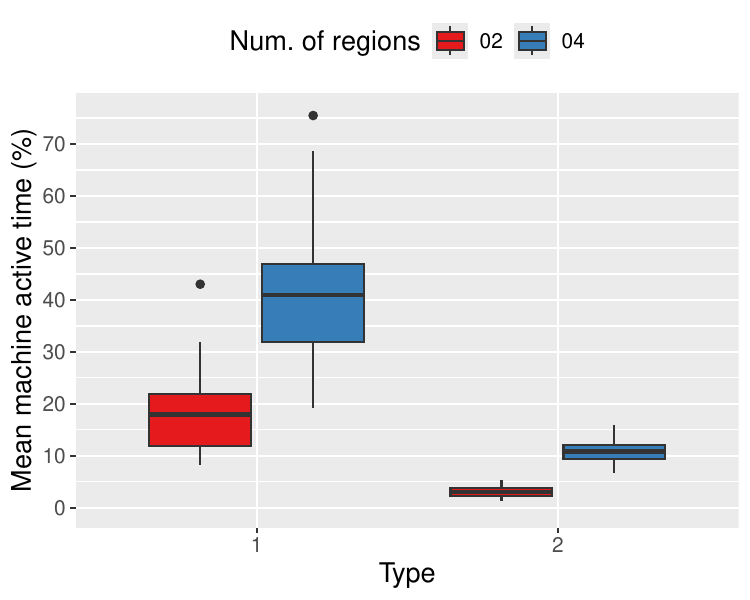}
        }
        \subfigure[Multi-floor instances (MIP+VI)]{
            \centering
            \includegraphics[scale=0.42]{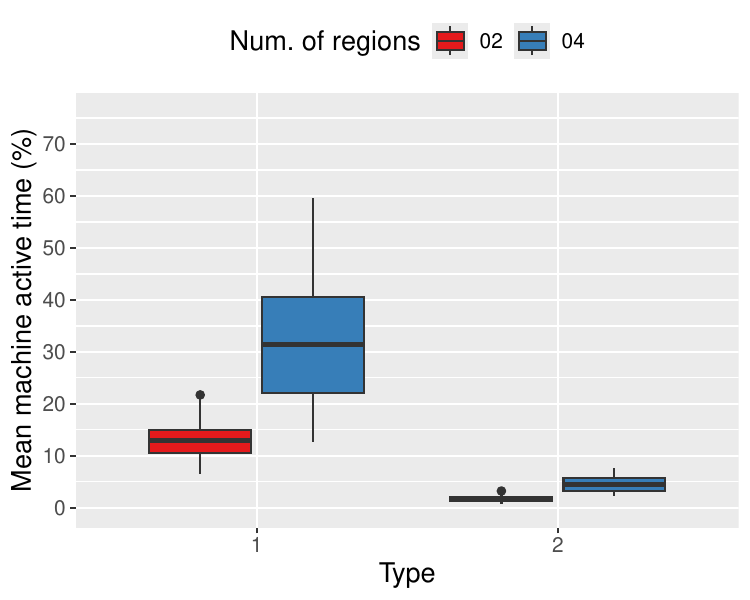}
        }
        \\
        \subfigure[Multi-island instances (MSLP)]{
            \centering
            \includegraphics[scale=0.42]{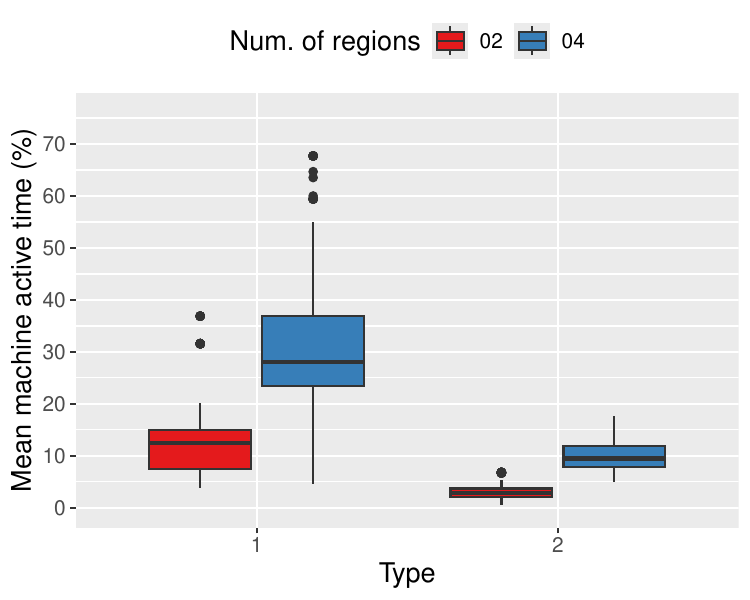}
        }
        \subfigure[Multi-floor instances (MSLP)]{
            \centering
            \includegraphics[scale=0.42]{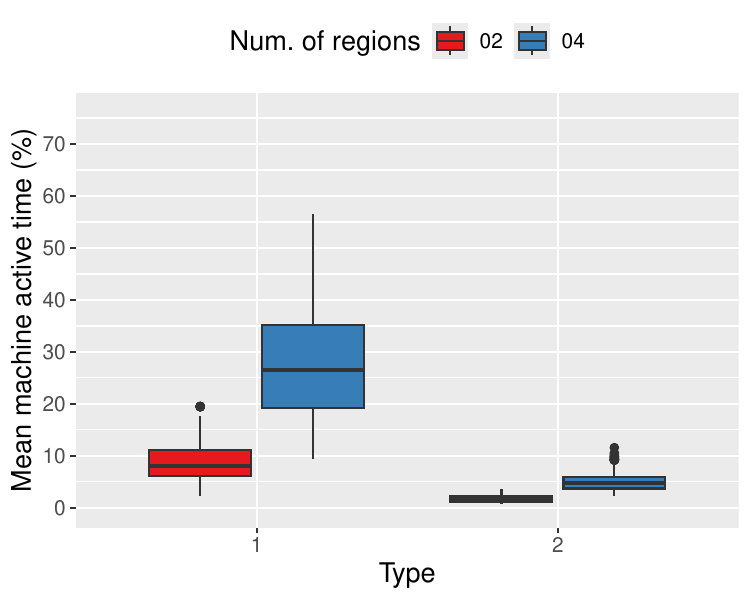}
        }
    \end{minipage}
    \caption{Mean Machine Active Time (\%) x Num. of Regions}
    \label{fig:mean_mach_act_time_vs_num_regions}
\end{figure}

The box plots of Figure \ref{fig:mean_mach_act_time_vs_num_regions}
depict the mean machine active time distribution for all instances, varying the number of regions.
It can be seen that increasing the number of regions also significantly increases the mean machine active times for both multi-island and multi-floor instances. 
Indeed, the median value in each type increased more than twice \okC{as the number of regions increased}.

\FloatBarrier
%\newpage

\end{document}